\documentclass[11pt,openbib]{article} 
\usepackage{amsfonts,amsthm}
\usepackage{amssymb,amsmath}
\usepackage{epic}
\newcommand{\R}{\mathbb{R}} 
\newcommand{\C}{\mathbb{C}}
 
\newcommand{\Z}{\mathbb{Z}}
 
\newcommand{\gl}{\mathop{\mathrm{gl}}\nolimits}

\newcommand{\ad}{\mathop{\mathrm{ad}}\nolimits}

\newcommand{\zetaminus}{\mbox{$\, {-}_{\zeta }$}}

\newcommand{\smallzetaminus}{\mbox{{\tiny ${-}_{\zeta }$}}}

\newcommand{\vvee}{\mbox{\tiny $\vee $}}

\newcommand{\setrule}{\, \mathop{\rule[-4pt]{.5pt}{13pt}\, }\nolimits}
\newcommand{\smallspace}{\smallskip\par\noindent}

\newcommand{\rowspace}{\rule{0pt}{16pt}}

\newcommand{\onehalf}{\mbox{$\frac{\scriptstyle 1}{\scriptstyle 2}\,$}} 
\newcommand{\spann}{\mathop{\rm span}\nolimits} 

\allowdisplaybreaks
\begin{document}
\begin{center}
{\bf \Large Very special sandwich algebras}
\vspace{.05in}\par\noindent \hspace{.75in} Richard Cushman\footnotemark 
\end{center} 

\footnotetext{\parbox[t]{4in}{version: \today \\
email: rcushman@ucalgary.ca \\ 
Department of Mathematics and Statistics \\ 
University of Calgary} } 

\noindent The goal of this paper is to investigate a class of algebras called sandwich algebras, which are certain complex Lie algebras with a nilpotent radical whose elements are sandwiches.\footnote{See Cohen, et al. \cite{cohen} 
who originated this terminology.}\! We present a classification of all very special sandwich algebras.  

\section{Definition of a sandwich algebra}

We say that a complex Lie algebra $\widetilde{\mathfrak{g}}$ is a \emph{sandwich algebra} if it has 
the following three properties. \medskip 

\indent 1. \parbox[t]{4.5in}{$\widetilde{\mathfrak{g}}$ has a nilpotent radical $\widetilde{\mathfrak{n}}$, which is 
a sandwich, that is, $[\widetilde{\mathfrak{n}}, [\widetilde{\mathfrak{n}}, \widetilde{\mathfrak{n}}]] =0$.} 
\smallspace
\indent 2. \parbox[t]{4.5in}{The complex Lie algebra $\widetilde{\mathfrak{g}}/\widetilde{\mathfrak{n}}$ is isomorphic 
to a semisimple complex Lie algebra $\mathfrak{g}$.} 
\smallspace
\indent 3. \parbox[t]{4.5in}{Let $\mathfrak{h}$ be a Cartan subalgebra of $\mathfrak{g}$. Then 
${\ad }_{\mathfrak{h}} = {\{ {\ad }_H \} }_{H \in \mathfrak{h}}$ is a maximal family of 
commuting semisimple linear maps of $\widetilde{\mathfrak{n}}$ into itself.}\medskip 

\noindent If $\mathfrak{g}$ is a semisimple Lie algebra, then $\mathfrak{g}$ is a sandwich algebra 
with nilpotent radical $\{ 0 \} $. From now on we assume that $\widetilde{\mathfrak{n}} \ne \{ 0 \} $. \medskip 

\noindent A sandwich algebra $\widetilde{\mathfrak{g}}= \mathfrak{g} \oplus \widetilde{\mathfrak{n}}$ 
is \emph{simple} if $\mathfrak{g}$ is a simple Lie algebra. \medskip 

\noindent Let $\widehat{\mathcal{R}} \subseteq {\mathfrak{h}}^{\ast }$ be a finite set of roots 
corresponding to the maximal toral subalgebra ${\ad }_{\mathfrak{h}}$ of $\gl \, (\widetilde{\mathfrak{n}}, \C)$. The root space ${\widehat{\mathfrak{g}}}_{\widehat{\alpha }}$ corresponding to the root $\widehat{\alpha } \in 
\widehat{\mathcal{R}}$ is spanned by a nonzero root vector $X_{\widehat{\alpha }} \in \widetilde{\mathfrak{n}}$. 
Then $\widetilde{\mathfrak{n}} = \sum_{\widehat{\alpha } \in \widehat{\mathcal{R}}} \oplus 
{\widehat{\mathfrak{g}}}_{\widehat{\alpha }}$ and $\widetilde{\mathfrak{g}} 
= \mathfrak{g} \oplus \widetilde{\mathfrak{n}}$. \medskip 

\noindent \textbf{Claim 1.1} Every sandwich algebra is a direct sum of simple sandwich \linebreak 
algebras. \medskip 

\noindent \textbf{Proof.} Suppose that $\widetilde{\mathfrak{g}} = \mathfrak{g} \oplus \widetilde{\mathfrak{n}}$ 
is a sandwich algebra. Since the Lie algebra $\mathfrak{g}$ is semisimple with Cartan subalgebra 
$\mathfrak{h}$ and associated root system $\Phi $, we may write $\mathfrak{g} = 
\sum^n_{i=1} \oplus {\mathfrak{g}}_i$, where ${\mathfrak{g}}_i$ is a simple Lie algebra with Cartan 
subalgebra ${\mathfrak{h}}_i$ and associated root system ${\Phi }_i$. Then $\mathfrak{h} = 
\sum^n_{i=1}\oplus {\mathfrak{h}}_i$, where ${\mathfrak{h}}_i$ are nondegenerate orthogonal 
subspaces of $\mathfrak{h}$ under the Killing form on $\mathfrak{g}$. Thus 
the we may write the collection $\widehat{\mathcal{R}}$ of roots associated to the nilpotent radical 
$\widetilde{\mathfrak{n}}$ of the sandwich algebra $\widetilde{\mathfrak{g}}$ as 
$\coprod_{1\le i \le n}{\widehat{\mathcal{R}}}_i$, where ${\widehat{\mathcal{R}}}_i = 
\widehat{\mathcal{R}} \cap {\mathfrak{h}}^{\ast }_i = 
\{ \widehat{\alpha }|{\mathfrak{h}}_i \in {\mathfrak{h}}^{\ast }_i \, \setrule \, \widehat{\alpha } 
\in \widehat{\mathcal{R}} \} $. 
Let ${\widetilde{\mathfrak{n}}}_i = \sum_{{\widehat{\alpha }}_i \in {\widehat{\mathcal{R}}}_i} \oplus 
{\widehat{\mathfrak{g}}}_{{\widehat{\alpha }}_i}$ and let ${\widetilde{\mathfrak{g}}}_i = {\mathfrak{g}}_i \oplus 
{\widetilde{\mathfrak{n}}}_i$. We get  
\begin{align}
\widetilde{\mathfrak{g}} & = \mathfrak{g} \oplus \widetilde{\mathfrak{n}} = 
\mathfrak{h} \oplus \sum_{\alpha \in \Phi } {\mathfrak{g}}_{\alpha } \oplus 
\sum_{\widehat{\alpha } \in \widehat{\mathcal{R}}}{\mathfrak{g}} {\widehat{\mathfrak{g}}}_{\widehat{\alpha}} 
\notag \\
& = \sum^n_{i=1}{\mathfrak{h}}_i \oplus \sum^n_{i=1}  \sum_{{\alpha }_i \in {\Phi }_i } {\mathfrak{g}}_{{\alpha }_i} 
\oplus \sum^n_{i=1} \sum_{{\widehat{\alpha }}_i \in {\widehat{\mathcal{R}}}_i} \oplus  
{\widehat{\mathfrak{g}}}_{{\widehat{\alpha}}_i}, 
\notag \\
&\hspace{.5in} \parbox[t]{3in}{since $\mathfrak{g} = 
\sum^n_{i=1}\oplus {\mathfrak{g}}_{{\alpha }_i}$ and ${\mathfrak{g}}_{\widehat{\alpha }} = 
\sum^n_{i=1}\oplus {\mathfrak{g}}_{{\widehat{\alpha }}_i}$ because 
$\alpha = \sum^n_{i=1}{\alpha }_i$ and $\widehat{\alpha } = \sum^n_{i=1}{\widehat{\alpha }}_i $} 
\notag \\
& = \sum^n_{i=1} \oplus \big( {\mathfrak{h}}_i \oplus \sum_{{\alpha }_i \in {\Phi }_i}{\mathfrak{g}}_{{\alpha }_i} 
\oplus \sum_{{\widehat{\alpha }}_i \in {\widehat{\mathcal{R}}}_i} {\mathfrak{g}}_{{\widehat{\alpha }}_i} \big) 
\notag \\
& = \sum^n_{i=1} \big( {\mathfrak{g}}_i \oplus {\widetilde{\mathfrak{n}}}_i \big) = 
\sum^n_{i=1}\oplus {\widetilde{\mathfrak{g}}}_i. 
\label{eq-sec1one}
\end{align}

We now verify that ${\widetilde{\mathfrak{g}}}_i$ is a sandwich algebra by showing that 
properties 1--3 hold. \medskip 

\noindent 1. The following argument shows that ${\widetilde{\mathfrak{n}}}_i$ is the nilpotent radical 
of ${\widetilde{\mathfrak{g}}}_i$. First we show that ${\widetilde{\mathfrak{n}}}_i$ is an ideal of 
${\widetilde{\mathfrak{g}}}_i$. Clearly ${\widetilde{\mathfrak{n}}}_i$ is a vector space. We have 
\begin{align}
[ {\mathfrak{g}}_i, {\widetilde{\mathfrak{n}}}_i] & = 
[ H_i + \sum_{\alpha \in {\Phi }_i} c_{\alpha } X_{\alpha }, 
\sum_{\widehat{\alpha } \in {\widehat{\mathcal{R}}}_i} c_{\widehat{\alpha }} X_{\widehat{\alpha }} ], 
\, \, \mbox{where the bracket is that of $\widetilde{\mathfrak{g}}$} 
\notag \\
& = \sum_{\widehat{\alpha } \in {\widehat{\mathcal{R}}}_i} c_{\widehat{\alpha }}\, \widehat{\alpha }(H_i)  X_{\widehat{\alpha }}  + \sum_{\alpha \in {\Phi }_i, \, \widehat{\alpha } \in {\widehat{\mathcal{R}}}_i} c_{\alpha }c_{\widehat{\alpha }} \, [X_{\alpha }, X_{\widehat{\alpha }} ] \notag \\
& = \sum_{\widehat{\alpha } \in {\widehat{\mathcal{R}}}_i} c_{\widehat{\alpha }}\widehat{\alpha }(H_i)  X_{\widehat{\alpha }}  + \sum_{\alpha + \widehat{\alpha } \in {\Phi }_i } c_{\alpha }c_{\widehat{\alpha }} \, c_{\alpha , \widehat{\alpha }} \,  X_{\alpha + \widehat{\alpha }} . 
\label{eq-sec1two} 
\end{align}
So the second term in (\ref{eq-sec1two}) lies in ${\mathfrak{g}}_i$, whereas the first term in (\ref{eq-sec1two}) lies 
in ${\widetilde{\mathfrak{n}}}_i$. Hence $[{\mathfrak{g}}_i, {\widetilde{\mathfrak{n}}}_i] \subseteq 
{\widetilde{\mathfrak{g}}}_i$. Also 
\begin{align}
[{\widetilde{\mathfrak{n}}}_i, {\widetilde{\mathfrak{n}}}_i] & = 
[\sum_{\widehat{\alpha } \in {\widehat{\mathcal{R}}}_i} b_{\widehat{\alpha}} X_{\widehat{\alpha }} , 
\sum_{\widehat{\beta } \in {\widehat{\mathcal{R}}}_i} c_{\widehat{\beta}} X_{\widehat{\beta }} ], 
\, \, \mbox{where the bracket is that of $\widetilde{\mathfrak{g}}$} \notag \\
& = \sum_{\widehat{\alpha }, \widehat{\beta } \in {\widehat{\mathcal{R}}}_i} b_{\widehat{\alpha }} c_{\widehat{\beta }} 
\, [X_{\widehat{\alpha }}, X_{\widehat{\beta }}] = 
\sum_{\widehat{\alpha } + \widehat{\beta } \in {\widehat{\mathcal{R}}}_i} b_{\widehat{\alpha }} c_{\widehat{\beta }} 
\, X_{\widehat{\alpha } +\widehat{\beta }} , \notag 
\end{align}
which lies in ${\widetilde{\mathfrak{n}}}_i$. So $[{\widetilde{\mathfrak{n}}}_i, {\widetilde{\mathfrak{n}}}_i ]
\subseteq {\widetilde{\mathfrak{n}}}_i$. Thus ${\widetilde{\mathfrak{n}}}_i$ is an ideal of 
${\widetilde{\mathfrak{g}}}_i$. Next we show that the ideal ${\widetilde{\mathfrak{n}}}_i$ is nilpotent. 
Let $n_i \in {\widetilde{\mathfrak{n}}}_i \subseteq \widetilde{\mathfrak{n}} $. Then ${\ad }_{n_i}$ is a nilpotent linear map of $\widetilde{\mathfrak{n}}$ into itself because $\widetilde{\mathfrak{n}}$ is a nilpotent ideal in 
$\widetilde{\mathfrak{g}}$. Hence ${\ad }_{n_i}$ is a nilpotent linear map of ${\widetilde{\mathfrak{n}}}_i$ 
into itself. So ${\widetilde{\mathfrak{n}}}_i$ is a nilpotent ideal. It is a sandwich because 
$\widetilde{\mathfrak{n}}$ is a sandwich. Finally, we show that the nilpotent ideal ${\widetilde{\mathfrak{n}}}_i$ 
is the maximal nilpotent ideal of ${\widetilde{\mathfrak{g}}}_i$. Suppose that there is a nilpotent ideal 
${\mathfrak{m}}_i$ of ${\widetilde{\mathfrak{g}}}_i$, which strictly contains ${\widetilde{\mathfrak{n}}}_i$. 
Then ${\mathfrak{m}}_i \cap {\mathfrak{g}}_i$  is a nonzero nilpotent ideal of ${\mathfrak{g}}_i$. But 
${\mathfrak{g}}_i$ is a simple Lie algebra. So either ${\mathfrak{m}}_i \cap {\mathfrak{g}}_i = \{ 0 \}$ or 
${\mathfrak{m}}_i \cap {\mathfrak{g}}_i = {\mathfrak{g}}_i$. Suppose that the second alternative holds. 
Then ${\mathfrak{g}}_i \subseteq {\mathfrak{m}}_i$. But ${\mathfrak{g}}_i$ contains a Cartan subalgebra 
${\mathfrak{h}}_i$, which consists of nonzero semisimple elements. This contradicts the fact 
that ${\mathfrak{m}}_i$ is consists of nilpotent elements. 
Thus ${\mathfrak{m}}_i \cap {\mathfrak{g}}_i = \{ 0 \}$, which contradicts the hypothesis that 
${\mathfrak{m}}_i \supsetneqq {\widetilde{\mathfrak{n}}}_i$. Hence the ideal ${\widetilde{\mathfrak{m}}}_i$ does not exist. So ${\widetilde{\mathfrak{n}}}_i$ is the maximal nilpotent ideal of ${\widetilde{\mathfrak{g}}}_i$. \medskip 

\noindent 2. Since ${\widetilde{\mathfrak{g}}}_i = {\mathfrak{g}}_i \oplus {\widetilde{\mathfrak{n}}}_i$, we 
get ${\widetilde{\mathfrak{g}}}_i /{\widetilde{\mathfrak{n}}}_i$ is isomorphic to ${\mathfrak{g}}_i$, which is 
a simple Lie algebra. \medskip 

\noindent 3. Since ${\mathfrak{h}}_i \subseteq \mathfrak{h}$ for every ${\widehat{\alpha }}_i \in 
{\mathcal{R}}_i$ there is an $\widehat{\alpha } \in \mathcal{R}$ such that ${\widehat{\alpha }}_i = 
\widehat{\alpha}|{\mathfrak{h}}_i$. Hence ${\mathfrak{g}}_{{\widehat{\alpha }}_i } \subseteq 
{\mathfrak{g}}_{\widehat{\alpha}}$. But $\dim {\mathfrak{g}}_{\widehat{\alpha}} =1$. So 
$\dim {\mathfrak{g}}_{{\widehat{\alpha }}_i } =1$. Consequently, ${\ad }_{{\mathfrak{h}}_i}$ is 
a maximal toral subalgebra of $\gl ({\widetilde{\mathfrak{n}}}_i, \C )$. \medskip 

\noindent Thus ${\widetilde{\mathfrak{g}}}_i$ is a sandwich algebra, which is simple because 
${\mathfrak{g}}_i$ is simple. From (\ref{eq-sec1one}) we see that we have proved the claim. \hfill $\square $ 

\section{Basic properties}

\noindent In this section we construct Heisenberg subalgbras of a given 
sandwich \linebreak 
algebra $\widetilde{\mathfrak{g}} = \mathfrak{g} \oplus \widetilde{\mathfrak{n}}$. \medskip 

\noindent Assume that $\widetilde{\mathfrak{n}}$ is nonabelian. Let $Z = [\widetilde{\mathfrak{n}}, \widetilde{\mathfrak{n}}]$. Since 
$\widetilde{\mathfrak{n}}$ is not abelian, $Z \ne \{ 0 \}$. Because $\widetilde{\mathfrak{n}}$ is a 
sandwich, $Z$ is the center of $\widetilde{\mathfrak{n}}$. Let $\mathcal{Z} = \{ \widehat{\alpha } 
\in \widehat{\mathcal{R}} \, \setrule \, X_{\widehat{\alpha }} \in Z \} $. For each 
$\zeta \in \mathcal{Z}$ let 
${\widetilde{\mathcal{R}}}_{\zeta } = \{ (\widehat{\alpha }, \widehat{\beta }) \in \widehat{\mathcal{R}} 
\times \widehat{\mathcal{R}} \, \setrule \, \widehat{\alpha } + \widehat{\beta } = \zeta \} $. \medskip 

\noindent \textbf{Lemma 2.1}. For each $\zeta \in \mathcal{Z}$, ${\widetilde{\mathcal{R}}}_{\zeta } \ne 
\varnothing$.  \medskip 

\noindent \textbf{Proof}. By definition, for each 
$\zeta \in \mathcal{Z}$ the nonzero root vector $X_{\zeta }$ lies in $Z$. Since $Z = 
[\widetilde{\mathfrak{n}}, \widetilde{\mathfrak{n}}]$,  there are nonzero 
$x = \sum_{\widehat{\alpha } \in \widehat{\mathcal{R}}} a_{\widehat{\alpha }} X_{\widehat{\alpha }}$ and 
$y = \sum_{\widehat{\beta } \in \widehat{\mathcal{R}}} b_{\widehat{\beta }} X_{\widehat{\beta }} \in 
\widetilde{\mathfrak{n}}$ such that 
\begin{align}
X_{\zeta } & = [x,y] = \sum_{\widehat{\alpha }, \widehat{\beta } \in \widehat{\mathcal{R}}} a_{\widehat{\alpha}}
b_{\widehat{\beta }} \, [X_{\widehat{\alpha }}, X_{\widehat{\beta }}] 
 = \sum_{\{ \widehat{\alpha }, \widehat{\beta } \in \widehat{\mathcal{R}} \, \setrule \, \widehat{\alpha } + \widehat{\beta } 
\in \widehat{\mathcal{R}} \}} \hspace{-24pt} a_{\widehat{\alpha}}
b_{\widehat{\beta }}\,  c_{\widehat{\alpha },\widehat{\beta }} X_{\widehat{\alpha } +\widehat{\beta }} , \notag 
\end{align} 
for some $c_{\widehat{\alpha},\widehat{\beta}} \in \C $, since 
$\left[ {\widehat{\mathfrak{g}}}_{\widehat{\alpha }}, {\widehat{\mathfrak{g}}}_{\widehat{\beta }} \right] = 
\left\{ \begin{array}{rl} 
{\widehat{\mathfrak{g}}}_{\widehat{\alpha }+\widehat{\beta }}, & 
\mbox{if $\widehat{\alpha }+\widehat{\beta } \in \widehat{\mathcal{R}}$} \\
0, & \mbox{otherwise} 
\end{array} \right. $. If $c_{\widehat{\alpha},\widehat{\beta}} =0 $ for all $\widehat{\alpha }$, 
$\widehat{\beta } \in \widehat{\mathcal{R}}$, then $X_{\zeta } =0$, which is a contradiction. 
Hence 
\begin{align} 
X_{\zeta } & = \sum_{\{ \widehat{\alpha } +\widehat{\beta } \ne \zeta \} }\hspace{-8pt} a_{\widehat{\alpha}}
b_{\widehat{\beta }} \, c_{\widehat{\alpha },\widehat{\beta }} X_{\widehat{\alpha } +\widehat{\beta }} +
\Big( \hspace{-12pt} \sum_{\{ \widehat{\alpha } +\widehat{\beta } = \zeta \}} \hspace{-8pt}a_{\widehat{\alpha}}b_{\widehat{\beta }}\,  c_{\widehat{\alpha },\widehat{\beta }} \Big) X_{\zeta } . \notag 
\end{align}
If $\{ \widehat{\alpha } +\widehat{\beta } = \zeta \} = \varnothing $, then the second term 
in the preceding equality is vacuous. Thus $X_{\zeta } $ is a nonzero vector in ${\spann }\{ X_{\widehat{\alpha } +\widehat{\beta }} \, \setrule \, \widehat{\alpha } +\widehat{\beta } \ne \zeta \} $, which contradicts the fact that the sum 
$\sum_{\widehat{\alpha } \in \widehat{\mathcal{R}}} \oplus 
{\widehat{\mathfrak{g}}}_{\widehat{\alpha }}$ is direct. Hence ${\widetilde{\mathcal{R}}}_{\zeta} \ne \varnothing $. \hfill {\footnotesize $\square $} \medskip  

\noindent From the definition of the root space 
${\widehat{\mathfrak{g}}}_{\widehat{\alpha }}$ and the fact that $\widetilde{\mathfrak{n}} = 
\sum_{\widehat{\alpha } \in \widehat{\mathcal{R}}}\oplus {\widehat{\mathfrak{g}}}_{\widehat{\alpha }}$, it follows that 
$Y = \sum_{\widehat{\alpha } \in \widehat{\mathcal{R}} \setminus \mathcal{Z}}\oplus 
{\widehat{\mathfrak{g}}}_{\widehat{\alpha }}$ is an ${\ad }_{\mathfrak{h}}$-invariant complementary subspace to $Z$ in $\widetilde{\mathfrak{n}}$. \medskip 

\noindent \textbf{Lemma 2.2} We have  
\begin{equation}
Z = [Y,Y].
\label{eq-zero}
\end{equation}
\textbf{Proof}. Because $Y \subseteq \widetilde{\mathfrak{n}}$ and 
$Z = [\widetilde{\mathfrak{n}}, \widetilde{\mathfrak{n}}]$ it follows that $[Y,Y] \subseteq Z$. Suppose that $[Y,Y] \ne Z$. Then there is a $z \in Z \setminus [Y,Y] $ such that 
$z = [n',n'' ]$ for some 
$n'$, $n'' \in \widetilde{\mathfrak{n}}$ both of which do not lie in $Y$. But 
$\widetilde{\mathfrak{n}} = Z \oplus Y$ by definition of $Y$. Hence there are 
$z'$, $z'' \in Z$, at least one of which is nonzero, and $y'$, $y'' \in Y$ such that $n' = z'+y'$ and $n'' = z''+y''$. So 
\begin{align}
z & = [n',n'' ] = [z' +y' , z'' +y'' ] \notag \\
& = [z',z'' ] +[z',y'' ] +[y',z'' ] + [y',y'' ] = [y',y'' ], \notag
\end{align}
since $Z$ is the center of $\widetilde{\mathfrak{n}}$. In other words, 
$z \in [Y,Y]$, which contradicts the definition of $z$. Therefore $Z = 
[Y,Y]$. \hfill {\footnotesize $\square $} \medskip 

\noindent For each $\zeta \in \mathcal{Z}$ we construct the set ${\mathcal{R}}^{\vvee}_{\zeta }$ as follows. 
We have shown that the set ${\widetilde{\mathcal{R}}}_{\zeta } \ne \varnothing $. Let 
${\widehat{\alpha }}_1 \in {\mathcal{R}}^{\vvee}_{\zeta }$ if $({\widehat{\alpha }}_1, {\widehat{\beta }}_1) 
\in {\widetilde{\mathcal{R}}}_{\zeta }$. Then $({\widehat{\beta }}_1, {\widehat{\alpha }}_1) 
\in {\widetilde{\mathcal{R}}}_{\zeta }$. Then recursively let ${\widehat{\alpha }}_{i+1} \in 
{\mathcal{R}}^{\vvee}_{\zeta }$ if $({\widehat{\alpha }}_{i+1}, {\widehat{\beta }}_{i+1}) \in 
{\widetilde{\mathcal{R}}}_{\zeta } \setminus \big\{ ({\widehat{\alpha }}_j, {\widehat{\beta }}_j ) \, \, \& \, \, 
({\widehat{\beta }}_j, {\widehat{\alpha }}_j ) \in {\widetilde{\mathcal{R}}}_{\zeta } \, \, 
\mbox{for $1 \le j \le i$} \big\} $. The recursion stops after a finite number of repetitions because 
${\widetilde{\mathcal{R}}}_{\zeta } \subseteq \widehat{\mathcal{R}} \times \widehat{\mathcal{R}}$, 
which is finite.  Given $\widehat{\alpha } \in {\mathcal{R}}^{\vvee}_{\zeta }$, there 
is a $\widehat{\beta } \in \widehat{\mathcal{R}}$ such that $\widehat{\alpha } + \widehat{\beta  } = \zeta $. This 
$\widehat{\beta }$ is uniquely determined and lies in ${\mathcal{R}}^{\vvee}_{\zeta }$ because the operation $+$ is commutative. Denote $\widehat{\beta }$ by $\zetaminus \widehat{\alpha }$. For each $\zeta \in \mathcal{Z}$ 
let ${\mathcal{R}}_{\zeta}  = \{ \widehat{\alpha } \in {\mathcal{R}}^{\vvee}_{\zeta }\, \setrule \, 
X_{\widehat{\alpha }} \in Y \, \, \, \mathrm{and} \, \, \, X_{\smallzetaminus \widehat{\alpha }} \in Y \}$. Note that if 
$\widehat{\alpha } \in {\mathcal{R}}_{\zeta}$ then so is $\zetaminus \widehat{\alpha }$. \medskip 

\noindent \textbf{Lemma 2.3} We have 
 \begin{equation}
\widehat{\mathcal{R}} \setminus \mathcal{Z} \subseteq 
\bigcup_{\zeta \in \mathcal{Z}}{\mathcal{R}}_{\zeta }. 
\label{eq-one}
\end{equation} 

\noindent \textbf{Proof}. Suppose that $\widehat{\alpha } \in \widehat{\mathcal{R}} \setminus \mathcal{Z}$. Then $X_{\widehat{\alpha }} \in Y$. Since $Z = [ \widetilde{\mathfrak{n}}, \widetilde{\mathfrak{n}} ]$, for every nonzero $n \in \widetilde{\mathfrak{n}}$ there is a ${\zeta}_n \in \mathcal{Z}$ such that  $X_{{\zeta}_n} = [X_{\widehat{\alpha }}, n ] $. Suppose that $X_{{\zeta }_n} =0$ for every $n \in \widetilde{\mathfrak{n}}$. Then 
$X_{\widehat{\alpha }} \in Z$. Consequently, 
$\widehat{\alpha } \in \mathcal{Z}$, which contradicts the definition of $\widehat{\alpha}$. So for some nonzero 
$n \in \widetilde{\mathfrak{n}}$, we have $X_{\zeta } = [X_{\widehat{\alpha }}, n]$ for some nonzero $X_{\zeta } \in Z$. Because $\mathfrak{n} = \sum_{\widehat{\gamma } \in \widehat{\mathcal{R}}} \oplus {\widehat{\mathfrak{g}}}_{\widehat{\gamma }}$, we may write $n =
\sum_{\widehat{\gamma } \in \widehat{\mathcal{R}}} a_{\widehat{\gamma }}
X_{\widehat{\gamma }}$. Since $X_{\zeta} = \sum_{\widehat{\gamma } \in \widehat{\mathcal{R}}} a_{\widehat{\gamma }}[X_{\widehat{\alpha }}, X_{\widehat{\gamma }}]$ and $X_{\zeta } \ne 0$, there is a $\widehat{\beta } \in \widehat{\mathcal{R}}$ such that 
$[X_{\widehat{\alpha }}, X_{\widehat{\beta }}] \ne 0$. In fact, $\widehat{\beta } \in \widehat{\mathcal{R}} \setminus 
\mathcal{Z}$; for if $\widehat{\beta} \in \mathcal{Z}$, then $X_{\widehat{\beta }} \in Z$, which implies $[X_{\widehat{\alpha }}, X_{\widehat{\beta }}] = 0$. This is a contradiction. Thus $X_{\widehat{\beta}} \in Y$. Since $[\widetilde{\mathfrak{n}}, \widetilde{\mathfrak{n}}] =Z$, there is a nonzero $X_{\widetilde{\zeta}} \in Z$ such that $X_{\widetilde{\zeta }} = [X_{\widehat{\alpha }}, X_{\widehat{\beta }}]$. In other words, $\widetilde{\zeta} = \widehat{\alpha } + 
\widehat{\beta }$. To see this we argue as follows. For every $H \in \mathfrak{h}$ we have 
\begin{align} 
\widetilde{\zeta} (H) X_{\widetilde{\zeta }} & = [H, X_{\widetilde{\zeta }} ] = 
[H, [X_{\widehat{\alpha }}, X_{\widehat{\beta }}]] \notag \\
& = [[H,X_{\widehat{\alpha }}], X_{\widehat{\beta }} ] + [ X_{\widehat{\alpha }}, [ H, X_{\widehat{\beta }} ]] 
= \big(  \widehat{\alpha } + \widehat{\beta } \big) (H) [X_{\widehat{\alpha }}, X_{\widehat{\beta }}] \notag \\
& = \big(  \widehat{\alpha } + \widehat{\beta } \big) (H) X_{\widetilde{\zeta }}. \notag 
\end{align}
The assertion follows because $X_{\widetilde{\zeta }}$ is nonzero. Therefore $\widehat{\alpha } \in 
{\mathcal{R}}_{\widetilde{\zeta }}$, because 1) $\widehat{\alpha } \in \widehat{\mathcal{R}} \setminus \mathcal{Z}$ implies $X_{\widehat{\alpha }} \in Y$; 2) $X_{\widehat{\beta }} \in Y$; 
and 3) $\widehat{\alpha } + \widehat{\beta } = \widetilde{\zeta }$. Consequently, (\ref{eq-one}) holds. 
\hfill {\footnotesize $\square $} \medskip 

\noindent  \textbf{Corollary 2.4} For each $\zeta \in \mathcal{Z}$ let $Y^{\zeta } = {\spann}_{\C } \{ X_{\widehat{\alpha }} \in Y \, \setrule \, \widehat{\alpha } \in {\mathcal{R}}_{\zeta} \} $. Then 
\begin{equation}
Y = \mbox{\Large $+$}\hspace{-13pt}\raisebox{-9pt}{$\scriptstyle \zeta \in \mathcal{Z}$} Y^{\zeta }. 
\label{eq-onestar}
\end{equation}  

\noindent \textbf{Proof}. Equation (\ref{eq-onestar}) follows immediately from 
(\ref{eq-one}) and the fact that $Y = {\spann }_{\C } \{ X_{\widehat{\alpha }} \in {\widehat{\mathfrak{g}}}_{\widehat{\alpha }} \, \setrule \, 
\widehat{\alpha } \in \widehat{\mathcal{R}} \setminus \mathcal{Z} \} $. \hfill {\footnotesize $\square $} \medskip 

\noindent On $Y^{\zeta}$ define the skew symmetric bilinear form ${\Omega }^{\zeta }$ by 
${\Omega }^{\zeta }(y,y') = \nu_{\zeta }([y,y'])$ for $y$, $y' \in Y^{\zeta }$. Here ${\nu }_{\zeta }$ is a complex linear function on $Z$, which is $1$ on $X_{\zeta }$ and $0$ on a complement to 
${\spann}_{\C}\{ X_{\zeta } \} $ in 
$Z$. \bigskip 

\noindent \textbf{Lemma 2.5}. For each $\zeta \in \mathcal{Z}$, $(Y^{\zeta }, {\Omega }_{\zeta })$ is a symplectic 
vector space. \medskip 

\noindent \textbf{Proof}. From the definition of ${\mathrm{R}}_{\zeta }$ and $Y^{\zeta }$ we have 
\begin{equation}
Y^{\zeta } = \bigoplus_{\widehat{\alpha } \in {\mathrm{R}}_{\zeta }} {\spann }_{\C} 
\{ X_{\widehat{\alpha }}, X_{\zetaminus \widehat{\alpha }} \} . 
\label{eq-onestardagger}
\end{equation} 
For $\widehat{\alpha }$, 
$\widehat{\beta } \in {\mathrm{R}}_{\zeta }$ with $\widehat{\beta } \ne \zetaminus \widehat{\alpha }$ 
we get 
\begin{align}
{\Omega }_{\zeta }(X_{\widehat{\alpha }}, X_{\widehat{\beta }} )  = 
{\nu }_{\zeta } ( [ X_{\widehat{\alpha }}, X_{\widehat{\beta }} ] ) 
& = {\nu }_{\zeta}(X_{{\zeta }^{\prime }}), \quad \parbox[t]{2in}{for some 
${\zeta }^{\prime } \in \mathcal{Z} \setminus \{ \zeta \} $, which implies $X_{{\zeta }^{\prime }} 
\not\in {\spann }_{\C}\{ X_{\zeta} \} $} 
\notag \\
& = 0, \quad \mbox{by definition of ${\nu }_{\zeta }$.} \notag 
\end{align}
In addition, 
\begin{displaymath}
{\Omega }_{\zeta }(X_{\widehat{\alpha }}, X_{\zetaminus \widehat{\alpha }} ) = 
{\nu }_{\zeta } ( [ X_{\widehat{\alpha }}, X_{\zetaminus \widehat{\alpha }} ] ) 
 = {\nu }_{\zeta}(X_{{\zeta }}) = 1;
\end{displaymath}
whereas ${\Omega }_{\zeta }(X_{\widehat{\alpha }}, X_{\widehat{\alpha }} ) = 
{\nu }_{\zeta } ( [ X_{\widehat{\alpha }}, X_{\widehat{\alpha }} ] ) = 0$ and 
\begin{align}
{\Omega }_{\zeta }(X_{\zetaminus \widehat{\alpha }}, X_{\zetaminus \widehat{\alpha }} ) & = 
{\nu }_{\zeta } ( [ X_{\zetaminus \widehat{\alpha }}, X_{\zetaminus \widehat{\alpha }} ] ) = 0. 
\notag 
\end{align}
Thus for every $\widehat{\alpha } \in {\mathcal{R}}_{\zeta }$ the space ${\pi }_{\widehat{\alpha }} = {\spann}_{\C} \{ X_{\widehat{\alpha }}, X_{\zetaminus \widehat{\alpha }} \} $ is a symplectic subspace of 
$(Y^{\zeta} , {\Omega }_{\zeta })$. From (\ref{eq-onestardagger}) it follows that $Y^{\zeta}$ is a direct sum of 
${\Omega }_{\zeta}$-symplectic subspaces ${\pi }_{\widehat{\alpha }}$, $\widehat{\alpha } \in {\mathcal{R}}_{\zeta}$. Indeed, these subspaces are ${\Omega }_{\zeta }$-perpendicular. Hence $(Y^{\zeta }, {\Omega }_{\zeta})$ is a symplectic vector space. \hfill $\square $ \medskip 

\vspace{-.15in}\noindent By construction for every $\zeta \in \mathcal{Z}$ 
and every $\widehat{\alpha } \in {\mathcal{R}}_{\zeta}$ the vectors 
$\{ X_{\widehat{\alpha }}, \, X_{\mbox{\tiny $\zetaminus$} \widehat{\alpha }},$ $X_{\zeta } \} $ 
span a $3$-dimensional Heisenberg subalgebra ${\mathrm{h}}^{\zeta }_3$ of $\widetilde{\mathfrak{g}}$, namely, 
$[X_{\widehat{\alpha }}, X_{\zeta }] = 0$, $[X_{\mbox{\tiny $\zetaminus$} \widehat{\alpha }}, X_{\zeta }] = 0$, 
and $[X_{\widehat{\alpha }}, X_{\mbox{\tiny $\zetaminus$} \widehat{\alpha }} ] = X_{\zeta }$. We say 
that $\onehalf \dim Y^{\zeta }$ is the \emph{multiplicity} of ${\mathrm{h}}^{\zeta }_3$ in 
$\widetilde{\mathfrak{g}}$. Since ${\pi }_{\widehat{\alpha }} = {\spann }_{\C} \{ X_{\widehat{\alpha }}, 
X_{\zetaminus \widehat{\alpha }} | \widehat{\alpha } \in {\mathcal{R}}_{\zeta } \} $ is a 
symplectic subspace of $(Y^{\zeta }, {\Omega }_{\zeta })$ and $Y^{\zeta }$ is the symplectically 
perpendicular direct sum of the subspaces ${\pi }_{\widehat{\alpha }}$ for each $\zeta \in \mathcal{Z}$, 
the root vectors $\{ X_{\widehat{\alpha }}, X_{\zetaminus \widehat{\alpha }},\, \alpha \in {\mathcal{R}}_{\zeta } ; 
\, X_{\zeta } \} $ are a basis for a Heisenberg subalgebra ${\mathrm{h}}_{2m+1}$ of the nilradical 
$\widetilde{\mathfrak{n}}$ of dimension $2m = \dim Y^{\zeta }$. Thus we have proved \medskip 

\noindent \textbf{Claim 2.6}  Suppose that $\widetilde{\mathfrak{g}} = \mathfrak{g} \oplus \widetilde{\mathfrak{n}}$ 
is a sandwich algebra with a nonabelian nilpotent radical $\widetilde{\mathfrak{n}}$. Let 
$\mathcal{Z} = {\{ {\zeta }_i \} }^p_{i=1}$. Then there is a $2m_i +1$-dimensional Heisenberg subalgebra 
${\mathrm{h}}^{{\zeta }_i}_{2m_i+1}$ of the sandwich algebra $\widetilde{\mathfrak{g}}$, where 
$m_i = \onehalf \dim Y^{{\zeta }_i}$ for $1 \le i \le p$.

\section{Special sandwich algebras}

In this section we look at a sandwich algebra which is a subalgebra of a semisimple Lie 
algebra. We call such a sandwich algebra \emph{special}. \medskip 
 
We now prove \medskip 

\noindent \textbf{Claim 3.1} Every special sandwich algebra $\widetilde{\mathfrak{g}} = \mathfrak{g} \oplus \widetilde{\mathfrak{n}}$ is a subalgebra of a semisimple algebra $\underline{\mathfrak{g}}$ whose rank is $1$ greater than the rank of $\mathfrak{g}$. \medskip 

\noindent \textbf{Proof.} First we construct the algebra $\underline{\mathfrak{g}}$. Let $\widetilde{\mathfrak{g}} = \mathfrak{g} \oplus \widetilde{\mathfrak{n}}$ be a special 
sandwich algebra, which is a subalgebra of a semisimple algebra $\overline{\mathfrak{g}}$. 
Let $\overline{\mathfrak{h}}$ be a Cartan subalgebra of $\overline{\mathfrak{g}}$ with associated root system 
$\overline{\Phi}$. Then 
$\mathfrak{h} = \overline{\mathfrak{h}} \cap \mathfrak{g} \subsetneqq \overline{\mathfrak{h}}$ 
is a Cartan subalgebra of $\mathfrak{g}$. Let $\Phi $ be the set of roots of $\mathfrak{g}$ associated to 
$\mathfrak{h}$. Because $\widetilde{\mathfrak{g}}$ is a sandwich algebra, ${\ad }_{\mathfrak{h}}$ is a 
maximal toral subalgebra of $\gl (\widetilde{\mathfrak{n}}, \C )$. Thus we can write $\widetilde{\mathfrak{n}} = 
\sum_{\widetilde{\alpha } \in \widetilde{\Phi}} \oplus {\widetilde{\mathfrak{g}}}_{\widetilde{\alpha }}$, 
where $\widetilde{\Phi}$ is the set of roots for $\widetilde{\mathfrak{n}}$ associated to 
$\mathfrak{h}$ and  
\begin{displaymath}
{\widetilde{\mathfrak{g}}}_{\widetilde{\alpha }} = \{ X \in \widetilde{\mathfrak{g}} \setrule 
[H,X] = \widetilde{\alpha }(H) X \, \, \mbox{for all $H \in \mathfrak{h}$} \}  \ne \{ 0 \} 
\end{displaymath}
is the root space associated to the root $\widetilde{\alpha } \in \widetilde{\Phi}$. Since ${\ad }_{\mathfrak{h}}$ is a maximal torus in $\gl (\widetilde{\mathfrak{n}}, \C)$, 
we have $\dim {\widetilde{\mathfrak{g}}}_{\widetilde{\alpha }} 
=1$ for every $\widetilde{\alpha } \in \widetilde{\Phi }$. Let $X_{\widetilde{\alpha }}$ be a basis of 
${\widetilde{\mathfrak{g}}}_{\widetilde{\alpha }}$. \medskip 

\noindent Because $\mathfrak{h} \subsetneqq \overline{\mathfrak{h}}$, we can 
write $\overline{\mathfrak{h}} = \mathfrak{h} \oplus \mathfrak{c}$. 
Let ${\alpha }^{\vvee}$ be the linear function on $\overline{\mathfrak{h}}$, which is equal 
to $\widetilde{\alpha }$ on $\mathfrak{h}$ and $0$ on the complement $\mathfrak{c}$ to $\mathfrak{h}$ in 
$\overline{\mathfrak{h}}$. Then ${\alpha }^{\vvee} \in {\overline{\mathfrak{h}}}^{\, \ast } = 
{\spann }_{\overline{\alpha} \in \overline{\Phi}}\{ \overline{\alpha } \} $. The root space 
${\overline{\mathfrak{g}}}_{{\alpha }^{\vvee}}$ in $\overline{\mathfrak{g}}$ with respect to 
the Cartan subalgebra $\overline{\mathfrak{h}}$ is equal to the root space 
${\widetilde{\mathfrak{g}}}_{\widetilde{\alpha }}$. Let ${\Phi }^{\vvee} = 
\{ {\alpha }^{\vvee} \in \overline{\Phi } \} $. \medskip 

Choose $\widetilde{\zeta } \in \widetilde{\Phi }$. Then ${\zeta }^{\vvee} \in 
\overline{\Phi }$ and ${\overline{\mathfrak{g}}}_{{\zeta }^{\vvee}} = {\widetilde{\mathfrak{g}}}_{\widetilde{\zeta }}$, 
which is spanned by $X_{{\zeta}^{\vvee}} = X_{\widetilde{\zeta}}$. Because $\overline{\Phi}$ is a root system 
for the semisimple algebra $\overline{\mathfrak{g}}$ and ${\zeta}^{\vvee}$ is a root, 
which is nonzero, we see that $-{\zeta }^{\vvee}$ is a root in $\overline{\Phi}$, whose rootspace 
${\overline{\mathfrak{g}}}_{-{\zeta}^{\vvee}}$ is spanned by 
the root vector $X_{-{\zeta}^{\vvee}}$. Since $\overline{\mathfrak{g}}$ is semisimple, the vector 
$H_{{\zeta }^{\vvee}} = [X_{{\zeta }^{\vvee}}, X_{-{\zeta }^{\vvee}}]$ lies in $\overline{\mathfrak{h}}$, see 
Humphreys \cite[p. 37]{humphreys}. Since $X_{{\zeta }^{\vvee}} \in \widetilde{\eta }$ and 
$\mathfrak{g} \cap \widetilde{\mathfrak{n}} = \{ 0 \} $, it follows that ${\zeta }^{\vvee} \notin 
\Phi $. Consequently, $H_{{\zeta }^{\vvee}} \notin \mathfrak{h}$. However, 
for every $H \in \mathfrak{h}$ we have 
\begin{align}
[H, H_{{\zeta }^{\vvee}} ] & = [H, [X_{{\zeta }^{\vvee}}, X_{-{\zeta }^{\vvee}} ] ] = 
[[H, X_{{\zeta }^{\vvee}}], X_{-{\zeta }^{\vvee}}] + [X_{{\zeta }^{\vvee}}, [ H, X_{-{\zeta}^{\vvee}}]]  \notag \\
& = {\zeta }^{\vvee}(H) [X_{{\zeta }^{\vvee}}, X_{-{\zeta }^{\vvee}}] - 
{\zeta }^{\vvee} (H) [X_{{\zeta }^{\vvee}}, X_{-{\zeta }^{\vvee}}] = 0. \notag 
\end{align}
So 
\begin{displaymath}
{\mathfrak{h}}^{\vvee} = \{ H, \, \, H \in \mathfrak{h}; \, H_{{\zeta }^{\vvee}} \} = \mathfrak{h} \oplus 
\spann \{ H_{{\zeta }^{\vvee}} \} 
\end{displaymath}
is an toral subalgebra of $\overline{\mathfrak{g}}$. Since ${\Phi }^{\vvee} = \Phi \amalg \{ \pm {\zeta }^{\vvee} \} $ 
is a root subsystem of the root system $\overline{\Phi }$ of $\overline{\mathfrak{g}}$, it follows that 
$({\mathfrak{h}}^{\vvee})^{\ast } = {\spann }_{\C }\{ {\alpha }^{\vvee} \in {\Phi }^{\vvee} \}$ and the 
Killing form on $\overline{\mathfrak{g}}$ restricted to ${\mathfrak{h}}^{\vvee}$ is nondegenerate. Because 
${\mathfrak{h}}^{\vvee}$ is a toral subalgebra of $\overline{\mathfrak{g}}$, we have 
\begin{equation}
\overline{\mathfrak{g}} = {\mathfrak{h}}^{\vvee } \oplus \sum_{{\alpha }^{\vvee} \in {\Phi }^{\vvee}} \oplus 
{\overline{\mathfrak{g}}}^{\vvee}_{{\alpha }^{\vvee}},
\label{eq-sec3one}
\end{equation} 
where 
\begin{displaymath}
{\overline{\mathfrak{g}}}^{\vvee}_{{\alpha}^{\vvee}} = \{ X \in \overline{\mathfrak{g}} \, \setrule \, [H^{\vvee},X] = 
{\alpha }^{\vvee }(H^{\vvee}) X, \, \, \mbox{for all $H^{\vvee} \in {\mathfrak{h}}^{\vvee} $} \} . 
\end{displaymath} 
Note that ${\overline{\mathfrak{g}}}^{\vvee}_{{\alpha }^{\vvee}} \ne \{ 0 \} $ for all 
${\alpha }^{\vvee } \in {\Phi }^{\vvee}$. Since 
${\mathfrak{h}}^{\vvee}$ is not necessarily a maximal toral subalgebra of 
$\overline{\mathfrak{g}}$, the subspaces ${\overline{\mathfrak{g}}}^{\vvee}_{{\alpha }^{\vvee}}$ may 
have dimension strictly greater than $1$. For each ${\alpha }^{\vvee} \in {\Phi }^{\vvee}$ choose a nonzero vector $X_{{\alpha }^{\vvee}} \in {\overline{\mathfrak{g}}}_{{\alpha }^{\vvee}}$, the root space of 
$\overline{\mathfrak{g}}$ corresponding to the root ${\alpha }^{\vvee}$ in $\overline{\Phi }$. Since 
${\mathfrak{h}}^{\vvee} \subseteq \overline{\mathfrak{h}}$, it follows that 
${\overline{\mathfrak{g}}}_{{\alpha }^{\vvee}} \subseteq {\overline{\mathfrak{g}}}^{\vvee}_{{\alpha }^{\vvee}}$. Let 
${\mathfrak{g}}^{\vvee}_{{\alpha }^{\vvee}} = {\spann }_{\C } \{ X_{{\alpha }^{\vvee}} \} $ and set 
\begin{equation}
{\mathfrak{g}}^{\vvee } = {\mathfrak{h}}^{\vvee} \oplus 
\sum_{{\alpha }^{\vvee} \in {\Phi }^{\vvee} } \oplus {\mathfrak{g}}^{\vvee}_{{\alpha }^{\vvee}}.
\label{eq-sec3two}
\end{equation}
Observe that ${\mathfrak{g}}^{\vvee}$ is a subalgebra of $\overline{\mathfrak{g}}$, because 
it inherits the bracket relations of $\overline{\mathfrak{g}}$ from the decompositon 
$\overline{\mathfrak{g}} = \overline{\mathfrak{h}} \oplus \sum_{\overline{\alpha } \in \overline{\Phi }} 
\oplus {\overline{\mathfrak{g}}}_{\overline{\alpha }}$. The algebra ${\mathfrak{g}}^{\vvee }$ clearly contains 
the special sandwich algebra $\widetilde{\mathfrak{g}}$ and hence is a candidate to be the desired 
algebra $\underline{\mathfrak{g}}$. To finish the proof of the claim we need to 
establish \medskip   

\noindent \textbf{Lemma 3.2} The Lie algebra ${\mathfrak{g}}^{\vvee }$ is semisimple with 
root system ${\Phi }^{\vvee}$ associated to the Cartan subalgebra ${\mathfrak{h}}^{\vvee}$. \medskip 

\noindent \textbf{Proof.} First we verify that the algebra ${\mathfrak{g}}^{\vvee}$ is semisimple. 
Our argument follows that of Humphreys \cite[p.100]{humphreys}. 
This amounts to showing that if $A$ is an abelian ideal $A$ of ${\mathfrak{g}}^{\vvee}$, then $A =\{ 0 \}$. 
Since ${\ad }_{{\mathfrak{h}}^{\vvee}} A \subseteq A$ using (\ref{eq-sec3two}) we get  
\begin{displaymath}
A = A \cap {\mathfrak{h}}^{\vvee}  \oplus \sum_{{\alpha }^{\vvee} \in {\Phi }^{\vvee}} 
({\mathfrak{g}}^{\vvee}_{{\alpha }^{\vvee}} \cap A). 
\end{displaymath}
Suppose that for some ${\alpha }^{\vvee} \in {\Phi }^{\vvee}$ we have ${\mathfrak{g}}^{\vvee}_{{\alpha }^{\vvee}} 
\cap A \ne \{0 \}$. Then ${\mathfrak{g}}^{\vvee}_{{\alpha }^{\vvee}}  \subseteq A$, 
since ${\mathfrak{g}}^{\vvee}_{{\alpha }^{\vvee}}  = \spann \{ X_{{\alpha }^{\vvee}} \}$. 
Because $A$ is an ideal in ${\mathfrak{g}}^{\vvee}$ we get 
$[{\mathfrak{g}}^{\vvee}_{{\alpha }^{\vvee}}, {\mathfrak{g}}^{\vvee}_{-{\alpha }^{\vvee}}] \subseteq 
[A, {\mathfrak{g}}^{\vvee}_{-{\alpha }^{\vvee}}] \subseteq A$. We now show that 
${\mathfrak{g}}^{\vvee}_{-{\alpha }^{\vvee}} \subseteq A$. Let $X_{-{\alpha }^{\vvee}} \in 
{\mathfrak{g}}^{\vvee}_{-{\alpha }^{\vvee}}$. Then $H_{{\alpha }^{\vvee}} = 
[X_{{\alpha }^{\vvee}}, X_{-{\alpha }^{\vvee}} ] \in {\mathfrak{h}}^{\vvee} $. 
But $[X_{{\alpha }^{\vvee}}, X_{-{\alpha }^{\vvee}}] \in A$ since $[{\mathfrak{g}}^{\vvee}_{{\alpha }^{\vvee}}, 
{\mathfrak{g}}^{\vvee}_{-{\alpha }^{\vvee}}] \subseteq A$. Thus 
$H_{{\alpha }^{\vvee}} \in A \cap {\mathfrak{h}}^{\vvee} $. Now 
${\alpha }^{\vvee} (H_{{\alpha }^{\vvee}}) X_{-{\alpha }^{\vvee}} = 
[H^{\vvee}_{{\alpha }^{\vvee}} , X_{-{\alpha }^{\vvee}} ] \in A$ and 
${\alpha }^{\vvee}(H_{{\alpha }^{\vvee}}) \ne 0$, since ${\Phi }^{\vvee}$ is a root system for 
${\mathfrak{g}}^{\vvee}$. Consequently, 
$X_{-{\alpha }^{\vvee}} \in A$. So ${\mathfrak{g}}^{\vvee}_{-{\alpha }^{\vvee}} \subseteq A$. 
But then $\{ H^{\vvee}_{{\alpha }^{\vvee}} , X_{{\alpha }^{\vvee}}, X_{-{\alpha }^{\vvee}} \}$ spans a nonabelian subalgebra of ${\mathfrak{g}}^{\vvee}$ contained in $A$, which is abelian. This is a contradiction. 
So ${\mathfrak{g}}^{\vvee}_{{\alpha }^{\vvee}} \cap A = \{ 0 \} $ for every ${\alpha }^{\vvee} \in {\Phi }^{\vvee}$. 
Hence $A = A \cap {\mathfrak{h}}^{\vvee} $. Whence $[{\mathfrak{g}}^{\vvee}_{{\alpha }^{\vvee}}, A] \subseteq
{\mathfrak{g}}^{\vvee}_{{\alpha }^{\vvee}} \cap A = \{ 0 \}$ for every 
${\alpha }^{\vvee} \in {\Phi }^{\vvee} $. So $A \subseteq \cap_{{\alpha }^{\vvee }\in {\Phi }^{\vvee} } \ker 
{\alpha }^{\vvee} $. But $\cap_{{\alpha }^{\vvee} \in {\Phi }^{\vvee}} \ker {\alpha }^{\vvee} = \{ 0 \} $, 
since $({\mathfrak{h}}^{\vvee })^{\ast } = \spann \{ {\alpha }^{\vvee } \setrule {\alpha }^{\vvee } \in {\Phi }^{\vvee} \}$ 
by construction of ${\mathfrak{h}}^{\vvee}$. Thus ${\mathfrak{g}}^{\vvee}$ is semisimple. \medskip 

Because $\dim {\mathfrak{g}}^{\vvee}_{{\alpha }^{\vvee}} = 1$ for all ${\alpha }^{\vvee} \in {\Phi }^{\vvee}$. 
it follows that ${\mathfrak{h}}^{\vvee }$ is a maximal toral subalgebra of ${\mathfrak{g}}^{\vvee}$ and 
hence is a Cartan subalgebra. \hfill {\footnotesize $\square $} \medskip

\section{Construction of very special sandwich algebras}

We say that a special sandwich Lie algebra $\widetilde{\mathfrak{g}}= \mathfrak{g} \oplus \widetilde{\mathfrak{n}}$, which is a subalgebra of the semisimple algebra $\underline{\mathfrak{g}}$ of rank $1$ more than the rank of 
$\mathfrak{g}$, is \emph{very special} if $\mathfrak{g}$ and $\underline{\mathfrak{g}}$ are simple Lie algebras. 
We now construct every very special sandwich algebra, which we list starting on page \pageref{list}. \medskip 

Because the rank of the simple Lie algebra $\underline{\mathfrak{g}}$ constructed in claim 3.1 is $1$ more than 
the rank of the simple Lie algebra $\mathfrak{g}$, the Dynkin diagram of a very special sandwich 
algebra $\mathfrak{g}$ of rank $\ell $ is obtained from the Dynkin diagram of the simple Lie 
algebra $\underline{\mathfrak{g}}$ of rank $\ell +1$ by removing a node, numbered $L$. Because 
$\mathfrak{g}$ is a simple Lie algebra its Dynkin diagram is connected. Thus the node $L$ must be an 
extremity of the Dynkin diagram of $\underline{\mathfrak{g}}$. Let $\underline{\mathfrak{h}}$ 
be a Cartan subalgebra of $\underline{\mathfrak{g}}$. Then 
$\mathfrak{h} = \underline{\mathfrak{h}} \cap \mathfrak{g}$ is a 
Cartan subalgebra of $\mathfrak{g}$. We now \emph{align} the root system $\mathcal{R}$ associated to 
Cartan subalgebra $\mathfrak{h}$ with the root system $\underline{\mathcal{R}}$ associated to the 
Cartan subalgebra $\underline{\mathfrak{h}}$. In other words, we find a vector 
$h^{\ast } \in \underline{\mathfrak{h}}$, which does not lie in $\mathfrak{h}$, such that 
the root system $\mathcal{R}$ is equal to ${\mathrm{R}}^0 = \{ \alpha \in \underline{\mathcal{R}} \setrule \, 
\alpha (h^{\ast }) =0 \} $.  \medskip 

\noindent \textbf{Claim 4.1} Let ${\{ h_i \}}^{\ell +1}_{i=1}$ be an ordered basis of the Cartan subalgebra 
$\underline{\mathfrak{h}}$ of $\underline{\mathfrak{g}}$. The vector $h^{\ast}$ is the generator 
of the $1$-dimensional $\Z $-submodule ${\underline{\mathfrak{h}}}_{\Z } = 
{\spann }_{\Z } {\{ h_i \} }^{\ell +1}_{i=1}$ defined by $\ker C^{(L)}_{\underline{\mathfrak{g}}}$. 
Here $C^{(L)}_{\underline{\mathfrak{g}}}$ is the $\ell \times (\ell +1)$ matrix formed by 
removing the $L^{\mathrm{th}}$ row of the $(\ell +1) \times (\ell +1) $ Cartan matrix $C_{\underline{\mathfrak{g}}}
\big( {\alpha }_i(h_j) \big)_{\! \stackrel{1 \le i \le \ell +1}{\mbox{$\scriptscriptstyle 1 \le j \le \ell +1$}}}$ of 
the Lie algebra $\underline{\mathfrak{g}}$ associated to the root system $\underline{\mathcal{R}}$, 
where $\underline{\Pi } = {\{ {\alpha }_i \}}^{\ell +1}_{i=1}$ is an ordered basis of $\underline{\mathfrak{h}}$ 
consisting of simple positive roots. \medskip 

\noindent \textbf{Proof.} First note that $\ker C^{(L)}_{\underline{\mathfrak{g}}}$ is a $\Z $-submodule 
of ${\underline{\mathfrak{h}}}_{\Z }$, since the entries of the Cartan matrix $C_{\underline{\mathfrak{g}}}$ 
are integers. Let $C^{(L,L)}_{\underline{\mathfrak{g}}}$ be the $\ell \times \ell $ matrix 
$C^{(L,L)}_{\underline{\mathfrak{g}}}$ formed by removing the $L^{\mathrm{th}}$ column of the $\ell \times (\ell +1)$ 
matrix $C^{(L)}_{\underline{\mathfrak{g}}}$. Then $C^{(L,L)}_{\underline{\mathfrak{g}}}$ is the Cartan matrix 
$C_{\mathfrak{g}}$ of the simple Lie algebra $\mathfrak{g}$. Since $C^{(L,L)}_{\underline{\mathfrak{g}}}$ is invertible, the rank of the matrix $C^{(L)}_{\underline{\mathfrak{g}}}$ is $\ell $. 
Thus $\dim \ker C^{(L)}_{\underline{\mathfrak{g}}} =1$. Hence there is a nonzero vector 
$h^{\ast } \in {\underline{\mathfrak{h}}}_{\Z }$ which generates 
the $\Z $-submodule $\ker C^{(L)}_{\underline{\mathfrak{g}}}$. We normalize $h^{\ast }$ so that 
${\alpha }_L (h^{\ast }) > 0$, where ${\alpha }_L$ is the $L^{\mathrm{th}}$ simple positive root in 
$\underline{\Pi }$.  \medskip 

Recall that ${\mathrm{R}}^0 = \{ \alpha \in \underline{\mathcal{R}} \setrule \, \alpha (h^{\ast }) =0 \} $. 
Let $\Pi = \{ {\alpha }_i \in \underline{\Pi }, \, i \ne L \} $ and let $\{ h_i, \, i \ne L \} $ be the corresponding 
basis of the Cartan subalgebra $\mathfrak{h}$ of $\mathfrak{g}$. Since $C_{\underline{\mathfrak{g}}}= 
{\big( {\alpha }_i (h_j) \big)}_{\! \stackrel{1 \le i \le \ell +1}{\mbox{$\scriptscriptstyle 1 \le j \le \ell +1$}}}$, 
by construction we have ${\alpha }_i(h^{\ast }) =0 $ for every $1 \le i \le \ell +1$ with $i \ne L$. In 
other words, $\Pi \subseteq {\mathrm{R}}^0$. Since ${\mathrm{R}}^0$ is closed under addition and 
$\Pi $ generates the root system $\mathcal{R}$ for $\mathfrak{g}$ corresponding to the Cartan subalgebra 
$\mathfrak{h} = \{ h_i, \, i \ne L \} $, it follows that $\mathcal{R} \subseteq {\mathrm{R}}^0$. But 
${\alpha }_L(h^{\ast }) >0$. So ${\alpha }_L \notin {\mathrm{R}}^0$. Hence ${\mathrm{R}}^0$ is generated 
by $\Pi $, that is, $\mathcal{R} = {\mathrm{R}}^0$, because $\underline{\Pi }$ generates 
$\underline{\mathcal{R}}$.    \hfill {\footnotesize $\square $} \medskip 

 \noindent \textbf{Fact 4.2} Let ${\mathrm{R}}^{-} = \{ \alpha \in \underline{\mathcal{R}} \setrule 
\alpha (h^{\ast }) < 0 \} $. Then ${\mathrm{R}}^{-}$ is the set of all roots in $\underline{\mathcal{R}}$ whose 
$L^{\mathrm{th}}$ component with respect to the ordered basis $\underline{\Pi }$ of $\underline{\mathcal{R}}$ is 
negative. \medskip
\vspace{-.1in} 

\noindent \textbf{Proof.} Let $\alpha = \sum^{\ell +1}_{i =1} a_i {\alpha }_i \in \underline{\mathcal{R}} 
\cap {\mathrm{R}}^{-} \subseteq {\underline{\mathfrak{h}}}_{\R }$ with $a_i \in \R $. Then 
\begin{equation}
0 > \alpha (h^{\ast }) = \sum_{i\ne L} a_i {\alpha }_i(h^{\ast }) + a_L {\alpha }_L(h^{\ast }) = a_L {\alpha }_L(h^{\ast }), 
\label{eq-sec4one}
\end{equation}
since ${\alpha }_i, \, i \ne L$ lie in ${\mathrm{R}}^0$. But $a_L \ne 0$. For if $a_L =0$, then 
$\alpha = \sum_{i\ne L}a_i {\alpha }_i \in {\mathrm{R}}^0$, that is, $\alpha (h^{\ast }) =0$. This 
contradicts the hypothesis that $\alpha \in {\mathrm{R}}^{-}$. If $a_L >0$, then from (\ref{eq-sec4one}) 
it follows that ${\alpha }_L(h^{\ast}) < 0$. But this contradicts the normalization of the vector 
$h^{\ast }$. Consequently, $a_L <0$. \hfill {\footnotesize $\square $} \medskip  

Let $\widetilde{\mathfrak{g}} = \mathfrak{h} \oplus \sum_{\underline{\alpha } 
\in {\mathrm{R}}^0\, \amalg \, {\mathrm{R}}^{-}} {{\underline{\mathfrak{g}}}_{\underline{\alpha }}}$ and 
$\widetilde{\mathfrak{n}} = \sum_{\underline{\alpha } \in {\mathrm{R}}^{-}} \oplus 
{{\underline{\mathfrak{g}}}_{\underline{\alpha }}}$. Then $\widetilde{\mathfrak{g}}$ is a 
subalgebra of $\underline{\mathfrak{g}}$ and $\widetilde{\mathfrak{g}} = 
\mathfrak{g} \oplus \widetilde{\mathfrak{n}}$.  Using the next few lemmas we will show that 
$\widetilde{\mathfrak{g}}$ is a very special sandwich algebra. \medskip 

\noindent \textbf{Lemma 4.3.} $\widetilde{\mathfrak{n}}$ is a subalgebra of $\widetilde{\mathfrak{g}}$. \medskip 

\noindent \textbf{Proof.} This amounts to showing that 
if ${\underline{\alpha }}_1$, ${\underline{\alpha }}_2 \in {\mathrm{R}}^{-}$ and 
${\underline{\alpha }}_1 + {\underline{\alpha }}_2 \in \underline{\mathcal{R}}$, 
then ${\underline{\alpha }}_1 + {\underline{\alpha }}_2 \in {\mathrm{R}}^{-}$, because 
\begin{displaymath}
[ \widetilde{\mathfrak{n}}, \widetilde{\mathfrak{n}} ] = \sum_{\stackrel{{\underline{\alpha }}_1 + 
{\underline{\alpha }}_2 \in \underline{\mathcal{R}}}{{\underline{\alpha }}_1, 
{\underline{\alpha }}_2 \in {\mathrm{R}}^{-}}} \oplus 
{\underline{\mathfrak{g}}}_{{\underline{\alpha }}_1 + {\underline{\alpha }}_2}. 
\end{displaymath}
To see this for every 
$H \in \mathfrak{h}$ we compute 
\begin{displaymath}
({\underline{\alpha }}_1+{\underline{\alpha }}_2)(H) = {\underline{\alpha }}_1(H) +{\underline{\alpha }}_2(H) 
< 0, 
\end{displaymath}
since ${\underline{\alpha }}_1(H) <0$ and ${\underline{\alpha }}_2(H) <0$ because 
${\underline{\alpha }}_1$ and ${\underline{\alpha }}_2$ lie in ${\mathrm{R}}^{-}$. \hfill {\footnotesize $\square $} \medskip 

\noindent \textbf{Lemma 4.4} $\widetilde{\mathfrak{n}}$ is an ideal of the Lie algebra $\widetilde{\mathfrak{g}}$. 
\medskip 

\noindent \textbf{Proof.} We need to show that 
\begin{equation}
[\widetilde{\mathfrak{g}}, \widetilde{\mathfrak{n}} ] 
\subseteq \widetilde{\mathfrak{n}} .
\label{eq-twoverynw}
\end{equation} 
Because $\widetilde{\mathfrak{n}}$ is a subalgebra, we have $[ \widetilde{\mathfrak{n}}, 
\widetilde{\mathfrak{n}} ]  \subseteq \widetilde{\mathfrak{n}}$. Thus to prove (\ref{eq-twoverynw}) we 
need only show that $[\mathfrak{g}, \widetilde{\mathfrak{n}} ] \subseteq \widetilde{\mathfrak{n}}$, 
because $\widetilde{\mathfrak{g}} = \mathfrak{g} \oplus \widetilde{\mathfrak{n}}$. From the definition of 
$\widetilde{\mathfrak{n}}$ we see that 
$[ \mathfrak{h}, \widetilde{\mathfrak{n}} ] \subseteq \widetilde{\mathfrak{n}}$. Now let 
$\underline{\alpha } \in {\mathfrak{h}}^{\ast } = {\mathrm{R}}^{\circ}$ and 
$\underline{\beta } \in {\mathrm{R}}^{-}$. Then 
\begin{displaymath}
(\underline{\alpha } +\underline{\beta })(H^{\ast }) = \underline{\alpha }(H^{\ast }) + \underline{\beta }(H^{\ast }) = 
\underline{\beta }(H^{\ast }) < 0.
\end{displaymath}
So $\underline{\alpha }+ \underline{\beta } \in {\mathrm{R}}^{-}$. This gives 
\begin{displaymath}
\Big[ \sum_{\underline{\alpha } \in {\mathrm{R}}^{\circ}} \oplus {{\underline{\mathfrak{g}}}_{\underline{\alpha }}} , 
\sum_{\underline{\beta } \in {\mathrm{R}}^{-}} \oplus {{\underline{\mathfrak{g}}}_{\underline{\beta }}} \Big] = 
\sum_{\underline{\alpha } + \underline{\beta } \in {\mathrm{R}}^{-} } \oplus 
{\underline{\mathfrak{g}}}_{\underline{\alpha } + \underline{\beta }} \subseteq \widetilde{\mathfrak{n}}. 
\end{displaymath}
Consequently, $[\mathfrak{g}, \widetilde{\mathfrak{n}} ] \subseteq \widetilde{\mathfrak{n}}$, 
since $\mathfrak{g} = \mathfrak{h} \oplus \sum_{\underline{\alpha } \in {\mathrm{R}}^{\circ}} \oplus {{\underline{\mathfrak{g}}}_{\underline{\alpha }}} $. \hfill {\footnotesize $\square $} 
\newpage 

\noindent \textbf{Lemma 4.5.} $\widetilde{\mathfrak{n}}$ is the nilpotent radical of 
$\widetilde{\mathfrak{g}}$. \medskip 

\noindent \textbf{Proof.} Suppose that $\mathfrak{l}$ is a nilpotent ideal of 
$\widetilde{\mathfrak{g}}$ which properly contains $\widetilde{\mathfrak{n}}$. 
Then $\mathfrak{l} \cap \mathfrak{g} $ is an ideal of $\mathfrak{g}$. Because 
$\mathfrak{g}$ is a simple Lie algebra, either $\mathfrak{l} \cap \mathfrak{g} = 
\mathfrak{g}$ or $\mathfrak{l} \cap \mathfrak{g} = \{ 0 \} $. If $\mathfrak{l} \cap \mathfrak{g} = 
\mathfrak{g}$, then $\mathfrak{g} \subseteq \mathfrak{l}$. But $\mathfrak{h} 
\subseteq \mathfrak{g}$ and $\mathfrak{h}$ has only semisimple elements. This contradicts the fact that 
$\mathfrak{l}$ is a nilpotent ideal. Therefore $\mathfrak{l} \cap \mathfrak{g} = \{ 0 \} $, which 
contradicts the fact that $\mathfrak{l}$ properly contains $\widetilde{\mathfrak{n}}$ and 
$\widetilde{\mathfrak{g}} = \mathfrak{g} \oplus \widetilde{\mathfrak{n}}$. Therefore 
the ideal $\mathfrak{l}$ does not exist. Hence $\widetilde{\mathfrak{n}}$ is a 
maximal nilpotent ideal in $\widetilde{\mathfrak{g}}$. \hfill {\footnotesize $\square $} \medskip 

Let $\widehat{\mathcal{R}} = \{ \alpha |\mathfrak{h} \setrule \alpha \in {\mathrm{R}}^{-} \} $. If 
$\widehat{\alpha }$ and $\widehat{\beta} \in \widehat{\mathcal{R}}$, the sum $\widehat{\alpha } 
+ \widehat{\beta }$ is an element of $\widehat{\mathcal{R}}$ if there exists an $\underline{\alpha}$ and 
$\underline{\beta }$ in $\underline{\mathcal{R}}$ with $\widehat{\alpha } =\underline{\alpha }|\mathfrak{h}$ and 
$\widehat{\beta } = \underline{\beta }|\mathfrak{h}$ such that 
$\widehat{\alpha } + \widehat{\beta } = (\underline{\alpha }+ \underline{\beta }) |\mathfrak{h}$. If no such 
$\underline{\alpha }$ and $\underline{\beta }$ in $\underline{\mathcal{R}}$ exist, then the sum 
$\widehat{\alpha } + \widehat{\beta }$ is not defined. In other words, $\widehat{\alpha } + \widehat{\beta } 
\notin \widehat{\mathcal{R}}$. \medskip 

\noindent \textbf{Lemma 4.6} $\widetilde{\mathfrak{n}}$ is a sandwich, if all additive relations 
in $\widehat{\mathcal{R}}$ are given by addition of pairs of roots in $\widehat{\mathcal{R}}$. \medskip 

\noindent \textbf{Proof.} Let $\widehat{\mathcal{R}} = {\{ {\widehat{\alpha }}_j \} }^N_{j=1}$.  
\par \noindent \textsc{case} 1. $\widehat{\mathcal{R}}$ has no additive relations. In other words, 
for every $i,j \in \{ 1, \ldots, N \} $ we have ${\widehat{\alpha }}_i + {\widehat{\alpha }}_j \notin 
\widehat{\mathcal{R}}$. Thus the bracket of the root vectors $X_{{\widehat{\alpha }}_i}$ and 
$X_{{\widehat{\alpha }}_j}$ vanishes. Consequently the nilradical $\widetilde{\mathfrak{n}} = 
{\spann }_{\C } \{ X_{{\widehat{\alpha }}_i}, \, 1 \le i \le N \} $ is abelian and hence is a sandwich. \medskip 

\noindent \textsc{case} 2. $\widehat{\mathcal{R}}$ has additive relations. Find all pairs of roots in 
$\widehat{\mathcal{R}}$ whose sum is a root in $\widehat{\mathcal{R}}$. Sort these pairs into a 
finite number of classes ${\widetilde{\mathcal{R}}}_{{\widehat{\zeta}}_i}$, $i \in F$ whose elements sum to 
${\widehat{\zeta }}_i$. By hypothesis this gives a complete list of the additive relations among the roots in 
$\widehat{\mathcal{R}}$. Thus the sum of ${\widehat{\zeta}}_i$ and 
any root ${\widehat{\alpha }}_j$ in $\widehat{\mathcal{R}}$ is not a root in $\widehat{\mathcal{R}}$. 
Let $\mathcal{Z} = \{ {\widehat{\zeta}}_i, \, i \in F \} $. Then $Z = {\spann }_{\C } \{ X_{\widehat{\zeta}} , \, 
\widehat{\zeta } \in \mathcal{Z} \} $ is the center of the nilpotent radical $\widetilde{\mathfrak{n}} = 
{\spann }_{\C} \{ X_{{\widehat{\alpha }}_j}, \, 1 \le j \le N \} $. Let $\mathcal{Y} = \{ {\widehat{\alpha }}_j \ne {\widehat{\zeta }}_i \, \, \mbox{for every $i \in F$} \} $ and set $Y = {\spann }_{\C } \{ X_{\widehat{\alpha }}, \, \alpha \in \mathcal{Y} \}$. Then $\widetilde{\mathfrak{n}} = Z \oplus Y$. The above argument shows that each root in a pair of roots in ${\widetilde{\mathcal{R}}}_{\widehat{\zeta}}$ lies in $\mathcal{Y}$. Moreover, for every 
$\widehat{\alpha} \in \mathcal{Y}$ there is a $\widehat{\beta } \in \mathcal{Y}$ such that $\widehat{\alpha } + \widehat{\beta } = \widehat{\zeta }$ for some $\widehat{\zeta } \in \mathcal{Z}$. Hence $[Y, Y] = Z$. Now 
$[\widetilde{\mathfrak{n}}, \widetilde{\mathfrak{n}} ] = [Y+Z, Y+Z] = [Y,Y]$, since $Z$ is the center of 
$\widetilde{\mathfrak{n}}$ and $Y \subseteq \widetilde{\mathfrak{n}}$. Thus 
\begin{displaymath}
[ \widetilde{\mathfrak{n}}, [\widetilde{\mathfrak{n}}, \widetilde{\mathfrak{n}} ] ] = 
[Y+Z, [\widetilde{\mathfrak{n}}, \widetilde{\mathfrak{n}} ]] = [ Y, [\widetilde{\mathfrak{n}}, \widetilde{\mathfrak{n}} ] ] 
= [Y, [Y,Y]] = [Y,Z] =0, 
\end{displaymath}
that is, $\widetilde{\mathfrak{n}}$ is a sandwich. \hfill {\footnotesize $\square $} \medskip   

We now prove 
\newpage 

\noindent \textbf{Claim 4.7} $\widetilde{\mathfrak{g}} = \mathfrak{g} \oplus \widetilde{\mathfrak{n}}$ is a 
very special sandwich algebra. \medskip   

\noindent \textbf{Proof.} Since $\mathfrak{h}$ is a Cartan subalgebra of 
$\mathfrak{g}$ and $\widetilde{\mathfrak{n}}$ is an ideal of $\widetilde{\mathfrak{g}}$, 
${\{ {\ad }_H \} }_{H \in \mathfrak{h}}$ is a family of commuting linear mappings 
of $\widetilde{\mathfrak{n}}$ into itself. For every $H \in \mathfrak{h}$, the linear 
map ${\ad }_H$ preserves ${\underline{\mathfrak{g}}}_{\underline{\alpha }} $ for every 
$\underline{\alpha } \in {\mathrm{R}}^{-}$. From the fact that 
${\dim }_{\C } \, {\underline{\mathfrak{g}}}_{\underline{\alpha }} = 1$, it follows that 
${\ad }_H$ is semisimple and the family ${\{ {\ad }_H \} }_{H \in \mathfrak{h}}$ is 
maximal. This together with the facts $\widetilde{\mathfrak{g}} = \mathfrak{g} \oplus \widetilde{\mathfrak{n}} 
\subseteq \underline{\mathfrak{g}}$, 
where $\widetilde{\mathfrak{n}}$ is its nilpotent radical, $\mathfrak{g}$ is a simple Lie algebra,  
$\widetilde{\mathfrak{n}}$ is a sandwich, and $\underline{\mathfrak{g}}$ is a simple Lie algebra 
shows that $\widetilde{\mathfrak{g}}$ is a very special sandwich algebra. \hfill {\footnotesize $\square $} \medskip 

\noindent The following example shows omitting the root number $2$ in the Dynkin diagram of 
${\mathbf{E}}_8$ and following the above construction does \emph{not} give a very special sandwich algebra.  \medskip 

\noindent \textbf{Example.} Let $\underline{\mathfrak{g}} = {\mathbf{E}}_8$ with root system 
\begin{align}
\underline{\mathcal{R}} & = \big\{ \pm ({\varepsilon}_i \pm {\varepsilon }_j), \, 1\le i < j \le 8; \, \, 
\onehalf \sum^8_{i=1}(-1)^{k(i)}{\varepsilon}_i, \notag \\ 
&\hspace{.5in} \mbox{\parbox[t]{3.25in}{where $k(i) = 0$ or $1$ and $\sum^8_{i=1} k(i) \in \{ 0,2,4,6, 8\} $} } \, \big\} .\notag
\end{align}
$\underline{\mathcal{R}}$ has a subset positive simple roots
\begin{displaymath}
\underline{\Pi } = \big\{ \onehalf ({\varepsilon }_1 - {\varepsilon }_2 - \cdots -{\varepsilon }_7 
+{\varepsilon }_8); \, \, {\varepsilon }_1+{\varepsilon }_2; \, \, 
{\varepsilon }_{i+1} - {\varepsilon }_i, \, 1 \le i \le 6 \big\} . 
\end{displaymath}
The Cartan subalgebra corresponding to the root system $\underline{\mathcal{R}}$ is 
\begin{displaymath}
\underline{\mathfrak{h}} = \big\{ \onehalf (H_1-H_2- \cdots -H_7+H_8); \, \, 
H_1+H_2; \, \, H_{i+1}-H_i, \, 1\le i \le 6 \big\} , 
\end{displaymath}
where ${\varepsilon }_i(H_j) = {\delta }_{ij}$. Removing the root ${\varepsilon }_1+
{\varepsilon}_2$ from $\underline{\Pi }$ corresponds to removing the node 
$2$ from the Dynkin diagram of ${\mathbf{E}}_8$. This gives the set of positive simple roots 
\begin{displaymath}
\Pi = \big\{ \onehalf ({\varepsilon }_1 - {\varepsilon }_2 - \cdots -{\varepsilon }_7 
+{\varepsilon }_8);  \, \, {\varepsilon }_{i+1} - {\varepsilon }_i, \, 1 \le i \le 6 \big\} ,
\end{displaymath}
whose Dynkin diagram is that of the simple Lie algebra $\mathfrak{g} = {\mathbf{A}}_7$. $\Pi $ generates the closed subsystem of the set of positive roots of $\underline{\mathcal{R}}$ given by 
\begin{align}
& \big\{ {\varepsilon }_i - {\varepsilon }_j, \, 1 \le i < j \le 7; \, \, 
\onehalf ({\varepsilon }_1 + {\varepsilon }_2 -{\varepsilon }_3 - \cdots -{\varepsilon }_7  +{\varepsilon }_8) , \notag \\
&\hspace{.15in} \, \onehalf ({\varepsilon }_1 - {\varepsilon }_2 + {\varepsilon }_3 - 
{\varepsilon }_4 - \cdots -{\varepsilon }_7 +{\varepsilon }_8) , \,  
 \ldots \, , \, \onehalf ({\varepsilon }_1 - {\varepsilon }_2 -  \cdots - 
{\varepsilon }_6 + {\varepsilon }_7 +{\varepsilon }_8 ) \big\} , \notag 
\end{align} 
which is the set of positive roots ${\mathcal{R}}^{+}$ for $\mathfrak{g}$. A computation following the proof 
of claim 4.1 gives $H^{\ast } = -\sum^7_{i=1}H_i - 5H_8$. Consequently, ${\mathrm{R}}^{\circ} = \big\{ \underline{\alpha } \in \underline{\mathcal{R}} \, \setrule \, \underline{\alpha }(H^{\ast }) = 0 \big\} = \mathcal{R}$ 
and ${\mathrm{R}}^{-}  = \big\{  \underline{\alpha } \in \underline{\mathcal{R}} \, 
\setrule \, \underline{\alpha }(H^{\ast }) < 0  \big\} $. So  
\clearpage
\mbox{}\vspace{-1.5in}
\noindent
\begin{align}
&{\mathrm{R}}^{-} = \big\{ \chi ; \, z_i, \, 1 \le i \le 7; \, {\zeta }_i, \, 1 \le j \le 7;  \, {\zeta }_{ij}, \, 1 \le i < j \le 7;  \notag \\
& \hspace{.5in} {\delta }_i, \, 2 \le i \le 7; \, {\delta }_{ij}, \, 2 \le i < j \le 7; \, 
{\delta }_{ijk}, \, 2 \le i < j < k \le 7; \notag \\
& \hspace{.75in} \, {\delta }_{ijkn}, \, 2 \le i < j < k < n  \le 7 \big\} . \notag 
\end{align}
Here 
\begin{displaymath}
\begin{array}{rl}
\chi & = \onehalf ( {\varepsilon }_1 + \sum^7_{\ell =1}{\varepsilon}_{\ell } + 
{\varepsilon }_8), \, \, 
z_i = {\varepsilon }_i + {\varepsilon }_8, \, \, 
{\zeta }_i  = -{\varepsilon }_i + {\varepsilon }_8, \, \, 
{\zeta }_{ij}  = {\varepsilon }_i + {\varepsilon }_j  \\
\rule{0pt}{12pt} {\delta}_i & = \onehalf ( -{\varepsilon }_1 + \sum^7_{\ell =1} a_{\ell } {\varepsilon}_{\ell } + 
{\varepsilon }_8) \\ 
&\hspace{.5in}\mbox{where $a_{\ell } =-1$ if $\ell = i$ and $=1$ otherwise} \\
{\delta }_{ij} & = \onehalf ( {\varepsilon }_1 + \sum^7_{\ell =1} a_{\ell } {\varepsilon}_{\ell } + 
{\varepsilon }_8) \\ 
&\hspace{.5in}\mbox{where $a_{\ell } =-1$ if $\ell = i \, \, \& \, \, \ell = j$ and $=1$ otherwise} \\
{\delta }_{ijk} & = \onehalf ( - {\varepsilon }_1 + \sum^7_{\ell =1} a_{\ell } {\varepsilon}_{\ell } + 
{\varepsilon }_8) \\
&\hspace{.5in}\parbox[t]{3.5in}{where $a_{\ell } =-1$ if $\ell = i \, \, \& \, \, \ell = j \, \, \& \, \, \ell =k$ and $=1$ \\ otherwise} \\ 
\rowspace {\delta }_{ijkn} & = \onehalf ( {\varepsilon }_1 + \sum^7_{\ell =1} a_{\ell } {\varepsilon}_{\ell } + 
{\varepsilon }_8) \\
&\hspace{.5in}\parbox[t]{3.5in}{where $a_{\ell } =-1$ if $\ell = i \, \, \& \, \, \ell = j \, \, \& 
\, \, \ell =k \, \, \& \, \, \ell = n$ and $=1$ otherwise.}
\end{array}
\end{displaymath} 
Then $\widehat{\mathcal{R}} = \{ \widehat{\alpha } = \underline{\alpha }|\mathfrak{h} \, \setrule 
\, \underline{\alpha } \in {\mathrm{R}}^{-} \} $. Since ${\widehat{\delta }}_{2345} + {\widehat{\delta }}_{67} = {\widehat{z}}_1$ and 
${\widehat{\delta }}_{2367} +{\widehat{\zeta }}_{23} = {\widehat{\delta}}_{67}$, 
we get ${\widehat{\delta }}_{2345} + ({\widehat{\delta }}_{2367} + {\widehat{\zeta }}_{23}) = 
{\widehat{z}}_1$. This shows that $[ {\widehat{\mathfrak{g}}}_{{\widehat{\delta }}_{2345}}, 
[ {\widehat{\mathfrak{g}}}_{{\widehat{\delta }}_{2367}} , 
{\widehat{\mathfrak{g}}}_{{\widehat{\zeta }}_{23}} ] ] = 
{\widehat{\mathfrak{g}}}_{{\widehat{z}}_1}$. Thus $\widetilde{\mathfrak{n}} = \sum_{\widehat{\alpha } \in \widehat{\mathcal{R}}} \oplus {\widehat{\mathfrak{g}}}_{\widehat{\alpha }}$ is \emph{not} a sandwich. Hence $\widetilde{\mathfrak{g}} = 
\mathfrak{g} \oplus \widetilde{\mathfrak{n}} $ is not a sandwich algebra. \hfill {\footnotesize $\square $} 

\section*{List of very special sandwich algebras}
\label{list} 

The above example shows that we must verify that the nilradical is a sandwich. A case by case check shows that the roots $\widehat{\zeta }$ with various indices in $\widehat{\mathcal{R}}$, whic appear in the entries of the list, 
comprise $\mathcal{Z}$, whose corresponding root vectors span $Z$, the center of the nilpotent radical. Hence the additive relations given in the entries of the following list is complete. From lemma 4.6 it follows that the nilradical is a sandwich. Hence we have listed all very special sandwich algebras. \medskip 

In the list of very special sandwich algebras we use the notation 
${\{ {\varepsilon }_i \}}^n_{i=1}$ for the standard basis of $({\R }^n)^{\ast} $ dual to the standard 
basis ${\{ e_i \} }^n_{i=1}$ of ${\R }^n$; ${\mathrm{h}}_{2\ell +1}$ for 
the Lie algebra of the $2p +1$ dimensional Heisenberg group; $Z_m$ is the $m$ dimensional 
center of the nilradical, which is spanned by the root vectors corresponding to the roots in 
$\widehat{\mathcal{R}}$. \medskip

We note that because of Dynkin diagram automorphisms of 
$\underline{\mathfrak{g}}$ the \linebreak 
following sandwich algebras are isomorphic: ${\widetilde{\mathbf{A}}}^1_{\ell +1}$ 
and ${\widetilde{\mathbf{A}}}^2_{\ell +1}$; ${\widetilde{\mathbf{B}}}^1_{2}$ and ${\widetilde{\mathbf{C}}}^2_{2}$; 
${\widetilde{\mathbf{B}}}^2_{2}$ and ${\widetilde{\mathbf{C}}}^1_{2}$; ${\widetilde{\mathbf{D}}}^1_{4}$, 
${\widetilde{\mathbf{D}}}^2_{4}$, and ${\widetilde{\mathbf{D}}}^3_{4}$; ${\widetilde{\mathbf{D}}}^2_{\ell +1}$ 
and ${\widetilde{\mathbf{D}}}^3_{\ell +1}$ when $\ell \ge 4$; and ${\widetilde{\mathbf{E}}}^1_{6}$ and 
${\widetilde{\mathbf{E}}}^2_{6}$. %

\begin{tabular}{lll}
$1$.  & {\large ${\mathbf{A}}_{\ell +1}$}: \hspace{-1.25in} \begin{tabular}{l}
\setlength{\unitlength}{.45mm}
\begin{picture}(0,15)(-60,-3)
\begin{drawjoin}
\jput(0,5){\circle*{2}}
\jput(20,5){\circle*{2}} 
\jput(40,5){\circle*{2}}
\jput(60,5){\circle*{.0001}}
\end{drawjoin}
\begin{dottedjoin}{3}
\jput(60,5){\circle*{.0001}}
\jput(80,5){\circle*{.0001}}
\end{dottedjoin}
\begin{drawjoin}
\jput(80,5){\circle*{.0001}}
\jput(100,5){\circle*{2}}
\jput(120,5){\circle*{2}}
\end{drawjoin}
\put(0,0){\makebox(0,0){\small $1$}}
\put(20,0){\makebox(0,0){\small $2$}} 
\put(40,0){\makebox(0,0){\small $3$}}
\put(100,0){\makebox(0,0){\small $\ell$}}
\put(120,0){\makebox(0,0){\small $\ell +1$}}
\end{picture} 
\end{tabular} &  \\
& positive simple roots: $\underline{\Pi } = \{ {\alpha }_i = {\varepsilon }_i -{\varepsilon }_{i+1}, \, 1 \le  i\le \ell +1 \} $ & \\
& Cartan subalgebra: $\underline{\mathfrak{h}} = {\spann }_{\C }\{ h_i = e_i -e_{i+1}, \, 1 \le i \le \ell +1 \}$ & \\
\rule{0pt}{26pt} & \hspace{.25in} $C_{{\mathbf{A}}_{\ell +1}} = 
{\big( {\alpha }_i(h_j) \big) }_{\! \stackrel{1 \le i \le \ell +1}{\mbox{$\scriptscriptstyle 1 \le j \le \ell +1$}}} =
\mbox{\tiny $\left( \begin{array}{rrrrr}
2         & -1          & 0             &\cdots            & 0        \\
-1        &  2          &       -1      & \cdots           & 0         \\
\vdots  & \ddots   & \ddots      & \ddots          & \vdots  \\
0         & \cdots    &   -1          & 2                   & -1       \\ 
0         & \cdots    &  0             & -1                 & 2  
\end{array} \right) $}$ &  \\
\rule{0pt}{12pt}& \parbox[t]{4.5in}{is the Cartan matrix for the root system 
$\underline{\mathcal{R}} = {\mathbf{A}}_{\ell +1}$ \\ 
\rule{.25in}{0in} generated by $\underline{\Pi }$.} & \\
& positive roots: ${\underline{\mathcal{R}}}^{+} = \{ \sum^{j-1}_{k=i} {\alpha }_k =  
{\varepsilon }_i -{\varepsilon }_j, \, 1 \le i< j \le \ell +2 \} $ & \\
\\ 
1.1 ${\widetilde{\mathbf{A}}}^1_{\ell +1}$ &  {\large ${\mathbf{A}}^{(1)}_{\ell +1}$}:\hspace{-1.1in}\begin{tabular}{l}
\setlength{\unitlength}{.45mm}
\begin{picture}(0,15)(-60,-3)
\begin{drawjoin}
\jput(0,5){\circle*{2}}
\jput(20,5){\circle*{2}} 
\jput(40,5){\circle*{2}}
\jput(60,5){\circle*{.0001}}
\end{drawjoin}
\begin{dottedjoin}{3}
\jput(60,5){\circle*{.0001}}
\jput(80,5){\circle*{.0001}}
\end{dottedjoin}
\begin{drawjoin}
\jput(80,5){\circle*{.0001}}
\jput(100,5){\circle*{2}}
\jput(120,5){\circle*{2}}
\end{drawjoin}
\put(0,5){\makebox(0,0){\large $\times$}}
\put(0,0){\makebox(0,0){\small $1$}}
\put(20,0){\makebox(0,0){\small $2$}} 
\put(40,0){\makebox(0,0){\small $3$}}
\put(100,0){\makebox(0,0){\small $\ell$}}
\put(120,0){\makebox(0,0){\small $\ell +1$}}
\end{picture} 
\end{tabular} & \\
& positive simple roots: $\Pi = \{ {\alpha }_i, \, 2 \le i \le \ell +1 \} $ & \\
& Cartan subalgebra: $\mathfrak{h} = {\spann }_{\C} \{ h_i, \, 2 \le i \le \ell +1 \} $ & \\
\rule{0pt}{12pt} & \parbox[t]{4.5in}{Cartan matrix $C^{(1,1)}_{{\mathbf{A}}_{\ell +1}} = 
{ \big( {\alpha }_i(h_j) \big) }_{\! \stackrel{2 \le i \le \ell +1}{\mbox{$\scriptscriptstyle 2 \le j \le \ell +1$} } } $ for the root system $\mathcal{R} = {\mathbf{A}}_{\ell }$ \vspace{-3pt} \\
\rule{.25in}{0in} generated by $\Pi $.} & \\ 
& Set $h^{\ast} = \sum^{\ell +1}_{j=1} (\ell +2 -j) h_j = (\ell +1)e_1 - \sum^{\ell +2}_{i=2}e_i$. Then 
${\mathrm{R}}^0 = \mathcal{R}$. & \\
& ${\mathrm{R}}^{-} = \{ -( {\varepsilon }_1 - {\varepsilon }_j) = -\sum^{j-1}_{i=1} {\alpha }_i, \, 
2 \le j \le \ell +2 \} $ & \\
& nilradical roots: & \\
& $\widehat{\mathcal{R}} = 
\left\{ {\widehat{\alpha }}_{j-1} =  - \big( \sum^{j-1}_{i=1} {\alpha }_i \big) |\mathfrak{h} 
= -( {\varepsilon }_1 - {\varepsilon }_j)|\mathfrak{h}, \, 2 \le j \le \ell +2 \right. $ & \\
& additive relations in $\widehat{\mathcal{R}}$: none & \\
& nilradical structure: $Z_{\ell +1}$ & \\
\\ 
1.2 ${\widetilde{\mathbf{A}}}^2_{\ell +1}$ & {\large ${\mathbf{A}}^{(\ell +1)}_{\ell +1}$}:\hspace{-1.1in}\begin{tabular}{l}
\setlength{\unitlength}{.45mm}
\begin{picture}(0,15)(-60,-3)
\begin{drawjoin}
\jput(0,5){\circle*{2}}
\jput(20,5){\circle*{2}} 
\jput(40,5){\circle*{2}}
\jput(60,5){\circle*{.0001}}
\end{drawjoin}
\begin{dottedjoin}{3}
\jput(60,5){\circle*{.0001}}
\jput(80,5){\circle*{.0001}}
\end{dottedjoin}
\begin{drawjoin}
\jput(80,5){\circle*{.0001}}
\jput(100,5){\circle*{2}}
\jput(120,5){\circle*{2}}
\end{drawjoin}
\put(120,5){\makebox(0,0){\large $\times$}}
\put(0,0){\makebox(0,0){\small $1$}}
\put(20,0){\makebox(0,0){\small $2$}} 
\put(40,0){\makebox(0,0){\small $3$}}
\put(100,0){\makebox(0,0){\small $\ell$}}
\put(120,0){\makebox(0,0){\small $\ell +1$}}
\end{picture} 
\end{tabular} & \\
& positive simple roots: $\Pi = \{ {\alpha }_i, \, 1 \le i \le \ell \} $ & \\
& Cartan subalgebra: $\mathfrak{h} = {\spann }_{\C} \{ h_i, \, 1 \le i \le \ell \} $ & \\
& \parbox[t]{4.5in}{Cartan matrix $C^{(\ell +1,\ell +1)}_{{\mathbf{A}}_{\ell +1}} = 
{ \big( {\alpha }_i(h_j) \big) }_{\! \stackrel{1 \le i \le \ell }{\mbox{$\scriptscriptstyle 1 \le j \le \ell $}}} $ for the root \vspace{-3pt} \\ 
\rule{.25in}{0in} system $\mathcal{R} = {\mathbf{A}}_{\ell }$ generated by $\Pi $.} & \\ 
& \parbox[t]{4.5in}{Set $h^{\ast} = \sum^{\ell +1}_{j=1} j h_j = \sum^{\ell +1}_{i=1}e_i - (\ell +1)e_{\ell +2}$. Then 
${\mathrm{R}}^0 = \mathcal{R}$.} & \\
& ${\mathrm{R}}^{-} = \{ -( {\varepsilon }_j - {\varepsilon }_{\ell +2} ) = -\sum^{\ell +1}_{i=j} {\alpha }_i, 
\, 1 \le j \le \ell  \} $ & \\
& nilradical roots: & \\
& $\widehat{\mathcal{R}} = 
\left\{ {\widehat{\alpha }}_j =  - \big( \sum^{\ell +1}_{i=j} {\alpha }_i \big) |\mathfrak{h} 
= -( {\varepsilon }_j - {\varepsilon }_{\ell +2} )|\mathfrak{h}, \, 1 \le j \le \ell +1 \right. $ & \\
& additive relations in $\widehat{\mathcal{R}}$: none & \\
& nilradical structure: $Z_{\ell +1}$ & 
\end{tabular}
\newpage
\begin{tabular}{lll}
\vspace{-1.in} & & \\
2.& {\large ${\mathbf{B}}_{\ell +1} \, (\ell \ge 1)$:} \hspace{-.65in} \begin{tabular}{lll}
\begin{picture}(0,15)(-40,-3)
\begin{drawjoin}
\jput(0,5){\circle*{2.25}}
\jput(20,5){\circle*{2.25}} 
\jput(40,5){\circle*{2.25}}
\jput(60,5){\circle*{.0001}}
\end{drawjoin}
\begin{dottedjoin}{3}
\jput(60,5){\circle*{.0001}}
\jput(80,5){\circle*{.0001}}
\end{dottedjoin}
\begin{drawjoin}
\jput(80,5){\circle{.0001}}
\jput(100,5){\circle*{2.25}}
\end{drawjoin}
\put(100,4.5){\line(1,0){20}} 
\put(100,5.5){\line(1,0){20}}
\drawline(112,6.5)(109,5)
\drawline(112,3.5)(109,5)
\jput(120,5){\circle{2.25}}
\put(0,0){\makebox(0,-1){\small $1$}}
\put(20,0){\makebox(0,-1){\small $2$}} 
\put(40,0){\makebox(0,-1){\small $3$}}
\put(100,0){\makebox(0,-1){\small $\ell$}}
\put(125,0){\makebox(0,-1){\small $\ell +1$}}
\end{picture}  
\end{tabular} & \\ 
& simple positive roots: & \\
& $\underline{\Pi } = \{ {\alpha }_i = {\varepsilon }_i - {\varepsilon }_{i+1}, \, 1 \le i \le \ell ;\, 
{\alpha }_{\ell +1} = {\varepsilon }_{\ell +1} \} $ & \\
& Cartan subalgebra: & \\
& $\underline{\mathfrak{h}} = {\spann }_{\C }\{ h_i =e_i-e_{i+1}, \, 1 \le i \le \ell ; \, h_{\ell +1} = 2e_{\ell +1} \} $ & \\ 
\rule{0pt}{30pt}& \hspace{.15in} $C_{{\mathbf{B}}_{\ell +1}} = 
{\big( {\alpha }_i(h_j) \big) }_{\! \stackrel{1 \le i \le \ell +1}{\mbox{$\scriptscriptstyle 1 \le j \le \ell +1$}}} =
\mbox{\tiny $\left( \begin{array}{rrrrrr}
2          & -1          & 0             & 0                         & \cdots           & 0        \\
-1         &  2          &       -1      & 0                        & \cdots          & 0         \\
\vdots  & \ddots   & \ddots      & \ddots                 & \cdots           & \vdots  \\
          & \cdots    &-1              & 2                        & -1                   & 0           \\
  0         & \cdots    &  0            & -1                       & 2                   & -2  \\
  0        &  \cdots            &  0            & 0                       & -1                   & 2
\end{array} \right) $} $ &  \\
\rule{0pt}{12pt}& \parbox[t]{4.5in}{Cartan matrix for the root system 
$\underline{\mathcal{R}} = {\mathbf{B}}_{\ell +1}$ generated by $\underline{\Pi }$.} & \\
& positive roots: & \\
& ${\underline{\mathcal{R}}}^{+} = \left\{ \begin{array}{rl}
{\varepsilon }_i = & \hspace{-5pt} \sum^{\ell +1}_{k=i}{\alpha }_k, \, 1 \le i \le \ell +1; \\ 
\rule{0pt}{12pt} {\varepsilon }_i-{\varepsilon }_j  = & \hspace{-5pt}  \sum^{j-1}_{k=i}{\alpha }_k, \, 1 \le i < j \le \ell +1; \\
\rule{0pt}{12pt} {\varepsilon }_i+{\varepsilon }_j = &\hspace{-5pt}  \sum^{\ell +1}_{k=i}{\alpha }_k + 
\sum^{\ell +1}_{k=j}{\alpha }_k, \, 1 \le i < j \le \ell +1  
\end{array} \right. $ & \\
\\ 
2.1 ${\widetilde{\mathbf{B}}}^1_{\ell +1}$ &{\large ${\mathbf{B}}^{(1)}_{\ell +1}$}:\hspace{-.6in} \begin{tabular}{lll}
\begin{picture}(0,15)(-40,-3)
\begin{drawjoin}
\jput(0,5){\circle*{2.25}}
\jput(20,5){\circle*{2.25}} 
\jput(40,5){\circle*{2.25}}
\jput(60,5){\circle*{.0001}}
\end{drawjoin}
\begin{dottedjoin}{3}
\jput(60,5){\circle*{.0001}}
\jput(80,5){\circle*{.0001}}
\end{dottedjoin}
\begin{drawjoin}
\jput(80,5){\circle{.0001}}
\jput(100,5){\circle*{2.25}}
\end{drawjoin}
\put(100,4.5){\line(1,0){20}} 
\put(100,5.5){\line(1,0){20}}
\drawline(112,6.5)(109,5)
\drawline(112,3.5)(109,5)
\jput(120,5){\circle{2.25}}
\put(0,5){\makebox(0,0){\large $\times$}}
\put(0,0){\makebox(0,-1){\small $1$}}
\put(20,0){\makebox(0,-1){\small $2$}} 
\put(40,0){\makebox(0,-1){\small $3$}}
\put(100,0){\makebox(0,-1){\small $\ell$}}
\put(125,0){\makebox(0,-1){\small $\ell +1$}}
\end{picture}  
\end{tabular} & \\
& positive simple roots: $\Pi = \{ {\alpha }_i , \, 2 \le i \le \ell +1 \} $ & \\
& Cartan subalgebra: $\mathfrak{h} = {\spann }_{\C} \{ h_i , \, 2 \le i \le \ell +1 \} $ & \\
& \parbox[t]{4.5in}{Cartan matrix $C^{(1,1)}_{{\mathbf{B}}_{\ell +1}} = 
{ \big( {\alpha }_i(h_j) \big) }_{\! \stackrel{2 \le i \le \ell +1 }{\mbox{$\scriptscriptstyle 2 \le j \le \ell +1$}}}$ for 
the root \vspace{-3pt} \\ 
\rule{0pt}{12pt} \hspace{-4pt} \rule{.25in}{0in} system $\mathcal{R} = {\mathbf{B}}_{\ell }$ generated by $\Pi $.} & \\
& Set $h^{\ast } = 2\sum^{\ell }_{i=1}h_i + h_{\ell +1} = 2e_1$. Then $\mathcal{R} = {\mathrm{R}}^0$.& \\
& ${\mathrm{R}}^{-} = \big\{ -{\varepsilon }_1; \, - ({\varepsilon }_1 \pm {\varepsilon }_i), \, 2 \le i \le \ell +1 \big\} $ & \\
& nilradical roots: & \\
& $\widehat{\mathcal{R}} = \left\{ \begin{array}{rl}
{\widehat{\alpha}}_0 & = -{\varepsilon}_1 | \mathfrak{h} = 0|\mathfrak{h} \\
{\widehat{\alpha }}_{j-1} & = (-{\varepsilon }_1 + {\varepsilon }_j) |\mathfrak{h}, \, 2 \le j \le \ell +1 \\
{\widehat{\alpha }}_{\ell + j-1} & = (-{\varepsilon }_1 - {\varepsilon }_j) | \mathfrak{h} , \, 2 \le j \le \ell +1 
\end{array} \right. $ & \\
& additive relations in $\widehat{\mathcal{R}}$: none & \\
& nilradical structure: $Z_{2\ell +1}$ & \\ 
\\ 
2.2 ${\widetilde{\mathbf{B}}}^2_{\ell +1}$ &{\large ${\mathbf{B}}^{(\ell +1)}_{\ell +1}$}:\hspace{-.6in} \begin{tabular}{lll}
\begin{picture}(0,15)(-40,-3)
\begin{drawjoin}
\jput(0,5){\circle*{2.25}}
\jput(20,5){\circle*{2.25}} 
\jput(40,5){\circle*{2.25}}
\jput(60,5){\circle*{.0001}}
\end{drawjoin}
\begin{dottedjoin}{3}
\jput(60,5){\circle*{.0001}}
\jput(80,5){\circle*{.0001}}
\end{dottedjoin}
\begin{drawjoin}
\jput(80,5){\circle{.0001}}
\jput(100,5){\circle*{2.25}}
\end{drawjoin}
\put(100,4.5){\line(1,0){20}} 
\put(100,5.5){\line(1,0){20}}
\drawline(112,6.5)(109,5)
\drawline(112,3.5)(109,5)
\jput(120,5){\circle{2.25}}
\put(120,5){\makebox(0,0){\large $\times$}}
\put(0,0){\makebox(0,-1){\small $1$}}
\put(20,0){\makebox(0,-1){\small $2$}} 
\put(40,0){\makebox(0,-1){\small $3$}}
\put(100,0){\makebox(0,-1){\small $\ell$}}
\put(125,0){\makebox(0,-1){\small $\ell +1$}}
\end{picture}  
\end{tabular} & \\
& positive simple roots: $\Pi = \{ {\alpha }_i , \, 1 \le i \le \ell \} $ & \\
& Cartan subalgebra: $\mathfrak{h} = {\spann }_{\C} \{ h_i , \, 1 \le i \le \ell \} $ & \\
& \parbox[t]{4.5in}{Cartan matrix $C^{(\ell +1,\ell +1)}_{{\mathbf{B}}_{\ell +1}} = 
{ \big( {\alpha }_i(h_j) \big) }_{\! \stackrel{1 \le i \le \ell}{\mbox{$\scriptscriptstyle 1 \le j \le \ell $}}}$ for the root 
\vspace{-2pt} \\
\rule{.25in}{0in} system $\mathcal{R} = {\mathbf{A}}_{\ell }$ generated by $\Pi $.} & \\
& Set $h^{\ast } = 2\sum^{\ell }_{i=1}ih_i + (\ell +1)h_{\ell +1} = 2\sum^{\ell +1}_{i=1} e_i$. 
Then $\mathcal{R} = {\mathrm{R}}^0$. & \\
& ${\mathrm{R}}^{-} = \big\{ -{\varepsilon }_j, \, 1 \le j \le \ell +1 ; 
\, - ({\varepsilon }_i + {\varepsilon }_j), \, 1 \le i < j \le \ell +1 \big\} $ & \\
& nilradical roots: & \\
& $\widehat{\mathcal{R}} = \left\{ \begin{array}{rl}
{\widehat{\zeta }}_{i < j} & = -({\varepsilon }_i + {\varepsilon }_j) | \mathfrak{h} , \, 1 \le i < j \le \ell +1 \\
{\widehat{\alpha}}_j & = -{\varepsilon }_j |\mathfrak{h}, \, 1 \le j \le \ell +1 
\end{array} \right. $ & \\
& additive relations in $\widehat{\mathcal{R}}$: ${\widehat{\alpha }}_i+{\widehat{\alpha }}_j = 
{\widehat{\zeta }}_{i<j}, \, \, 1 \le i < j \le \ell +1$ & \\
& nilradical structure: $\sum_{1 \le i < j \le \ell +1} + {\mathrm{h}}^{{\widehat{\zeta }}_{i<j}}_3$ 
\end{tabular}
\newpage
\begin{tabular}{lll}
\vspace{-1.25in} &  & \\ 
3.& {\large ${\mathbf{C}}_{\ell +1} \, (\ell \ge 2)$:} \hspace{-.65in} \begin{tabular}{lll}
\begin{picture}(0,15)(-40,-3)
\begin{drawjoin}
\jput(0,5){\circle*{2.25}}
\jput(20,5){\circle*{2.25}} 
\jput(40,5){\circle*{2.25}}
\jput(60,5){\circle*{.0001}}
\end{drawjoin}
\begin{dottedjoin}{3}
\jput(60,5){\circle*{.0001}}
\jput(80,5){\circle*{.0001}}
\end{dottedjoin}
\begin{drawjoin}
\jput(80,5){\circle{.0001}}
\jput(100,5){\circle*{2.25}}
\end{drawjoin}
\put(100,4.5){\line(1,0){20}} 
\put(100,5.5){\line(1,0){20}}
\drawline(109,6.5)(112,5)
\drawline(109,3.5)(112,5)
\jput(120,5){\circle{2.25}}
\put(0,0){\makebox(0,-1){\small $1$}}
\put(20,0){\makebox(0,-1){\small $2$}} 
\put(40,0){\makebox(0,-1){\small $3$}}
\put(100,0){\makebox(0,-1){\small $\ell$}}
\put(125,0){\makebox(0,-1){\small $\ell +1$}}
\end{picture}  
\end{tabular} & \\ 
& positive simple roots: & \\
& \rule{.25in}{0pt} $\underline{\Pi } = \{ {\alpha }_i ={\varepsilon }_i-{\varepsilon }_{i+1}, \, 1 \le i \le \ell ; 
\, {\alpha }_{\ell +1} = 2{\varepsilon }_{\ell +1} \} $ & \\
& Cartan subalgebra: & \\
&\rule{.25in}{0pt}  $\underline{\mathfrak{h}} = {\spann }_{\C } \big\{ h_i = e_i-e_{i+1}, \, 1 \le i \le \ell ; \, h_{\ell +1} = e_{\ell +1} \big\} $ & \\ 
\rule{0pt}{30pt}& \hspace{.15in} $C_{{\mathbf{C}}_{\ell +1}} = 
{\big( {\alpha }_i(h_j) \big) }_{\! \stackrel{1 \le i \le \ell +1}{\mbox{$\scriptscriptstyle 1 \le j \le \ell +1$}}} =
\mbox{\tiny $\left( \begin{array}{rrrrrr}
2          & -1          & 0             & 0                         & \cdots           & 0        \\
-1         &  2          &       -1      & 0                        & \cdots          & 0         \\
\vdots   & \ddots   & \ddots      & \ddots                 & \cdots           & \vdots  \\
            & \cdots    &-1              & 2                        & -1                   & 0           \\
  0         & \cdots    &  0            & -1                       & 2                   & -1  \\
  0        &  \cdots            &  0            & 0                       & -2                   & 2
\end{array} \right) $} $ &  \\
\rule{0pt}{12pt}& \parbox[t]{4.5in}{Cartan matrix for the root system 
$\underline{\mathcal{R}} = {\mathbf{C}}_{\ell +1}$ \\ 
\rule{.25in}{0in} generated by $\underline{\Pi }$.} & \\
& positive roots: ${\underline{\mathcal{R}}}^{+} = $ & \\
& \hspace{.15in} $ \left\{ \begin{array}{rl}
2{\varepsilon }_i = & \hspace{-5pt} 2\sum^{\ell +1}_{k=i}{\alpha }_k, \, 1 \le i \le \ell +1  \\
\rule{0pt}{12pt} {\varepsilon }_i - {\varepsilon }_j = & \hspace{-5pt} \sum^{j-1}_{k=i} {\alpha }_k, \, 1 \le i < j \le \ell +1 \\
\rule{0pt}{12pt} {\varepsilon }_i + {\varepsilon }_j = & \hspace{-5pt} \sum^{j-1}_{k=i} {\alpha }_k + 
2 \sum^{\ell +1}_{k=j}{\alpha }_k, \, 1 \le i < j \le \ell +1
\end{array} \right. $ & \\
\\ 
3.1 ${\widetilde{\mathbf{C}}}^1_{\ell +1}$ &{\large ${\mathbf{C}}^{(1)}_{\ell +1}$:}\hspace{-.5in}\begin{tabular}{l}
\begin{picture}(0,15)(-40,-3)
\begin{drawjoin}
\jput(0,5){\circle*{2.25}}
\jput(20,5){\circle*{2.25}} 
\jput(40,5){\circle*{2.25}}
\jput(60,5){\circle*{.0001}}
\end{drawjoin}
\begin{dottedjoin}{3}
\jput(60,5){\circle*{.0001}}
\jput(80,5){\circle*{.0001}}
\end{dottedjoin}
\begin{drawjoin}
\jput(80,5){\circle{.0001}}
\jput(100,5){\circle*{2.25}}
\end{drawjoin}
\put(100,4.5){\line(1,0){20}} 
\put(100,5.5){\line(1,0){20}}
\drawline(109,6.5)(112,5)
\drawline(109,3.5)(112,5)
\jput(120,5){\circle{2.25}}
\put(0,5){\makebox(0,0){\large $\times$}}
\put(0,0){\makebox(0,-1){\small $1$}}
\put(20,0){\makebox(0,-1){\small $2$}} 
\put(40,0){\makebox(0,-1){\small $3$}}
\put(100,0){\makebox(0,-1){\small $\ell$}}
\put(125,0){\makebox(0,-1){\small $\ell +1$}}
\end{picture}  
\end{tabular} & \\
& positive simple roots: $\Pi = \{ {\alpha }_i , \, 2 \le i \le \ell +1 \} $ & \\
& Cartan subalgebra: $\mathfrak{h} = {\spann }_{\C} \{ h_i , \, 2 \le i \le \ell +1 \} $ & \\
& \parbox[t]{4.5in}{Cartan matrix $C^{(1,1)}_{{\mathbf{C}}_{\ell +1}} = 
{ \big( {\alpha }_i(h_j) \big) }_{\! \stackrel{2 \le i \le \ell +1 }{\mbox{$\scriptscriptstyle 2 \le j \le \ell +1$}}}$ for the root system \vspace{-3pt} \\ 
\rule{.25in}{0in} $\mathcal{R} = {\mathbf{C}}_{\ell }$ generated by $\Pi $.} & \\
\rule{0pt}{12pt} & Set $h^{\ast } = \sum^{\ell +1}_{i=1}h_i  = e_1$. Then $\mathcal{R} = {\mathrm{R}}^0$. & \\
& ${\mathrm{R}}^{-} = \big\{ -2{\varepsilon }_1, \, ; 
\, - {\varepsilon }_1 \pm {\varepsilon }_j, \, 2 \le j \le \ell +1 \big\} $ & \\
& nilradical roots: & \\
& $\widehat{\mathcal{R}} = \left\{ \begin{array}{rl}
\widehat{\zeta } & = -2{\varepsilon }_1|\mathfrak{h} \\
{\widehat{\alpha }}_j & = -({\varepsilon }_1 -{\varepsilon }_{j+1}) |\mathfrak{h}, \, 1 \le j \le \ell \\
{\widehat{\alpha }}^{\, j} & = -({\varepsilon }_1 + {\varepsilon }_{j+1}) | \mathfrak{h} , \, 1 \le j \le \ell 
\end{array} \right. $ & \\
& additive relations in $\widehat{\mathcal{R}}$: ${\widehat{\alpha }}_k +{\widehat{\alpha }}^{\, k} = \widehat{\zeta }, 
\, \, 1 \le k \le \ell $ & \\
& nilradical structure: ${\mathrm{h}}^{\widehat{\zeta }}_{2\ell +1}$ & \\
\\ 
3.2 ${\widetilde{\mathbf{C}}}^2_{\ell +1}$ &{\large ${\mathbf{C}}^{(\ell +1)}_{\ell +1}$:}\hspace{-.5in}\begin{tabular}{l}
\begin{picture}(0,15)(-40,-3)
\begin{drawjoin}
\jput(0,5){\circle*{2.25}}
\jput(20,5){\circle*{2.25}} 
\jput(40,5){\circle*{2.25}}
\jput(60,5){\circle*{.0001}}
\end{drawjoin}
\begin{dottedjoin}{3}
\jput(60,5){\circle*{.0001}}
\jput(80,5){\circle*{.0001}}
\end{dottedjoin}
\begin{drawjoin}
\jput(80,5){\circle{.0001}}
\jput(100,5){\circle*{2.25}}
\end{drawjoin}
\put(100,4.5){\line(1,0){20}} 
\put(100,5.5){\line(1,0){20}}
\drawline(109,6.5)(112,5)
\drawline(109,3.5)(112,5)
\jput(120,5){\circle{2.25}}
\put(120,5){\makebox(0,0){\large $\times$}}
\put(0,0){\makebox(0,-1){\small $1$}}
\put(20,0){\makebox(0,-1){\small $2$}} 
\put(40,0){\makebox(0,-1){\small $3$}}
\put(100,0){\makebox(0,-1){\small $\ell$}}
\put(125,0){\makebox(0,-1){\small $\ell +1$}}
\end{picture}  
\end{tabular} & \\
& positive simple roots: $\Pi = \{ {\alpha }_i , \, 1 \le i \le \ell \} $ & \\
& Cartan subalgebra: $\mathfrak{h} = {\spann }_{\C} \{ h_i , \, 1 \le i \le \ell \} $ & \\
& \parbox[t]{4.5in}{Cartan matrix $C^{(\ell +1, \ell +1)}_{{\mathbf{C}}_{\ell +1}} = 
{ \big( {\alpha }_i(h_j) \big) }_{\! \stackrel{1\le i \le \ell }{\mbox{$\scriptscriptstyle 1 \le j \le \ell $}}} $  
for the \\ 
\rule{.25in}{0in} root system $\mathcal{R} = {\mathbf{A}}_{\ell }$ generated by 
$\Pi $.} & \\
\rule{0pt}{12pt} & Set $h^{\ast } = \sum^{\ell +1}_{i=1}ih_i  = \sum^{\ell +1}_{i=1}e_i$. Then $\mathcal{R} = {\mathrm{R}}^0$. & \\
& ${\mathrm{R}}^{-} = \big\{ -2{\varepsilon }_i, \, 1 \le i \le \ell +1;\, 
 - ({\varepsilon }_i + {\varepsilon }_j), \, 1 \le i < j \le \ell +1  \big\} $ & \\
& nilradical roots: & \\
& $\widehat{\mathcal{R}} = \left\{ \begin{array}{rl}
{\widehat{\alpha }}_j & = -2{\varepsilon }_j|\mathfrak{h}, \, 1 \le j \le \ell +1 \\
{\widehat{\alpha }}_{i<j} & = -({\varepsilon }_i  + {\varepsilon }_j) |\mathfrak{h}, \, 1 \le i < j \le \ell +1  
\end{array} \right. $ & \\
& additive relations in $\widehat{\mathcal{R}}$: none & \\
& nilradical structure: $Z_{\mbox{$\scriptscriptstyle \frac{1}{2}$}(\ell +2)(\ell +1)}$ & 
\end{tabular}
\newpage
\begin{tabular}{lll}
\vspace{-1.5in} &  & \\ 
4.& {\large ${\mathbf{D}}_{\ell +1}\, (\ell \ge 3)$:} \hspace{-.65in} \begin{tabular}{l}
\begin{picture}(0,15)(-40,-3)
\begin{drawjoin}
\jput(0,5){\circle*{2}}
\jput(20,5){\circle*{2}} 
\jput(40,5){\circle*{2}}
\jput(60,5){\circle*{.0001}}
\end{drawjoin}
\begin{drawjoin}
\jput(80,5){\circle*{.0001}}
\jput(100,5){\circle*{2}}
\end{drawjoin}
\begin{dottedjoin}{3}
\jput(60,5){\circle*{.0001}}
\jput(80,5){\circle*{.0001}}
\end{dottedjoin}
\begin{drawjoin}[30]
\jput(100,5){\circle*{2}}
\jput(110,10){\circle*{2}}
\end{drawjoin}
\begin{drawjoin}[30]
\jput(100,5){\circle*{2}}
\jput(110,0){\circle*{2}}
\end{drawjoin}
\put(0,-1){\makebox(0,0){\small $1$}}
\put(20,-1){\makebox(0,0){\small $2$}} 
\put(40,-1){\makebox(0,0){\small $3$}}
\put(95,-2){\makebox(0,0){\footnotesize $\ell -1$}}
\put(110,16){\makebox(0,0){\small $\ell $}}
\put(115,-8){\makebox(0,0){\small $\ell +1$}}
\end{picture}
\end{tabular} & \\ 
\rule{0pt}{14pt} & positive simple roots: & \\
& $\underline{\Pi } = \{ {\alpha }_i = {\varepsilon }_i-{\varepsilon }_{i+1}, 
\, 1 \le i \le \ell ; \, {\alpha }_{\ell +1} = {\varepsilon }_{\ell } + {\varepsilon }_{\ell +1} \} $ & \\
& Cartan subalgebra: & \\
& $\underline{\mathfrak{h}} = {\spann }_{\C} \{ h_i = e_i -e_{i+1}, \, 1 \le i \le \ell ; \, 
h_{\ell +1} = e_{\ell } + e_{\ell +1} \} $ \\
\rule{0pt}{30pt}& \hspace{.15in} $C_{{\mathbf{D}}_{\ell +1}} = 
{\big( {\alpha }_i(h_j) \big) }_{\! \stackrel{1 \le i \le \ell +1}{\mbox{$\scriptscriptstyle 1 \le j \le \ell +1$}}} =
\mbox{\tiny $\left( \begin{array}{rrrrrr}
2          & -1          & 0             & 0                         & \cdots           & 0        \\
-1         &  2          &   -1         & 0                         & \cdots           & 0         \\
0           & -1         & 2            & -1                         & \cdots          & 0          \\
\vdots  & \ddots   & \ddots      & \ddots                 & \cdots           & \vdots  \\      
            & \cdots    &-1            & 2                        & -1                 & -1           \\
0          & \cdots    &  0            & -1                       & 2                    & 0  \\
0          &   \cdots  &  0            & -1                       & 0                  & 2
\end{array} \right) $} $ &  \\
\rule{0pt}{12pt}& \parbox[t]{4.5in}{Cartan matrix for the root system 
$\underline{\mathcal{R}} = {\mathbf{D}}_{\ell +1}$ generated by $\underline{\Pi }$} & \\
& positive roots: ${\underline{\mathcal{R}}}^{+} =$ & \\
& $ \left\{ \begin{array}{rl}
{\varepsilon }_i - {\varepsilon }_j = & \hspace{-5pt} \sum^{j-1}_{k=i} {\alpha }_k, \, 1 \le i < j \le \ell +1 \\
\rule{0pt}{12pt} {\varepsilon }_i + {\varepsilon }_j = & \hspace{-5pt} \left\{ 
\begin{array}{l}
\sum^{\ell}_{k=j} {\alpha }_k +  \sum^{\ell -1}_{k=i}{\alpha }_k+{\alpha }_{\ell +1}, \, \, \mbox{if $1 \le i < j \le \ell $} \\
\rule{0pt}{12pt} \hspace{-4pt} \sum^{\ell -1}_{k=i} {\alpha }_k +{\alpha }_{\ell +1}, \, \,  \mbox{if $1 \le i \le \ell -1 
\, \, \& \, \, j = \ell +1$} \\
{\alpha }_{\ell +1}, \, \,  \mbox{if $i = \ell \, \, \& \, \, j = \ell +1$}
\end{array} \right. 
\end{array} \right. $ & \\
\\ 
4.1 ${\widetilde{\mathbf{D}}}^1_{\ell+1}$ & {\large ${\mathbf{D}}^{(1)}_{\ell +1}$:}\hspace{-.5in}\begin{tabular}{l}
\begin{picture}(0,15)(-40,-3)
\begin{drawjoin}
\jput(0,5){\circle*{2}}
\jput(20,5){\circle*{2}} 
\jput(40,5){\circle*{2}}
\jput(60,5){\circle*{.0001}}
\end{drawjoin}
\begin{drawjoin}
\jput(80,5){\circle*{.0001}}
\jput(100,5){\circle*{2}}
\end{drawjoin}
\begin{dottedjoin}{3}
\jput(60,5){\circle*{.0001}}
\jput(80,5){\circle*{.0001}}
\end{dottedjoin}
\begin{drawjoin}[30]
\jput(100,5){\circle*{2}}
\jput(110,10){\circle*{2}}
\end{drawjoin}
\begin{drawjoin}[30]
\jput(100,5){\circle*{2}}
\jput(110,0){\circle*{2}}
\end{drawjoin}
\put(0,5){\makebox(0,0){\large $\times$}}
\put(0,-1){\makebox(0,0){\small $1$}}
\put(20,-1){\makebox(0,0){\small $2$}} 
\put(40,-1){\makebox(0,0){\small $3$}}
\put(95,-2){\makebox(0,0){\footnotesize $\ell -1$}}
\put(110,16){\makebox(0,0){\small $\ell $}}
\put(115,-8){\makebox(0,0){\small $\ell +1$}}
\end{picture}
\end{tabular} & \\
\rule{0pt}{19pt} & positive simple roots: $\Pi  = {\{ {\alpha }_i \} }^{\ell +1}_{i=2}$ & \\
& Cartan subalgebra: $\mathfrak{h} = {\spann }_{\C} {\{ h_i \} }^{\ell +1}_{i=2} $ \\
& \parbox[t]{4.5in}{Cartan matrix 
$C^{(1,1)}_{{\mathbf{D}}_{\ell +1}}=
{ \big( {\alpha }_i(h_j) \big) }_{\! \stackrel{2 \le i \le \ell +1}{\mbox{$\scriptscriptstyle 2 \le j \le \ell +1$} } }$ 
for the root system \\  
\rule{.25in}{0in} $\mathcal{R} ={\mathbf{D}}_{\ell }$ generated by $\Pi $.} & \\ 
\rule{0pt}{12pt}& Set $h^{\ast } = 2\sum^{\ell -1}_{i=1}h_i +h_{\ell } +h_{\ell +1} = 2e_1$. Then $\mathcal{R} = {\mathrm{R}}^0$. & \\
& ${\mathrm{R}}^{-} = \big\{ - ({\varepsilon }_1 \pm  {\varepsilon }_j), \, 2 \le j \le \ell +1  \big\} $ & \\
& nilradical roots: & \\
& $\widehat{\mathcal{R}} = \left\{ \begin{array}{rl} 
{\widehat{\alpha }}_{j-1} = & \hspace{-5pt}  (-{\varepsilon }_1  + {\varepsilon }_j) |\mathfrak{h}, \, 2 \le j \le \ell +1 \\
{\widehat{\alpha }}_{\ell + j-1} = & \hspace{-5pt}  -({\varepsilon }_1  + {\varepsilon }_j) |\mathfrak{h}, 
\, 2 \le j \le \ell +1
\end{array} \right. $ & \\
& additive relations in $\widehat{\mathcal{R}}$: none & \\
& nilradical structure: $Z_{2\ell }$ & \\
\\ 
4.2 ${\widetilde{\mathbf{D}}}^2_{\ell +1}$ &{\large ${\mathbf{D}}^{(\ell )}_{\ell +1}$:}\hspace{-.5in}\begin{tabular}{l}
\begin{picture}(0,15)(-40,-3)
\begin{drawjoin}
\jput(0,5){\circle*{2}}
\jput(20,5){\circle*{2}} 
\jput(40,5){\circle*{2}}
\jput(60,5){\circle*{.0001}}
\end{drawjoin}
\begin{drawjoin}
\jput(80,5){\circle*{.0001}}
\jput(100,5){\circle*{2}}
\end{drawjoin}
\begin{dottedjoin}{3}
\jput(60,5){\circle*{.0001}}
\jput(80,5){\circle*{.0001}}
\end{dottedjoin}
\begin{drawjoin}[30]
\jput(100,5){\circle*{2}}
\jput(110,10){\circle*{2}}
\end{drawjoin}
\begin{drawjoin}[30]
\jput(100,5){\circle*{2}}
\jput(110,0){\circle*{2}}
\end{drawjoin}
\put(110,10){\makebox(0,0){\large $\times$}}
\put(0,-1){\makebox(0,0){\small $1$}}
\put(20,-1){\makebox(0,0){\small $2$}} 
\put(40,-1){\makebox(0,0){\small $3$}}
\put(95,-2){\makebox(0,0){\footnotesize $\ell -1$}}
\put(110,16){\makebox(0,0){\small $\ell $}}
\put(115,-8){\makebox(0,0){\small $\ell +1$}}
\end{picture}
\end{tabular} & \\
\rule{0pt}{19pt} & positive simple roots: $\Pi  = \{ {\alpha }_i, \, 1 \le i \le \ell -1; {\alpha }_{\ell +1} \} $ & \\
& Cartan subalgebra: $\mathfrak{h} = {\spann }_{\C} \{ h_i, \, 1 \le i \le \ell -1; \, h_{\ell +1}  \} $ \\
& \parbox[t]{4.5in}{Cartan matrix $C^{(\ell , \ell )}_{{\mathbf{D}}_{\ell +1}}=
{ \big( {\alpha }_i(h_j) \big) }_{\! \stackrel{1 \le i \le \ell -1 \& i = \ell +1}{\mbox{$\scriptscriptstyle 1 \le j \le \ell -1 
\& j = \ell +1$} } }$ for the root \\  
\rule{0pt}{9pt}\hspace{-4pt}\rule{.25in}{0in} system $\mathcal{R} = {\mathbf{A}}_{\ell }$ generated by $\Pi $.} & \\ 
& Set & \\
& $h^{\ast } = 2\sum^{\ell -1}_{i=1}i h_i + (\ell +1) h_{\ell } + (\ell -1)h_{\ell +1} = 
2\sum^{\ell }_{i=1} e_i - 2e_{\ell +1}$. & \\
& Then $\mathcal{R} = {\mathrm{R}}^0$. & \\
& ${\mathrm{R}}^{-} = \big\{ - {\varepsilon }_i+  {\varepsilon }_{\ell +1} , \, 1 \le i \le \ell ; 
\, -({\varepsilon }_i + {\varepsilon }_j), \, 1 \le i < j \le \ell   \big\} $ & \\
& nilradical roots: $\widehat{\mathcal{R}} = \left\{ \begin{array}{rl} 
{\widehat{\alpha }}_{j} = & \hspace{-5pt}  (-{\varepsilon }_j  + {\varepsilon }_{\ell +1} ) |\mathfrak{h}, 
\, 1 \le  j \le \ell  \\
{\widehat{\alpha }}_{i<j } = & \hspace{-5pt}  -({\varepsilon }_i  + {\varepsilon }_j) |\mathfrak{h}, 
\, 1 \le i < j \le \ell 
\end{array} \right. $ & \\
& additive relations in $\widehat{\mathcal{R}}$: none & \\
& nilradical structure: $Z_{\mbox{$\scriptscriptstyle \frac{1}{2}$} \ell (\ell +1)}$ & 
\end{tabular}
\newpage 
\begin{tabular}{lll}
4.3 ${\widetilde{\mathbf{D}}}^3_{\ell +1}$ &{\large ${\mathbf{D}}^{(\ell +1)}_{\ell +1}$:}\hspace{-.5in}\begin{tabular}{l}
\begin{picture}(0,15)(-40,-3)
\begin{drawjoin}
\jput(0,5){\circle*{2}}
\jput(20,5){\circle*{2}} 
\jput(40,5){\circle*{2}}
\jput(60,5){\circle*{.0001}}
\end{drawjoin}
\begin{drawjoin}
\jput(80,5){\circle*{.0001}}
\jput(100,5){\circle*{2}}
\end{drawjoin}
\begin{dottedjoin}{3}
\jput(60,5){\circle*{.0001}}
\jput(80,5){\circle*{.0001}}
\end{dottedjoin}
\begin{drawjoin}[30]
\jput(100,5){\circle*{2}}
\jput(110,10){\circle*{2}}
\end{drawjoin}
\begin{drawjoin}[30]
\jput(100,5){\circle*{2}}
\jput(110,0){\circle*{2}}
\end{drawjoin}
\put(110,0){\makebox(0,0){\large $\times$}}
\put(0,-1){\makebox(0,0){\small $1$}}
\put(20,-1){\makebox(0,0){\small $2$}} 
\put(40,-1){\makebox(0,0){\small $3$}}
\put(95,-2){\makebox(0,0){\footnotesize $\ell -1$}}
\put(110,16){\makebox(0,0){\small $\ell $}}
\put(115,-8){\makebox(0,0){\small $\ell +1$}}
\end{picture}
\end{tabular} & \\
\rule{0pt}{19pt} & positive simple roots: $\Pi  = \{ {\alpha }_i, \, 1 \le i \le \ell  \} $ & \\
& Cartan subalgebra: $\mathfrak{h} = {\spann }_{\C} \{ h_i, \, 1 \le i \le \ell  \} $ \\
& \parbox[t]{4.5in}{Cartan matrix 
$C^{(\ell +1, \ell +1 )}_{{\mathbf{D}}_{\ell +1}}=
{ \big( {\alpha }_i(h_j) \big) }_{\stackrel{\! \! 1 \le i \le \ell }{\mbox{$\scriptscriptstyle 1 \le j \le \ell $ } } }$  
for the root \\
\rule{.25in}{0in} system $\mathcal{R} = {\mathbf{A}}_{\ell }$ generated by $\Pi $.} & \\ 
& Set & \\
& $h^{\ast } = 2\sum^{\ell -1}_{i=1}i h_i + (\ell -1) h_{\ell } + (\ell +1) h_{\ell +1} = 
2\sum^{\ell +1}_{i=1} e_i $. & \\
& Then $\mathcal{R} = {\mathrm{R}}^0$. & \\
& ${\mathrm{R}}^{-} = \big\{ - ( {\varepsilon }_i +  {\varepsilon }_j ), \, 1 \le i < j \le \ell +1 \big\} $  & \\
& nilradical roots: & \\
& $\widehat{\mathcal{R}} = \big\{ {\widehat{\alpha }}_{i<j } = 
  -({\varepsilon }_i  + {\varepsilon }_j) |\mathfrak{h}, \, 1 \le i < j \le \ell +1   $ & \\
& additive relations in $\widehat{\mathcal{R}}$: none & \\
& nilradical structure: $Z_{\mbox{$\scriptscriptstyle \frac{1}{2}$} \ell (\ell +1)}$ & 
\\ 
5. & {\large ${\mathbf{G}}_2$}: \hspace{-.65in} \begin{tabular}{l}
\begin{picture}(0,15)(-40,-3)
\begin{drawjoin}
\jput(0,5){\circle*{2}}
\jput(20,5){\circle*{2}} 
\end{drawjoin}
\drawline(0,6)(20,6)
\drawline(0,5)(20,5)
\drawline(0,4)(20,4)
\drawline(10,7.5)(8,5)
\drawline(10,2.5)(8,5)
\put(0,-1){\makebox(0,0){\small $1$}}
\put(20,-1){\makebox(0,0){\small $2$}} 
\end{picture}
\end{tabular} & \\ 
& positive simple roots: $\underline{\Pi } = \big\{ {\alpha }_1 = {\varepsilon }_1 -{\varepsilon }_2, \, 
{\alpha }_2 = -2{\varepsilon }_1 +{\varepsilon}_2 +{\varepsilon }_3 \big\} $ & \\
& Cartan subalgebra: $\underline{\mathfrak{h}} = \big\{ h_1 = e_1-e_2, \, 
h_2 = -2e_1-e_2-e_3 \big\} $ & \\
\rule{0pt}{17pt} & \hspace{.75in} $C_{{\mathbf{G}}_2} = 
{\big( {\alpha }_i(h_j) \big) }_{\stackrel{\! \! 1 \le i \le 2 }{\mbox{$\scriptscriptstyle 1 \le j \le 2$ } }} = 
\mbox{\footnotesize $\left( \! \! \begin{array}{rr}  2 & -1 \\ -3 & 2 \end{array} \! \! \right) $}$ & \\
& \parbox[t]{4.5in}{Cartan matrix for root system $\mathcal{R} = {\mathbf{G}}_2$ 
generated by $\underline{\Pi }$.} & \\
& positive roots: ${\underline{\mathcal{R}}}^{+} = $ & \\ 
& $\left\{ \begin{array}{l}
{\alpha }_1 =  {\varepsilon }_1 - {\varepsilon }_2, \, \, \, 
{\alpha }_2 = -2{\varepsilon }_1 +{\varepsilon }_2 +{\varepsilon }_3, \\
{\alpha }_1 + {\alpha }_2 =   -{\varepsilon }_1 +{\varepsilon }_3, \, \, \, 
2{\alpha }_1+{\alpha }_2 =  -{\varepsilon }_2 +{\varepsilon }_3, \\
3{\alpha }_1+{\alpha }_2 =   {\varepsilon }_1 -2{\varepsilon }_2 + {\varepsilon }_3, \, \, \, 
3{\alpha }_1 + 2{\alpha }_2 =   -{\varepsilon }_1 -{\varepsilon }_2 +2 {\varepsilon }_3 
\end{array} \right. $ \\
\\ 
5.1 ${\widetilde{\mathbf{G}}}^1_2$ & {\large ${\mathbf{G}}^{(1)}_2$}: \hspace{-.65in} \begin{tabular}{l}
\begin{picture}(0,15)(-40,-3)
\begin{drawjoin}
\jput(0,5){\circle*{2}}
\jput(20,5){\circle*{2}} 
\end{drawjoin}
\drawline(0,6)(20,6)
\drawline(0,5)(20,5)
\drawline(0,4)(20,4)
\drawline(10,7.5)(8,5)
\drawline(10,2.5)(8,5)
\put(0,5){\makebox(0,0){\large $\times$}}
\put(0,-1){\makebox(0,0){\small $1$}}
\put(20,-1){\makebox(0,0){\small $2$}} 
\end{picture}
\end{tabular} & \\ 
& positive simple roots: $\Pi  = \{ {\alpha }_2 \} $ & \\
& Cartan subalgebra: $\mathfrak{h} = {\spann }_{\C} \{ h_2 \} $ & \\
& \parbox[t]{4.5in}{Cartan matrix $C^{(1,1)}_{{\mathbf{G}}_2} = 
{\big( {\alpha }_i(h_j) \big) }_{\stackrel{\! \! i = 2 }{\mbox{$\scriptscriptstyle j = 2$ } }}$ 
for root system $\mathcal{R} = {\mathbf{A}}_1$ \vspace{-3pt} \\ 
\rule{.25in}{0in} generated by $\Pi $.} & \\
& Set $h^{\ast } = 2h_1+3h_2 = -4e_1-5e_2 - 3e_3$. Then $\mathcal{R} = {\mathrm{R}}^0$. & \\ 
& ${\mathrm{R}}^{-} = \big\{ {\varepsilon }_1 - {\varepsilon }_3, \, {\varepsilon }_2 -{\varepsilon }_3, \,  
{\varepsilon }_2-{\varepsilon }_1, \, 
-{\varepsilon }_1+2{\varepsilon }_2 -{\varepsilon }_3, \, {\varepsilon }_1+{\varepsilon }_2 -2{\varepsilon }_3
 \big\} $ & \\
& nilradical roots: & \\
& $\widehat{\mathcal{R}} = \left\{ \begin{array}{rl}
\widehat{\zeta} = & \hspace{-5pt} ({\varepsilon }_2 - {\varepsilon }_3)|\mathfrak{h} = 0|\mathfrak{h} \\
{\widehat{\alpha }}_1 = & \hspace{-5pt} ({\varepsilon }_1 - {\varepsilon }_3)|\mathfrak{h} = 
({\varepsilon }_1 + {\varepsilon }_2 - 2{\varepsilon }_3)|\mathfrak{h}  \\  
{\widehat{\alpha }}^{\, 1} = & \hspace{-5pt}  ({\varepsilon }_2 - {\varepsilon }_1)|\mathfrak{h} = 
(-{\varepsilon }_1+2{\varepsilon }_2 -{\varepsilon }_3)|\mathfrak{h}
\end{array} \right. $ & \\
& additive relations in $\widehat{\mathcal{R}}$: ${\widehat{\alpha }}_1 +{\widehat{\alpha }}^1 = \widehat{\zeta }$ & \\
& nilradical structure: ${\mathrm{h}}^{\widehat{\zeta}}_3$
\end{tabular}
\newpage
\begin{tabular}{lll} 
5.2 ${\widetilde{\mathbf{G}}}^2_2$ & {\large ${\mathbf{G}}^{(2)}_2$}: \hspace{-.65in} \begin{tabular}{l}
\begin{picture}(0,15)(-40,-3)
\begin{drawjoin}
\jput(0,5){\circle*{2}}
\jput(20,5){\circle*{2}} 
\end{drawjoin}
\drawline(0,6)(20,6)
\drawline(0,5)(20,5)
\drawline(0,4)(20,4)
\drawline(10,7.5)(8,5)
\drawline(10,2.5)(8,5)
\put(20,5){\makebox(0,0){\large $\times$}}
\put(0,-1){\makebox(0,0){\small $1$}}
\put(20,-1){\makebox(0,0){\small $2$}} 
\end{picture}
\end{tabular} & \\ 
& positive simple roots: $\Pi  = \{ {\alpha }_1 \} $ & \\
& Cartan subalgebra: $\mathfrak{h} = {\spann }_{\C }\{ h_1 \} $ & \\
& \parbox[t]{4.5in}{Cartan matrix $C^{(2,2)}_{{\mathbf{G}}_2} = 
{\big( {\alpha }_i(h_j) \big) }_{\stackrel{\! \! i = 1 }{\mbox{$\scriptscriptstyle j = 1$ } }}$ 
for root system $\mathcal{R} = {\mathbf{A}}_1$ \\ 
\rule{.25in}{0in} generated by $\Pi $.} & \\
& Set $h^{\ast } = h_1+2h_2 = -3e_1-3e_2 - 2e_3$. Then $\mathcal{R} = {\mathrm{R}}^0$. & \\ 
& ${\mathrm{R}}^{-} = \big\{ {\varepsilon }_1 - {\varepsilon }_3, \, {\varepsilon }_2 -{\varepsilon }_3, \, 
2{\varepsilon }_1-{\varepsilon }_2 -{\varepsilon}_3, \, 
-{\varepsilon }_1+2{\varepsilon }_2 -{\varepsilon }_3, \, {\varepsilon }_1+{\varepsilon }_2 -2{\varepsilon }_3 \big\} $ & \\
& nilradical roots: & \\
& $\widehat{\mathcal{R}} = \left\{ \begin{array}{rl}
\widehat{\zeta} = & \hspace{-5pt}({\varepsilon }_1+{\varepsilon }_2 -2{\varepsilon }_3)|\mathfrak{h}, \, \, \, 
{\widehat{\alpha }}_1 =(-{\varepsilon }_1 +2{\varepsilon }_2 - {\varepsilon }_3)|\mathfrak{h}  \\ 
{\widehat{\alpha }}^{\, 1} = & \hspace{-5pt} (2{\varepsilon }_1-{\varepsilon}_2 - {\varepsilon }_3)|\mathfrak{h},  \, \, \, 
{\widehat{\alpha }}_2 = ({\varepsilon }_2 -{\varepsilon }_3)|\mathfrak{h}, \, \, \, 
{\widehat{\alpha }}^{\, 2} = ({\varepsilon }_1 - {\varepsilon }_3)|\mathfrak{h}
\end{array} \right. $ & \\
& additive relations in $\widehat{\mathcal{R}}$: ${\widehat{\alpha }}_k +{\widehat{\alpha }}^{\, k} = \widehat{\zeta }, 
\, \, k =1,2 $ & \\
& nilradical structure: ${\mathrm{h}}^{\widehat{\zeta}}_5$
\\ 
6. & {\large ${\mathbf{F}}_4$}: \hspace{-.65in} \begin{tabular}{l} 
\begin{picture}(0,15)(-40,-3)
\begin{drawjoin}
\jput(0,5){\circle*{2}}
\jput(20,5){\circle*{2}} 
\end{drawjoin}
\put(20,5){\circle*{2}}
\put(20,5.5){\line(1,0){20}}
\put(20,4.5){\line(1,0){20}} 
\begin{drawjoin}
\jput(40,5){\circle*{2}}
\jput(60,5){\circle*{2}}
\end{drawjoin}
\drawline(28,6.5)(30,5)
\drawline(28,3.5)(30,5)
\put(0,0){\makebox(0,0){\small $1$}}
\put(20,0){\makebox(0,0){\small $2$}} 
\put(40,0){\makebox(0,0){\small $3$}}
\put(60,0){\makebox(0,0){\small $4$}}
\end{picture}
\end{tabular}  \\ 
& positive simple roots: & \\
& $\underline{\Pi } = \{ {\alpha }_1 = {\varepsilon }_2-{\varepsilon }_3, \, 
{\alpha }_2 = {\varepsilon }_3-{\varepsilon }_4, \, {\alpha }_3 = {\varepsilon }_4, \,  
{\alpha }_4= \onehalf ({\varepsilon }_1 -  {\varepsilon }_2 - {\varepsilon }_3 - {\varepsilon }_4)  \} $ & \\
& Cartan subalgebra: & \\
& $\underline{\mathfrak{h}} = {\spann }_{\C } \{ h_1 = e_2-e_3, \, h_2 = e_3-e_4, \, h_3 = 2e_4, \, 
h_4 = e_1-e_2-e_3-e_4 \} $ & \\
\rule{0pt}{17pt} & \hspace{.75in} $C_{{\mathbf{F}}_4} = 
{\big( {\alpha }_i(h_j) \big) }_{\stackrel{\! \! 1 \le i \le 4}{\mbox{$\scriptscriptstyle 1 \le j \le 4 $ } }} = 
\mbox{\tiny $ \left( \begin{array}{rrrr} 
2  & -1 &  0 & 0 \\
-1&   2 & -2 & 0 \\
0 & -1 & 2 & -1 \\
0 & 0 & -1 & 2 \end{array} \right) $} $ & \\
& \parbox[t]{4.5in}{Cartan matrix for root system $\underline{\mathcal{R}} = {\mathbf{F}}_4$ 
generated by $\underline{\Pi } $.} & \\
& positive roots: & \\
& ${\underline{\mathcal{R}}}^{+} = \big\{ {\varepsilon }_i, \, 1 \le i \le 4; \, 
{\varepsilon }_i \pm {\varepsilon }_j, \, 1 \le i < j \le 4; \, 
\onehalf ({\varepsilon }_1 \pm  {\varepsilon }_2 \pm {\varepsilon }_3 \pm {\varepsilon }_4) \big\} $ & \\
\\ 
6.1 ${\widetilde{\mathbf{F}}}^1_4$ & {\large ${\mathbf{F}}^{(4)}_4$}: \hspace{-.65in} \begin{tabular}{l} 
\begin{picture}(0,15)(-40,-3)
\begin{drawjoin}
\jput(0,5){\circle*{2}}
\jput(20,5){\circle*{2}} 
\end{drawjoin}
\put(20,5){\circle*{2}}
\put(20,5.5){\line(1,0){20}}
\put(20,4.5){\line(1,0){20}} 
\begin{drawjoin}
\jput(40,5){\circle*{2}}
\jput(60,5){\circle*{2}}
\end{drawjoin}
\drawline(28,6.5)(30,5)
\drawline(28,3.5)(30,5)
\put(60,5){\makebox(0,0){\large $\times$}}
\put(0,0){\makebox(0,0){\small $1$}}
\put(20,0){\makebox(0,0){\small $2$}} 
\put(40,0){\makebox(0,0){\small $3$}}
\put(60,0){\makebox(0,0){\small $4$}}
\end{picture}
\end{tabular}  \\ 
& positive simple roots: & \\
& $\Pi  = \{ {\alpha }_i, \, 1 \le i \le 3  \} $ & \\
& Cartan subalgebra: $\mathfrak{h} = {\spann}_{\C} \{ h_i, \, 1 \le i \le 3 \} $ & \\
& \parbox[t]{4.5in}{Cartan matrix $C^{(4,4)}_{{\mathbf{F}}_4} = 
{\big( {\alpha }_i(h_j) \big) }_{\stackrel{\! \! 1 \le i \le 3}{\mbox{$\scriptscriptstyle 1 \le j \le 3 $ } }}$ 
for root system $\mathcal{R} = {\mathbf{B}}_3$ \\ 
\rule{.25in}{0in} generated by $\Pi $.} & \\
& Set $h^{\ast } = 2h_1+3h_2+2h_3+h_4 = 2e_1$. Then $\mathcal{R} = {\mathrm{R}}^0$. & \\
& ${\mathrm{R}}^{-} = \big\{ -{\varepsilon}_1; \, -{\varepsilon }_1 \pm {\varepsilon }_j, \, 2 \le j \le 4; \, 
\onehalf (-{\varepsilon }_1\pm {\varepsilon }_2\pm {\varepsilon }_3 \pm {\varepsilon }_4) \big\} $ & \\ 
& nilradical roots: $\widehat{\mathcal{R}} = $ & \\
& $ \left\{ \begin{array}{l} 
{\widehat{\zeta}}_1 =  -{\varepsilon}_1|\mathfrak{h} = 0|\mathfrak{h} \\
{\widehat{\zeta}}_j = (- {\varepsilon }_1 +{\varepsilon }_j) |\mathfrak{h}, \, 2 \le j \le 4 \\
{\widehat{\zeta}}_{3 +j} = (- {\varepsilon }_1 - {\varepsilon }_j) |\mathfrak{h}, \, 2 \le j \le 4  \\ 
{\widehat{\alpha }}_1 =  \onehalf (-{\varepsilon }_1+ {\varepsilon }_2 +{\varepsilon }_3 +{\varepsilon }_4)|\mathfrak{h}, 
\, \, \, 
{\widehat{\alpha }}_2 = \onehalf (-{\varepsilon }_1 - {\varepsilon }_2 +{\varepsilon }_3 +{\varepsilon }_4 )|\mathfrak{h} 
\end{array} \right. $ &  
\end{tabular}
\newpage
\begin{tabular}{lll}
\vspace{-1.25in} & & \\
& $\left\{ \begin{array}{l} 
{\widehat{\alpha }}_3 = \onehalf (-{\varepsilon }_1+ {\varepsilon }_2 -{\varepsilon }_3 +
{\varepsilon }_4)|\mathfrak{h}, \, \, \, 
{\widehat{\alpha }}_{4} =  \onehalf (-{\varepsilon }_1 + {\varepsilon }_2 +{\varepsilon }_3 -
{\varepsilon }_4)|\mathfrak{h} \\
\rule{0pt}{12pt} {\widehat{\alpha }}_{5} =  \onehalf (-{\varepsilon }_1- {\varepsilon }_2 -{\varepsilon }_3 +
{\varepsilon }_4)|\mathfrak{h}, \, \, \, 
{\widehat{\alpha }}_{6} = \onehalf (-{\varepsilon }_1 - {\varepsilon }_2 +{\varepsilon }_3 -
{\varepsilon }_4)|\mathfrak{h} \\
\rule{0pt}{12pt} {\widehat{\alpha }}_{7} =   \onehalf (-{\varepsilon }_1+ {\varepsilon }_2 - {\varepsilon }_3 -
{\varepsilon }_4)|\mathfrak{h}, \, \, \, 
{\widehat{\alpha }}_{8} =  \onehalf (-{\varepsilon }_1 - {\varepsilon }_2 -{\varepsilon }_3 -
{\varepsilon }_4)|\mathfrak{h} 
\end{array} \right. $ & \\
& additive relations in $\widehat{\mathcal{R}}$: & \\
&\mbox{\footnotesize $ \left\{ \begin{array}{l}
{\widehat{\alpha }}_1+{\widehat{\alpha }}_8 = {\widehat{\zeta }}_1 \\
{\widehat{\alpha }}_2+{\widehat{\alpha }}_7 = {\widehat{\zeta }}_1 \\
{\widehat{\alpha }}_3+{\widehat{\alpha }}_6 = {\widehat{\zeta }}_1 \\
{\widehat{\alpha }}_4+{\widehat{\alpha }}_5 = {\widehat{\zeta }}_1 
\end{array} \right. $} $ \, \, \, \begin{array}{l}
\mbox{\footnotesize $\left\{ \begin{array}{l}
{\widehat{\alpha }}_1+{\widehat{\alpha }}_7 = {\widehat{\zeta }}_2 \\
{\widehat{\alpha }}_3+{\widehat{\alpha }}_4 = {\widehat{\zeta }}_2 \end{array} \right. $} \, \, \,  
\mbox{\footnotesize $\left\{ \begin{array}{l}
{\widehat{\alpha }}_1+{\widehat{\alpha }}_6 = {\widehat{\zeta }}_3 \\
{\widehat{\alpha }}_2+{\widehat{\alpha }}_4 = {\widehat{\zeta }}_3 \end{array} \right. $} \, \, \,  
\mbox{\footnotesize $\left\{ \begin{array}{l}
{\widehat{\alpha }}_1+{\widehat{\alpha }}_5 = {\widehat{\zeta }}_4 \\
{\widehat{\alpha }}_2+{\widehat{\alpha }}_3 = {\widehat{\zeta }}_4 \end{array} \right. $} \\
\rule{0pt}{20pt} \mbox{\footnotesize $\left\{ \begin{array}{l}
{\widehat{\alpha }}_2+{\widehat{\alpha }}_8 = {\widehat{\zeta }}_5 \\
{\widehat{\alpha }}_5+{\widehat{\alpha }}_6 = {\widehat{\zeta }}_5 \end{array} \right. $} \, \, \,  
\mbox{\footnotesize $\left\{ \begin{array}{l}
{\widehat{\alpha }}_3+{\widehat{\alpha }}_8 = {\widehat{\zeta }}_6 \\
{\widehat{\alpha }}_5+{\widehat{\alpha }}_7 = {\widehat{\zeta }}_6 \end{array} \right. $} \, \, \,  
\mbox{\footnotesize $\left\{ \begin{array}{l}
{\widehat{\alpha }}_4+{\widehat{\alpha }}_8 = {\widehat{\zeta }}_7 \\
{\widehat{\alpha }}_6+{\widehat{\alpha }}_7 = {\widehat{\zeta }}_7 \end{array} \right. $} 
\end{array}  $ & \\
& nilradical structure: ${\mathrm{h}}^{{\widehat{\zeta}}_1}_9 + 
\sum^7_{k=2}+{\mathrm{h}}^{{\widehat{\zeta}}_k}_5$ & \\ 
\\ 
6.2 ${\widetilde{\mathbf{F}}}^2_4$ & {\large ${\mathbf{F}}^{(1)}_4$}: \hspace{-.65in} \begin{tabular}{l} 
\begin{picture}(0,15)(-40,-3)
\begin{drawjoin}
\jput(0,5){\circle*{2}}
\jput(20,5){\circle*{2}} 
\end{drawjoin}
\put(20,5){\circle*{2}}
\put(20,5.5){\line(1,0){20}}
\put(20,4.5){\line(1,0){20}} 
\begin{drawjoin}
\jput(40,5){\circle*{2}}
\jput(60,5){\circle*{2}}
\end{drawjoin}
\drawline(28,6.5)(30,5)
\drawline(28,3.5)(30,5)
\put(0,5){\makebox(0,0){\large $\times$}}
\put(0,0){\makebox(0,0){\small $1$}}
\put(20,0){\makebox(0,0){\small $2$}} 
\put(40,0){\makebox(0,0){\small $3$}}
\put(60,0){\makebox(0,0){\small $4$}}
\end{picture}
\end{tabular}  \\ 
& positive simple roots: & \\
& $\Pi  = {\{ {\alpha }_i  \} }^4_{i=2} $ & \\
& Cartan subalgebra: $\mathfrak{h} = {\spann }_{\C} \{ h_i, \, 2 \le i \le 4 \} $ & \\
& \parbox[t]{4.5in}{Cartan matrix $C^{(1,1)}_{{\mathbf{F}}_4} = 
{\big( {\alpha }_i(h_j) \big) }_{\stackrel{\! \! 2 \le i \le 4}{\mbox{$\scriptscriptstyle 2 \le j \le 4 $ } }}$ 
for root system $\mathcal{R} = {\mathbf{C}}_3$ \vspace{-3pt} \\ 
\rule{.25in}{0in} generated by $\Pi $.} & \\
& Set $h^{\ast } = 2h_1+3h_2+2h_3+h_4 = e_1+e_2$. Then $\mathcal{R} = {\mathrm{R}}^0$. & \\
& ${\mathrm{R}}^{-} = \left\{ -{\varepsilon}_1; \, - {\varepsilon}_2; \, 
- ({\varepsilon }_1+{\varepsilon }_2); \,  -{\varepsilon }_1 \pm {\varepsilon }_j, \, j = 3, 4; 
\, -{\varepsilon }_2 \pm {\varepsilon }_j, \, j = 3, 4; \right. $ & \\ 
& $\left. \onehalf (-{\varepsilon }_1 - {\varepsilon }_2\pm {\varepsilon }_3 \pm {\varepsilon }_4) \right\} $ & \\
& nilradical roots: $\widehat{\mathcal{R}} = $& \\
& $\left\{ \begin{array}{l} 
\widehat{\zeta } = -({\varepsilon }_1 +{\varepsilon }_2 )| \mathfrak{h},  \\
{\widehat{\alpha }}_1 = -{\varepsilon }_1|\mathfrak{h}, \\ 
{\widehat{\alpha }}^{\, 1} = -{\varepsilon }_2|\mathfrak{h} \\ 
{\widehat{\alpha }}_{j-1} = (-{\varepsilon }_1 + {\varepsilon }_j ) |\mathfrak{h}, \, j =3, 4 \\
{\widehat{\alpha }}^{\, j-1} = (-{\varepsilon }_2 -{\varepsilon }_j )| \mathfrak{h}, \, j=3,4 \\
{\widehat{\alpha }}_{j+1} = (-{\varepsilon }_1 - {\varepsilon }_j ) |\mathfrak{h}, \, j =3, 4 \\
{\widehat{\alpha }}^{\, j+1} = (-{\varepsilon }_2 +{\varepsilon }_j )| \mathfrak{h}, \, j=3,4 \\
{\widehat{\alpha }}_{6}=  \onehalf (-{\varepsilon }_1- {\varepsilon }_2 +{\varepsilon }_3 +
{\varepsilon }_4)|\mathfrak{h}, \, \, \, 
{\widehat{\alpha }}_{7} = \onehalf (-{\varepsilon }_1 - {\varepsilon }_2 -{\varepsilon }_3 +
{\varepsilon }_4)|\mathfrak{h} \\
\rule{0pt}{12pt} {\widehat{\alpha }}^{\, 6} =  \onehalf (-{\varepsilon }_1 - {\varepsilon }_2 -{\varepsilon }_3 -
{\varepsilon }_4)|\mathfrak{h} , \, \, \, {\widehat{\alpha }}^{\, 7} =  \onehalf (-{\varepsilon }_1 - {\varepsilon }_2 +{\varepsilon }_3 - {\varepsilon }_4)|\mathfrak{h}  
\end{array} \right. $ & \\
& additive relations in $\widehat{\mathcal{R}}$: ${\widehat{\alpha }}_k+{\widehat{\alpha }}^k = \widehat{\zeta }, 
\, \, 1 \le k \le 7$ & \\
& nilradical structure: ${\mathrm{h}}^{\widehat{\zeta }}_{15} $ & \\
&& \\ 
7. & {\large ${\mathbf{E}}_6$}: \hspace{-.65in} \begin{tabular}{l}
\begin{picture}(0,15)(-40,-3)
\begin{drawjoin}
\jput(0,5){\circle*{2}}
\jput(20,5){\circle*{2}} 
\jput(40,5){\circle*{2}}
\jput(60,5){\circle*{2}}
\jput(80,5){\circle*{2}}
\end{drawjoin}
\put(40,5){\line(0,1){15}}
\put(40,20){\circle*{2}}
\put(0,-1){\makebox(0,0){\small $1$}}
\put(20,-1){\makebox(0,0){\small $3$}} 
\put(40,-1){\makebox(0,0){\small $4$}}
\put(44,20){\makebox(0,0){\small $2$}}
\put(60,-1){\makebox(0,0){\small $5$}}
\put(80,-1){\makebox(0,0){\small $6$}}
\end{picture}
\end{tabular} & \\
& simple positive roots: $\underline{\Pi } = \big\{ {\alpha }_1 = \onehalf ( {\varepsilon }_1 -\sum^5_{i=2}{\varepsilon }_i +{\varepsilon }_6), \, {\alpha }_2 = {\varepsilon }_2+{\varepsilon }_1, $ & \\
\rule{0pt}{12pt} & $\hspace{.5in} {\alpha }_i = {\varepsilon }_{i-1}-{\varepsilon }_{i-2}, \, 3 \le i \le 6 \big\} $ & \\
& Cartan subalgebra: $\underline{\mathfrak{h}} = {\spann }_{\C } \big\{  
h_1 = \onehalf (e_1-\sum^5_{i=2}e_i +3e_6), $ &\\ 
&\rule{.25in}{0in} $ \, h_2 = e_2+e_1, \, h_i = e_{i-1}-e_{i-2}, \, 3 \le i \le 6 \big\} $ & \\
\rule{0pt}{25pt} &\hspace{.25in} $C_{{\mathbf{E}}_6} =
{\big( {\alpha }_i(h_j) \big) }_{\stackrel{\! \! 1 \le i \le 6}{\mbox{$\scriptscriptstyle 1 \le j \le 6$ } }} = 
\mbox{\tiny $\left( \begin{array}{rrrrrr} 
2 & 0 & -1 & 0 & 0 & 0   \\
0 & 2 & 0 & -1 & 0 & 0 \\
-1 & 0 & 2 & -1 & 0 & 0 \\
0 & -1 & -1 & 2 & -1 & 0 \\
0 & 0 & 0 & -1 & 2 & -1 \\
0 & 0 & 0 & 0 & -1 & 2 \end{array} \right) $} $ & \\ 
& Cartan matrix for root system $\underline{\mathcal{R}}$ generated by $\underline{\Pi}$ & \\ 
\rule{0pt}{12pt}& positive roots: ${\underline{\mathcal{R}}}^{+} = \big\{ {\varepsilon }_i \pm {\varepsilon }_j, \, 
1 \le i < j \le 5; \, \onehalf (\sum^5_{i=1}(-1)^{k(i)}{\varepsilon }_i +{\varepsilon }_6), $ & \\ 
& \rule{0pt}{12pt} $\hspace{.5in} k(i) = 0$ or $1$ and $\sum^5_{k=1}k(i) \in \{ 0,2,4 \} \big\} $ & 
\end{tabular}
\newpage
\begin{tabular}{lll}
\vspace{-.25in} & & \\
7.1 ${\widetilde{\mathbf{E}}}^1_6$ &\rule{0pt}{25pt}\hspace{-.06in} {\large ${\mathbf{E}}^{(1)}_6$}: \hspace{-.65in} \begin{tabular}{l}
\begin{picture}(0,15)(-40,-3)
\begin{drawjoin}
\jput(0,5){\circle*{2}}
\jput(20,5){\circle*{2}} 
\jput(40,5){\circle*{2}}
\jput(60,5){\circle*{2}}
\jput(80,5){\circle*{2}}
\end{drawjoin}
\put(40,5){\line(0,1){15}}
\put(40,20){\circle*{2}}
\put(0,5){\makebox(0,0){\large $\times$}} 
\put(0,-1){\makebox(0,0){\small $1$}}
\put(20,-1){\makebox(0,0){\small $3$}} 
\put(40,-1){\makebox(0,0){\small $4$}}
\put(44,20){\makebox(0,0){\small $2$}}
\put(60,-1){\makebox(0,0){\small $5$}}
\put(80,-1){\makebox(0,0){\small $6$}}
\end{picture}
\end{tabular} & \\
& simple positive roots: $\Pi = \{ {\alpha }_i, \, 2 \le i \le 6 \} $ & \\
& Cartan subalgebra: $\mathfrak{h} = \{ h_i, 2 \le i \le 6 \} $ & \\
&\parbox[t]{4.5in}{Cartan matrix $C^{(1,1)}_{{\mathbf{E}}_6} =
{\big( {\alpha }_i(h_j) \big) }_{\stackrel{\! \! 2 \le i \le 6}{\mbox{$\scriptscriptstyle 2 \le j \le 6$ } }} $ 
for root system $\mathcal{R} = {\mathbf{D}}_5$ \vspace{-5pt} \\ 
\rule{0pt}{9pt} \hspace{-4pt} \rule{.25in}{0in} generated by $\Pi $.} & \\
\rule{0pt}{10pt} & Set $h^{\ast } = 4h_1+3h_2+5h_3+6h_4+4h_5 + 2h_6 = 6e_6$. Then $\mathcal{R} ={\mathrm{R}}^0$ & \\
& ${\mathrm{R}}^{-} = \{ \onehalf \big( \sum^5_{i=1} (-1)^{k(i)} {\varepsilon }_i -{\varepsilon }_6 \big) $, \, 
\mbox{where $k(i) =0$ or $1$ and} & \\
\rule{0pt}{11pt} &\hspace{.75in} \mbox{$\sum^5_{i=1}k(i) \in \{ 1,3,5 \} $} & \\
& nilradical roots: $\widehat{\mathcal{R}} = $ & \\ 
& $\left\{ \begin{array}{rl}
{\widehat{\alpha }}_j & = \big( \onehalf \sum^5_{i=1} (-1)^{k(i)} {\varepsilon }_i - {\varepsilon }_6 \big) 
|\mathfrak{h}, \\
& \hspace{.25in} \mbox{where $k(j) =1, \, k(i) =0$ if $i \notin \{ j \} , \, j \in \{ 1, \ldots , 5 \} $} \\
{\widehat{\alpha }}_{j\ell m} & = \big( \onehalf \sum^5_{i=1} (-1)^{k(i)} {\varepsilon }_i -{\varepsilon }_6 \big) |\mathfrak{h}, \\
& \hspace{.25in} \mbox{where $k(j) = k(\ell ) = k(m) = 1, \, k(i) =0$ if $i \notin \{ j, \ell ,m \} $} \\
& \hspace{.25in}\mbox{$\, j,\ell , m \in \{ 1, \ldots , 5 \} $ distinct } \\
{\widehat{\alpha }}_{-} & = - \onehalf \big( \sum^6_{i=1}{\varepsilon }_i \big) |\mathfrak{h}
\end{array} \right. $ & \\
& additive relations in $\widehat{\mathcal{R}}$: none & \\
& nilradical structure: $Z_{16}$ & \\
\\
7.2 ${\widetilde{\mathbf{E}}}^2_6$ &\rule{0pt}{25pt}\hspace{-.06in} {\large ${\mathbf{E}}^{(6)}_6$}: \hspace{-.65in} \begin{tabular}{l}
\begin{picture}(0,15)(-40,-3)
\begin{drawjoin}
\jput(0,5){\circle*{2}}
\jput(20,5){\circle*{2}} 
\jput(40,5){\circle*{2}}
\jput(60,5){\circle*{2}}
\jput(80,5){\circle*{2}}
\end{drawjoin}
\put(40,5){\line(0,1){15}}
\put(40,20){\circle*{2}}
\put(80,5){\makebox(0,0){\large $\times$}} 
\put(0,-1){\makebox(0,0){\small $1$}}
\put(20,-1){\makebox(0,0){\small $3$}} 
\put(40,-1){\makebox(0,0){\small $4$}}
\put(44,20){\makebox(0,0){\small $2$}}
\put(60,-1){\makebox(0,0){\small $5$}}
\put(80,-1){\makebox(0,0){\small $6$}}
\end{picture}
\end{tabular} & \\
& simple positive roots: $\Pi = \{ {\alpha }_i, \,  1 \le i \le 5  \} $ & \\
& Cartan subalgebra: $\mathfrak{h} = \{ h_i, 1 \le i \le 5 \} $ & \\
& \parbox[t]{4.5in}{Cartan matrix: $C^{(6,6)}_{{\mathbf{E}}_6} =
{\big( {\alpha }_i(h_j) \big) }_{\stackrel{\! \! 1 \le i \le 5}{\mbox{$\scriptscriptstyle 1 \le j \le 5$ } }}$ for root system 
$\mathcal{R} = {\mathbf{D}}_5$ \vspace{-5pt} \\ 
\rule{0pt}{10pt} \hspace{-4pt} \rule{.25in}{0in} generated by $\Pi $.} & \\
\rule{0pt}{10pt} & Set $h^{\ast } = 2h_1+3h_2+4h_3 + 6h_4+5h_5 + 4h_6 = 3(e_5+e_6)$. & \\
& Then $\mathcal{R} = {\mathrm{R}}^0$. & \\
& ${\mathrm{R}}^{-}  = \big\{ \pm {\varepsilon }_i - {\varepsilon }_5, \, 1 \le i \le 4; \, 
\onehalf \big( \sum^4_{i=1} (-1)^{k(i)}{\varepsilon }_i - {\varepsilon }_5 - {\varepsilon }_6 \big) |\mathfrak{h}$; & \\
& \hspace{5pt} \mbox{$k(i) = 0$ or $1$; $\sum^4_{i=1}k(i) \in \{ 0,2,4 \} $} & \\ 
\rule{0pt}{12pt} & nilradical roots: $\widehat{\mathcal{R}}  = $ & \\
& $ \left\{ \begin{array}{rl}
{\widehat{\alpha }}_j = & \hspace{-5pt} ({\varepsilon }_j - {\varepsilon }_5)|\mathfrak{h}, \, 
2 \le j \le 4 \\
{\widehat{\alpha }}_{j+4}= & \hspace{-5pt} (-{\varepsilon }_j - {\varepsilon }_5)|\mathfrak{h}, \, 
2 \le j \le 4 \\
{\widehat{\alpha }}_{+} = & \hspace{-5pt}\onehalf \big( \sum^4_{i=1}{\varepsilon }_i - {\varepsilon }_5
- {\varepsilon }_6 \big) |\mathfrak{h} \\
\rule{0pt}{12pt} {\widehat{\alpha }}_{j \ell } = & \hspace{-5pt} \onehalf \big( \sum^4_{i=1} (-1)^{k(i)}{\varepsilon }_i - {\varepsilon }_5
- {\varepsilon }_6 \big) |\mathfrak{h}, \, \mbox{$k(i) =0$ or $1$};  \\
& \hspace{.25in}\parbox[t]{3in}{$k(j) =k(\ell ) = 1$ or $1$; $k(i) =0$ if $i \notin \{ j , \ell \}$; \\
$j, \ell \in \{ 1,2,3,4 \}$ distinct;} \\
{\widehat{\alpha }}_{-} = & \hspace{-5pt} -\onehalf \big( \sum^6_{i=1}{\varepsilon }_i \big) |\mathfrak{h} 
\end{array} \right. $ & \\
& additive structure in $\widehat{\mathcal{R}}$: none & \\
& nilradical structure: $Z_{16}$ &  
\end{tabular} 
\newpage
\begin{tabular}{lll}
\vspace{-.5in} & & \\ 
7.3 ${\widetilde{\mathbf{E}}}^3_6$ &\rule{0pt}{25pt}\hspace{-.06in} {\large ${\mathbf{E}}^{(2)}_6$}: \hspace{-.65in} \begin{tabular}{l}
\begin{picture}(0,15)(-40,-3)
\begin{drawjoin}
\jput(0,5){\circle*{2}}
\jput(20,5){\circle*{2}} 
\jput(40,5){\circle*{2}}
\jput(60,5){\circle*{2}}
\jput(80,5){\circle*{2}}
\end{drawjoin}
\put(40,5){\line(0,1){15}}
\put(40,20){\circle*{2}}
\put(40,20){\makebox(0,0){\large $\times$}} 
\put(0,-1){\makebox(0,0){\small $1$}}
\put(20,-1){\makebox(0,0){\small $3$}} 
\put(40,-1){\makebox(0,0){\small $4$}}
\put(46,20){\makebox(0,0){\small $2$}}
\put(60,-1){\makebox(0,0){\small $5$}}
\put(80,-1){\makebox(0,0){\small $6$}}
\end{picture}
\end{tabular} & \\
& simple positive roots: $\Pi = \big\{ {\alpha }_i, \, i \in I = \{ 1,3,4,5,6\} \big\} $ & \\
& Cartan subalgebra: $\mathfrak{h} = {\spann }_{\C} \big\{ h_i, \, i \in I \big\} $ & \\
& \parbox[t]{4,5in}{Cartan matrix: $C^{(2,2)}_{{\mathbf{E}}_6} = {\big( {\alpha }_i(h_j) \big)}_{i,j \in I}$ for root system \\  
\rule{.25in}{0in} $\mathcal{R} = {\mathbf{A}}_5$ generated by $\Pi $.} & \\ 
& Set & \\
& $h^{\ast } = 2h_1+4h_2+4h_3+6h_4+4h_5+2h_6 = \sum^5_{i=1}e_i + 3e_6$. & \\
& Then $\mathcal{R} = {\mathrm{R}}^0$. & \\
& ${\mathrm{R}}^{-} = \left\{ \begin{array}{l}
-({\varepsilon }_i +{\varepsilon }_j), \, 1 \le i < j \le 5; \, 
\onehalf \big( \sum^5_{i=1}(-1)^{k(i)}{\varepsilon }_i -{\varepsilon }_6 \big), \\
\parbox[t]{2.75in}{\hspace{.25in} $k(i) =0$ or $1$; $\sum^5_{i=1}k(i) \in \{ 3,5 \}$.} 
\end{array} \right. $ & \\
& nilradical roots: $\widehat{\mathcal{R}} = $ &\\
&\hspace{.25in} $\left\{ \begin{array}{rl}
\widehat{\zeta} = & \hspace{-5pt} -\onehalf \big( \sum^6_{i=1}{\varepsilon }_i \big) |\mathfrak{h} \\ 
{\widehat{\alpha}}_{i<j} = & \hspace{-5pt} -({\varepsilon }_i +{\varepsilon }_j)|\mathfrak{h}, \, 1 \le i < j \le 5; \\ 
{\widehat{\alpha }}_{j \ell m} = &\hspace{-5pt}  \onehalf \big( \sum^5_{i=1}(-1)^{k(i)}{\varepsilon }_i -{\varepsilon }_6 \big) |\mathfrak{h}  \\
&\hspace{.15in} \parbox[t]{2.75in}{$k(j) =k(\ell ) = k(m) =1$; $k(i) =0$ if \\
$i \notin \{ j , \ell , m \}$; $j, \ell, m \in \{ 1, \ldots , 5 \}$ distinct}  
\end{array} \right. $ & \\ 
& additive relations in $\widehat{\mathcal{R}}$: ${\widehat{\alpha }}_{i<j} + {\widehat{\alpha }}_{k\ell m} = 
\widehat{\zeta }$, \, \, if and only if & \\
& \rule{.25in}{0in} $\{ i,j,k, \ell , m \} = \{ 1, \ldots , 5 \} $  & \\
& nilradical structure: ${\mathrm{h}}^{\widehat{\zeta }}_{21}$ \\
\\
8. & {\large ${\mathbf{E}}_7$}: \hspace{-.65in} \begin{tabular}{l} 
\begin{picture}(0,15)(-40,-3)
\begin{drawjoin}
\jput(0,5){\circle*{2}}
\jput(20,5){\circle*{2}} 
\jput(40,5){\circle*{2}}
\jput(60,5){\circle*{2}}
\jput(80,5){\circle*{2}}
\jput(100,5){\circle*{2}}
\end{drawjoin}
\put(40,5){\line(0,1){15}}
\put(40,20){\circle*{2}}
\put(0,-1){\makebox(0,0){\small $1$}}
\put(20,-1){\makebox(0,0){\small $3$}} 
\put(40,-1){\makebox(0,0){\small $4$}}
\put(44,20){\makebox(0,0){\small $2$}}
\put(60,-1){\makebox(0,0){\small $5$}}
\put(80,-1){\makebox(0,0){\small $6$}}
\put(100,-1){\makebox(0,0){\small $7$}}
\end{picture}
\end{tabular} & \\ 
& simple positive roots: $\underline{\Pi } = \big\{ {\alpha }_1 = 
\onehalf (-\sum^6_{i=1}{\varepsilon }_i  + {\varepsilon }_7)$ ,& \\ 
& \hspace{.2in} ${\alpha }_2 = {\varepsilon }_2 -{\varepsilon }_1, \, 
{\alpha }_3 = {\varepsilon }_2 + {\varepsilon }_1, \, 
{\alpha }_i = {\varepsilon }_{i-1}-{\varepsilon }_{i-2}, \, 4 \le i \le 7 \big\} $ & \\
& Cartan subalgebra: $\underline{\mathfrak{h}} 
= {\spann }_{\C } \big\{ h_1 = \onehalf ( - \sum^6_{i=1}e_i + 2e_7 )$, & 
\end{tabular}
\begin{tabular}{lll}  
& \hspace{.25in}$h_2 = e_2-e_1, \,  h_3 = e_2+e_1, \, h_i = e_{i-1}-e_{i-2}, \, 4 \le i \le 7 \big\} $ & \\
\rule{0pt}{30pt} & \hspace{.25in} $C_{{\mathbf{E}}_7} = 
{\big( {\alpha }_i(h_j) \big) }_{\stackrel{\! \! 1 \le i \le 7}{\mbox{$\scriptscriptstyle 1 \le j \le 7$ } }} = 
\mbox{\tiny $\left( \begin{array}{rrrrrrr} 
2 & 0 & -1 & 0 & 0 & 0 & 0    \\
0 & 2 & 0 & -1 & 0 & 0 & 0 \\
-1 & 0 & 2 & -1 & 0 & 0 & 0  \\ 
0 & -1 & -1 & 2 & -1 & 0 & 0  \\
0 & 0 & 0 & -1 & 2 & -1 & 0 \\
0 & 0 & 0 & 0 & -1 & 2 & -1  \\ 
0 & 0 & 0 & 0 & 0 & -1 & 2 \end{array} \right) $} $ & \\ 
&\parbox[t]{4.5in}{Cartan matrix for root system $\underline{\mathcal{R}}$ generated by 
$\underline{\Pi }$.} & \\ 
& positive roots: & \\
& ${\underline{\mathcal{R}}}^{+} =  \{ {\varepsilon }_7; \, {\varepsilon }_i \pm {\varepsilon }_j, \, 
1 \le i < j \le  6; \onehalf \big( \sum^6_{i=1}(-1)^{k(i)}{\varepsilon }_i + {\varepsilon }_7 \big) $,  & \\
\rule{0pt}{12pt} &\hspace{.5in} \mbox{where $k(i) = 0$ or $1$ and $\sum^6_{i=1}k(i) \in \{ 0,2,4,6\} \} $} &  \\   
& & \\ 
8.1 ${\widetilde{\mathbf{E}}}^1_7$ & \rule{0pt}{25pt} {\large ${\mathbf{E}}^{(1)}_7$}: \hspace{-.65in} \begin{tabular}{l} 
\begin{picture}(0,15)(-40,-3)
\begin{drawjoin}
\jput(0,5){\circle*{2}}
\jput(20,5){\circle*{2}} 
\jput(40,5){\circle*{2}}
\jput(60,5){\circle*{2}}
\jput(80,5){\circle*{2}}
\jput(100,5){\circle*{2}}
\end{drawjoin}
\put(40,5){\line(0,1){15}}
\put(40,20){\circle*{2}}
\put(0,5){\makebox(0,0){\large $\times$}} 
\put(0,-1){\makebox(0,0){\small $1$}}
\put(20,-1){\makebox(0,0){\small $3$}} 
\put(40,-1){\makebox(0,0){\small $4$}}
\put(44,20){\makebox(0,0){\small $2$}}
\put(60,-1){\makebox(0,0){\small $5$}}
\put(80,-1){\makebox(0,0){\small $6$}}
\put(100,-1){\makebox(0,0){\small $7$}}
\end{picture}
\end{tabular} & \\
& simple positive roots: $\Pi = \{ {\alpha }_i, \, 2 \le i \le 7 \} $ & \\
& Cartan subalgebra: $\mathfrak{h} = {\spann }_{\C} \{ h_i, \, 2 \le i \le 7 \} $ & 
\end{tabular}
\begin{tabular}{lll}
\vspace{-.75in} & & \\
& \parbox[t]{4,5in}{Cartan matrix: 
$C^{(1,1)}_{{\mathbf{E}}_7} = 
{\big( {\alpha }_i(h_j) \big) }_{\stackrel{\! \! 2 \le i \le 7}{\mbox{$\scriptscriptstyle 2 \le j \le 7$ } }}$ for 
root system $\mathcal{R} = {\mathbf{D}}_6$ \vspace{-5pt} \\ 
\rule{0pt}{10pt}\hspace{-4pt} \rule{.25in}{0in} generated by $\Pi $} & \\ 
& Set $h^{\ast } = 2h_1+2h_2+3h_3+4h_4+3h_5+2h_6+h_7 = 2e_7 $. & \\
& Then $\mathcal{R} = {\mathrm{R}}^0 $ & \\
& ${\mathrm{R}}^{-} = \big\{ -{\varepsilon }_7;  \, 
\onehalf \big( \sum^6_{i=1} (-1)^{k(i)}{\varepsilon }_i - {\varepsilon }_7 \big)$ where $k(i)=0$ or $1$ & \\
&\rule{.5in}{0in} and $\sum^6_{i=1} k(i) \in \{ 0,2,4,6 \} \big\} $ & \\
& nilradical roots: $\widehat{\mathcal{R}} = $ & \\
& $ \left\{ \begin{array}{rl} 
\widehat{\zeta} & = - {\varepsilon }_7 |\mathfrak{h} = 0|\mathfrak{h} \\ 
{\widehat{\alpha }}_{-} & = -\onehalf \big( \sum^7_{i=1}{\varepsilon }_i \big) |\mathfrak{h} \\ 
\rule{0pt}{12pt} {\widehat{\alpha }}_{+} & = \onehalf \big(  \sum^6_{i=1}{\varepsilon }_i - {\varepsilon }_7 \big) |\mathfrak{h} \\ 
\rule{0pt}{12pt} {\widehat{\alpha }}_{j\ell } & = 
\onehalf \big(  \sum^6_{i=1}(-1)^{k(i)}{\varepsilon }_i - {\varepsilon }_7 \big) |\mathfrak{h}, \, 
\mbox{where $k(j) =k(\ell ) = 1$;} \\ 
&\hspace{.25in}\mbox{$k(i) =0$ if $j, \ell  \ne \{ i \}$ and $j, \ell  \in \{ 1, \ldots , 6 \} $ distinct } \\
\rule{0pt}{12pt} {\widehat{\alpha }}^{\, j\ell} & = \onehalf ( \sum^6_{i=1} (-1)^{k(i)}{\varepsilon }_i 
- {\varepsilon }_7) |\mathfrak{h}, \, \mbox{where $k(j) = k(\ell ) = 0$,}  \\
&\rule{.25in}{0in} \mbox{$k(i) =1$, if $i \notin \{ j, \ell\} $ and $j, \ell \in \{ 1, \ldots ,6 \}$ distinct} 
\end{array} \right. $ & \\
& additive relations in $\widehat{\mathcal{R}}$: ${\widehat{\alpha }}_{+} + {\widehat{\alpha }}_{-} = 
\widehat{\zeta }$ and ${\widehat{\alpha }}_{j\ell } + {\widehat{\alpha }}^{\, j\ell } = \widehat{\zeta }$ & \\
&\rule{.25in}{0in} $j, \ell \in \{ 1, \ldots , 6 \} $ distinct & \\
& nilradical structure: ${\mathrm{h}}^{\widehat{\zeta }}_{87}$ & \\ 
\vspace{.1in} && \\ 
\noindent 8.2 ${\widetilde{\mathbf{E}}}^2_7$ & {\large ${\mathbf{E}}^{(7)}_7$}: \hspace{-.65in} \begin{tabular}{l} 
\begin{picture}(0,15)(-40,-3)
\begin{drawjoin}
\jput(0,5){\circle*{2}}
\jput(20,5){\circle*{2}} 
\jput(40,5){\circle*{2}}
\jput(60,5){\circle*{2}}
\jput(80,5){\circle*{2}}
\jput(100,5){\circle*{2}}
\end{drawjoin}
\put(40,5){\line(0,1){15}}
\put(40,20){\circle*{2}}
\put(100,5){\makebox(0,0){\large $\times$}} 
\put(0,-2){\makebox(0,0){\small $1$}}
\put(20,-2){\makebox(0,0){\small $3$}} 
\put(40,-2){\makebox(0,0){\small $4$}}
\put(44,20){\makebox(0,0){\small $2$}}
\put(60,-2){\makebox(0,0){\small $5$}}
\put(80,-2){\makebox(0,0){\small $6$}}
\put(100,-2){\makebox(0,0){\small $7$}}
\end{picture}
\end{tabular} & \\ 
& simple positive roots: $\Pi = {\{ {\alpha }_i \} }^6_{i=1}$ & \\
& Cartan subalgebra: $\mathfrak{h} = {\spann }_{\C} \{ h_i, \, 1 \le i \le 6 \} $ & \\
& \parbox[t]{4,5in}{Cartan matrix: 
$C^{(7,7)}_{{\mathbf{E}}_7} = 
{\big( {\alpha }_i(h_j) \big) }_{\stackrel{\! \! 1 \le i \le 6}{\mbox{$\scriptscriptstyle 1 \le j \le 6$ } }}$ for 
root system $\mathcal{R} = {\mathbf{E}}_6$ \vspace{-5pt} \\ 
\rule{0pt}{12pt} \rule{.25in}{0in} generated by $\Pi $.} & \\ 
& Set $h^{\ast } = 2h_1+3h_2+4h_3+6h_4+5h_5+4h_6+3h_7 = 2(e_6+e_7) $. & \\
& Then $\mathcal{R} = {\mathrm{R}}^0 $ & \\
& $\begin{array}{l}
{\mathrm{R}}^{-} = \big\{ -{\varepsilon }_7; -{\varepsilon }_6 \pm {\varepsilon }_i, \, 1 \le i \le 5; \, 
\onehalf \big( \sum^5_{i=1} (-1)^{k(i)}{\varepsilon }_i - {\varepsilon }_6- {\varepsilon }_7 \big) , \\ 
\rule{0pt}{12pt} \hspace{.5in} \mbox{where $k(i)=0$ or $1$ and $\sum^5_{i=1} k(i) \in \{ 1,3,5 \} $}  \big\} 
\end{array} $ & \\
& nilradical roots: $\widehat{\mathcal{R}} = $ & \\
& $\left\{ \begin{array}{rl} 
{\widehat{\alpha }}^{\, 0} & = -{\varepsilon }_7| \mathfrak{h} \\
{\widehat{\alpha }}^{\, j}& = ({\varepsilon }_j - {\varepsilon }_6) |\mathfrak{h}, \, 1 \le j \le 5 \\
{\widehat{\alpha }}^{\, j+5} & =  (-{\varepsilon }_j - {\varepsilon }_6) |\mathfrak{h}, \, 1 \le j \le 5 \\
{\widehat{\alpha }}_{-} & = -\onehalf \big( \sum^7_{i=1}{\varepsilon }_i \big)|\mathfrak{h}  \\
\rule{0pt}{12pt} \hspace{-4pt} {\widehat{\alpha }}_j & = 
\onehalf \big( \sum^5_{i=1} (-1)^{k(i)}{\varepsilon }_i - {\varepsilon }_6- {\varepsilon }_7 \big) |\mathfrak{h} , \\ 
&\hspace{.25in} \parbox[t]{3in}{where $k(j) = 1$; $k(i) =i$, if $i \notin \{ j \} $ \\
and $j \in \{ 1, \ldots ,5 \}$}  \\
\rule{0pt}{12pt} \hspace{-4pt} {\widehat{\alpha }}_{j\ell m} & = \onehalf \big( 
\sum^5_{i=1} (-1)^{k(i)}{\varepsilon }_i - {\varepsilon }_6- {\varepsilon }_7 \big) |\mathfrak{h} , \\
&\hspace{.25in} \parbox[t]{3in}{where $k(j) =k(\ell ) = k(m) = 1$; $k(i) =i$, if \\ $i \notin \{ j , \ell , m \} $ 
and $j, \ell , m \in \{ 1, \ldots ,5 \}$ distinct}
\end{array} \right. $ & \\
& additive relations in $\widehat{\mathcal{R}}$: none & \\
& nilradical structure: $Z_{16}$
\end{tabular}
\begin{tabular}{lll}
\vspace{.25in} & & \\
8.3 ${\widetilde{\mathbf{E}}}^3_7$ & {\large ${\mathbf{E}}^{(2)}_7$}: \hspace{-.65in} \begin{tabular}{l} 
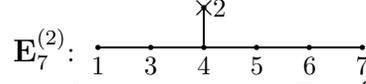
\begin{picture}(0,15)(-40,-3)
\begin{drawjoin}
\jput(0,5){\circle*{2}}
\jput(20,5){\circle*{2}} 
\jput(40,5){\circle*{2}}
\jput(60,5){\circle*{2}}
\jput(80,5){\circle*{2}}
\jput(100,5){\circle*{2}}
\end{drawjoin}
\put(40,5){\line(0,1){15}}
\put(40,20){\circle*{2}}
\put(40,20){\makebox(0,0){\large $\times$}} 
\put(0,-2){\makebox(0,0){\small $1$}}
\put(20,-2){\makebox(0,0){\small $3$}} 
\put(40,-2){\makebox(0,0){\small $4$}}
\put(46,20){\makebox(0,0){\small $2$}}
\put(60,-2){\makebox(0,0){\small $5$}}
\put(80,-2){\makebox(0,0){\small $6$}}
\put(100,-2){\makebox(0,0){\small $7$}}
\end{picture}%
\end{tabular}  & \\   
& simple positive roots: $\Pi = \big\{ {\alpha }_i, \, i \in I =\{ 1,3,4,5,6 \} \big\} $ & \\
& Cartan subalgebra: $\mathfrak{h} = {\spann }_{\C} \{ h_i, \, i \in I  \} $ & \\
& \parbox[t]{4,5in}{Cartan matrix: 
$C^{(2,2)}_{{\mathbf{E}}_7} = {\big( {\alpha }_i(h_j) \big) }_{i,j \in I }$ for 
root system $\mathcal{R} = {\mathbf{A}}_6$  \\ 
\rule{.25in}{0in} generated by $\Pi $.} & \\ 
& Set & \\
& $\begin{array}{rl}
h^{\ast } & = 4h_1+7h_2+8h_3+12h_4+9h_5+6h_6+3h_7 \\
& = -e_1 +\sum^6_{i=2}e_i+4e_7  \end{array} $. & \\
& Then $\mathcal{R} = {\mathrm{R}}^0 $ & \\
& \parbox[t]{4.5in}{${\mathrm{R}}^{-} = \big\{ -{\varepsilon }_7; {\varepsilon }_1 -{\varepsilon }_i, \, 2 \le i \le 6 ;\,
-( {\varepsilon }_i+{\varepsilon }_j), \, 2\le i < j \le 6; \\ 
\rule{0pt}{12pt} \hspace{.5in} \onehalf \big( \sum^6_{i=2} (-1)^{k(i)}{\varepsilon }_i 
-{\varepsilon }_1- {\varepsilon }_7 \big) , \, \mbox{where $k(i)=0$ or $1$ and} \\
\rule{0pt}{12pt} \hspace{.5in} \rule{.25in}{0in} \sum^6_{i=2}k(i) \in \{ 3,5 \} ;  \\
\rule{0pt}{12pt} \hspace{.5in} \onehalf \big( \sum^6_{i=2} (-1)^{k(i)}{\varepsilon }_i +
{\varepsilon }_1- {\varepsilon }_7\big) , 
\mbox{where $k(i)=0$ or $1$ and} \\
\rule{0pt}{12pt} \rule{.25in}{0in} \sum^6_{i=2}k(i) \in \{ 2, 4 \}  \big\} $ }& \\
& nilradical roots: $\widehat{\mathcal{R}} = $ & \\
& $\left\{ \begin{array}{rl} 
{\widehat{\zeta }}_1 & = -{\varepsilon }_7  |\mathfrak{h} \\
{\widehat{\zeta}}_j & = \onehalf \big( \sum^6_{i=2}(-1)^{k(i)}{\varepsilon }_i +{\varepsilon }_1 - 
{\varepsilon }_7 \big) | \mathfrak{h},  \\ 
& \hspace{.25in} \parbox[t]{5in}{where $k(j) = 1$,  $k(i) =0$, if $i \notin \{ j  \} $   
and \\ $j \in \{ 2, \ldots ,6 \}$ distinct}  \\
\rule{0pt}{12pt} \hspace{-4pt} {\widehat{\zeta }}_{-} & = 
-\onehalf \big( \sum^7_{i=1}{\varepsilon }_i \big) |\mathfrak{h} \\
{\widehat{\alpha }}_j & = -({\varepsilon }_j - {\varepsilon }_1)|\mathfrak{h}, \, 2 \le j \le 6 \\
{\widehat{\alpha }}_{i<j} & = -({\varepsilon }_i+{\varepsilon }_j)|\mathfrak{h}, \, 2 \le i < j \le 6 \\
\rule{0pt}{12pt} \hspace{-4pt} {\widehat{\alpha }}_{ j\ell } & = 
 \onehalf ( \sum^6_{i=2} (-1)^{k(i)}{\varepsilon }_i -{\varepsilon }_1- {\varepsilon }_7)  |\mathfrak{h}, \\
&\rule{0pt}{12pt}  \hspace{.25in} \parbox[t]{5in}{where $k(j) = k(\ell ) = 0$,  $k(i) =1$, if $i \notin \{ j, \ell  \} $ \\  
and $j, \ell \in \{ 2, \ldots ,6 \}$ distinct}  \\
\rule{0pt}{12pt} {\widehat{\alpha }}^{\, j\ell} & = 
\onehalf (\sum^6_{i=2} (-1)^{k(i)}{\varepsilon }_i +{\varepsilon }_1 - {\varepsilon }_7 ) |\mathfrak{h},  \\
&\rule{0pt}{12pt}  \hspace{.25in} \parbox[t]{5in}{where $k(j) = k(\ell ) = 1$;  
$k(i) =0$, if $i \notin \{ j, \ell \} $ and \\ $j, \ell   \in \{ 2, \ldots ,6 \}$ distinct}    
\end{array} \right. $ & \\
& additive relations in $\widehat{\mathcal{R}}$: & \\
& $\left\{ \begin{array}{rl} 
{\widehat{\alpha }}_{j\ell } + {\widehat{\alpha }}^{\, j\ell } & = {\widehat{\zeta }}_1, \, \, \mbox{ $ j, \ell \in 
\{ 2, \ldots , 6 \} $ distinct} \\
{\widehat{\alpha }}_{i < j } + {\widehat{\alpha }}_{i j} & = {\widehat{\zeta }}_{-}, \, \,  2 \le i < j \le 6  \\
{\widehat{\alpha }}_j + {\widehat{\alpha }}_{j\ell } & = {\widehat{\zeta }}_j, \, \, \mbox{ $ j, \ell \in 
\{ 2, \ldots , 6 \} $ distinct} 
\end{array} \right. $ & \\
& nilradical structure: ${\mathrm{h}}^{{\widehat{\zeta }}_1}_{43} + {\mathrm{h}}^{{\widehat{\zeta }}_{-}}_{11} 
+ \sum^6_{k=2} +{\mathrm{h}}^{{\widehat{\zeta }}_k}_{11}$ & 
\end{tabular}
\begin{tabular}{lll} 
\vspace{.25in} & & \\
9. & {\large ${\mathbf{E}}_8$}: \hspace{-.65in} \begin{tabular}{l} 
\begin{picture}(0,15)(-40,-3)
\begin{drawjoin}
\jput(0,5){\circle*{2}}
\jput(20,5){\circle*{2}} 
\jput(40,5){\circle*{2}}
\jput(60,5){\circle*{2}}
\jput(80,5){\circle*{2}}
\jput(100,5){\circle*{2}}
\jput(120,5){\circle*{2}}
\end{drawjoin}
\put(40,5){\line(0,1){15}}
\put(40,20){\circle*{2}}
\put(0,-1){\makebox(0,0){\small $1$}}
\put(20,-1){\makebox(0,0){\small $3$}} 
\put(40,-1){\makebox(0,0){\small $4$}}
\put(44,20){\makebox(0,0){\small $2$}}
\put(60,-1){\makebox(0,0){\small $5$}}
\put(80,-1){\makebox(0,0){\small $6$}}
\put(100,-1){\makebox(0,0){\small $7$}}
\put(120,-1){\makebox(0,0){\small $8$}}
\end{picture}
\end{tabular} & \\ 
& positive simple roots: $\underline{\Pi} = \{
{\alpha }_1 = \onehalf ( {\varepsilon }_1 -\sum^7_{i=2}{\varepsilon }_i +{\varepsilon }_8)$ & \\
& \hspace{.25in} ${\alpha }_2 = {\varepsilon }_2+{\varepsilon }_1, \,
{\alpha }_i ={\varepsilon }_{i-1}-{\varepsilon }_{i-2}, \, 3 \le i \le 8 \} $ & \\
& Cartan subalgebra: $\underline{\mathfrak{h}} = {\spann }_{\C } \big\{ h_1 = 
\onehalf ( e_1 -\sum^7_{i=2}e_i +e_8) $, & \\ 
&\hspace{.25in} $h_2 = e_2+e_1, \, h_i = e_{i-1}-e_{i-2}, \, 3 \le i \le 8 \big\} $ & \\ 
\rule{0pt}{30pt} & \hspace{.25in} $C_{{\mathbf{E}}_8} = 
{\big( {\alpha }_i(h_j) \big) }_{\stackrel{\! \! 1 \le i \le 8}{\mbox{$\scriptscriptstyle 1 \le j \le 8$ } }} = 
\mbox{\tiny $\left( \begin{array}{rrrrrrrr} 
2 & 0 & -1 & 0 & 0 & 0 & 0 & 0   \\
0 & 2 & 0 & -1 & 0 & 0 & 0 & 0\\
-1 & 0 & 2 & -1 & 0 & 0 & 0 & 0 \\ 
0 & -1 & -1 & 2 & -1 & 0 & 0 & 0 \\
0 & 0 & 0 & -1 & 2 & -1 & 0 & 0 \\
0 & 0 & 0 & 0 & -1 & 2 & -1 & 0 \\ 
0 & 0 & 0 & 0 & 0 & -1 & 2 & -1 \\
0 & 0 & 0 & 0 & 0 & 0 & -1 & 2 
\end{array} \right) $} $ & \\ 
&\parbox[t]{4.5in}{Cartan matrix for root system $\underline{\mathcal{R}}$ generated by 
$\underline{\Pi }$.} & \\ 
& positive roots: ${\underline{\mathcal{R}}}^{+} =
\big\{ {\varepsilon }_i\pm{\varepsilon }_j, \, 1 \le i < j \le 7; \, \onehalf \sum^8_{i=1}(-1)^{k(i)} {\varepsilon }_i $, & \\
&\rule{.25in}{0in} where $k(i) =0$ or $1$ and $\sum^8_{i=1}k(i) \in \{0,2,4,6,8 \} \big\} $. & \\
&\rule{0pt}{16pt}  & \\
& & \\ 
9.1 ${\widetilde{\mathbf{E}}}^1_8$ & {\large ${\mathbf{E}}^{(1)}_8$}: \hspace{-.65in} \begin{tabular}{l} 
\begin{picture}(0,15)(-40,-3)
\begin{drawjoin}
\jput(0,5){\circle*{2}}
\jput(20,5){\circle*{2}} 
\jput(40,5){\circle*{2}}
\jput(60,5){\circle*{2}}
\jput(80,5){\circle*{2}}
\jput(100,5){\circle*{2}}
\jput(120,5){\circle*{2}}
\end{drawjoin}
\put(40,5){\line(0,1){15}}
\put(40,20){\circle*{2}}
\put(0,5){\makebox(0,0){\large $\times$}} 
\put(0,-1){\makebox(0,0){\small $1$}}
\put(20,-1){\makebox(0,0){\small $3$}} 
\put(40,-1){\makebox(0,0){\small $4$}}
\put(44,20){\makebox(0,0){\small $2$}}
\put(60,-1){\makebox(0,0){\small $5$}}
\put(80,-1){\makebox(0,0){\small $6$}}
\put(100,-1){\makebox(0,0){\small $7$}}
\put(120,-1){\makebox(0,0){\small $8$}}
\end{picture}  
\end{tabular} & \\
& positive simple roots: $\Pi  = \{ {\alpha }_i, \, 2 \le i \le 8 \} $ & \\
& Cartan subalgebra: $\mathfrak{h} = \{ h_i, \, 2 \le i \le 8 \} $ & \\
& Cartan matrix $C^{(1,1)}_{{\mathbf{E}}_8}$ for root system $\mathcal{R} = {\mathbf{D}}_7$ generated by 
$\Pi $. & \\
& Set $h^{\ast } = 4h_1+5h_2+7h_3+10h_4+8h_5+6h_6+4h_7 + 2h_8 = e_8$ & \\
& Then $\mathcal{R} = {\mathrm{R}}^0$ & \\
& ${\mathrm{R}}^{-} = \big\{  \onehalf \big( \sum^7_{i=1}(-1)^{k(i)}{\varepsilon}_i - {\varepsilon }_8 \big) $ 
where $k(i) =0$ or $1$ and & \\
&\rule{0pt}{13pt} \rule{.5in}{0in} $\sum^7_{i=1}k(i) \in \{ 1,3,5,7 \} \big\} $ & \\
& nilradical roots: $\widehat{\mathcal{R}} =$ & \\
& $ \left\{ \begin{array}{rl}
{\widehat{\alpha }}_{-} & = -\onehalf \big( \sum^8_{i=1}{\varepsilon }_i \big) |\mathfrak{h}; \\
\rule{0pt}{12pt} {\widehat{\alpha }}_j & = 
\onehalf \big( \sum^7_{i=1}(-1)^{k(i)}{\varepsilon}_i - {\varepsilon }_8 \big) \mathfrak{h},  \\
&\hspace{.25in}\parbox{3in}{where $k(j)=1$; $k(i) = 0$ for $i \notin \{ j \}$ and $j \in \{ 1,2,\ldots , 7 \} $}   \\
\rule{0pt}{12pt} {\widehat{\alpha }}_{j\ell m} & = 
\onehalf \big( \sum^7_{i=1}(-1)^{k(i)}{\varepsilon}_i - {\varepsilon }_8 \big) \mathfrak{h},  \\
&\hspace{.25in}\parbox{3in}{where $k(j)=k(\ell ) = k(m)=1$; $k(i) = 0$ for $i \notin \{ j,\ell ,m \}$ and 
$j,\ell ,m \in \{ 1,2,\ldots , 7 \} $ distinct } \\
\rule{0pt}{12pt} {\widehat{\alpha }}^{\, j \ell }& = 
\onehalf \big( \sum^7_{i=1}(-1)^{k(i)}{\varepsilon}_i - {\varepsilon }_8 \big) \mathfrak{h},  \\
&\hspace{.25in}\parbox{3in}{where $k(j)=k(\ell ) =0$; $k(i) = 1$ for $i \notin \{ j, \ell \}$ and 
$j, \ell  \in \{ 1,2,\ldots , 7 \} $ distinct }   
\end{array} \right. $ &  \\
& additive relations in $\widehat{\mathcal{R}}$: none & \\
& nilradical structure: $Z_{78}$
\end{tabular}
\begin{tabular}{lll} 
\vspace{-.5in} & & \\
9.2 ${\widetilde{\mathbf{E}}}^2_8$ & {\large ${\mathbf{E}}^{(8)}_8$}: \hspace{-.65in} \begin{tabular}{l} 
\begin{picture}(0,15)(-40,-3)
\begin{drawjoin}
\jput(0,5){\circle*{2}}
\jput(20,5){\circle*{2}} 
\jput(40,5){\circle*{2}}
\jput(60,5){\circle*{2}}
\jput(80,5){\circle*{2}}
\jput(100,5){\circle*{2}}
\jput(120,5){\circle*{2}}
\end{drawjoin}
\put(40,5){\line(0,1){15}}
\put(40,20){\circle*{2}}
\put(120,5){\makebox(0,0){\large $\times$}} 
\put(0,-1){\makebox(0,0){\small $1$}}
\put(20,-1){\makebox(0,0){\small $3$}} 
\put(40,-1){\makebox(0,0){\small $4$}}
\put(44,20){\makebox(0,0){\small $2$}}
\put(60,-1){\makebox(0,0){\small $5$}}
\put(80,-1){\makebox(0,0){\small $6$}}
\put(100,-1){\makebox(0,0){\small $7$}}
\put(120,-1){\makebox(0,0){\small $8$}}
\end{picture}  
\end{tabular} & \\ 
& positive simple roots: $\Pi  = \{ {\alpha }_i, \, 1 \le i \le 7 \} $ & \\
& Cartan subalgebra: $\mathfrak{h} = \{ h_i, \, 1 \le i \le 7 \} $ & \\
& Cartan matrix $C^{(8,8)}_{{\mathbf{E}}_8}$ for root system $\mathcal{R} = {\mathbf{E}}_7$ generated by 
$\Pi $. & \\
& Set $h^{\ast } = 2h_1+3h_2+4h_3+6h_4+5h_5+4h_6+3h_7 + 2h_8 = e_7 + e_8$ & \\
& Then $\mathcal{R} = {\mathrm{R}}^0$ & \\
& ${\mathrm{R}}^{-} = \big\{ \pm {\varepsilon }_i -{\varepsilon }_7, \, 1 \le i \le 6; \,  \onehalf \big( \sum^6_{i=1}(-1)^{k(i)}{\varepsilon }_i - ({\varepsilon }_7 + {\varepsilon }_8) \big) $, & \\ 
& \hspace{.5in} where $k(i) = 0$ or $1$ and $\sum^6_{i=1}k(i) \in \{ 0,2,4,6 \}  \big\} $ & \\ 
& nilradical roots: $\widehat{\mathcal{R}} = $ & \\
& $ \left\{ \begin{array}{rl}
{\widehat{\alpha}}_j & = ( {\varepsilon }_j - {\varepsilon }_7)|\mathfrak{h}, \, 1 \le j \le 6 \\
{\widehat{\alpha }}_{6+j} & = ( -{\varepsilon }_j - {\varepsilon }_7)|\mathfrak{h}, \, 1 \le j \le 6 \\
{\widehat{\alpha }}_{+} & =
\big( \onehalf \sum^6_{i=1}{\varepsilon }_i -({\varepsilon}_7 +{\varepsilon }_8) \big) |\mathfrak{h}, \\
\rule{0pt}{12pt} {\widehat{\alpha }}_{-} & = -\onehalf \big( \sum^8_{i=1}{\varepsilon }_i \big) |\mathfrak{h} \\ 
\rule{0pt}{12pt} {\widehat{\alpha }}_{j\ell} & = \big( \onehalf \sum^6_{i=1}(-1)^{k(i)}{\varepsilon }_i -
({\varepsilon }_7 +{\varepsilon }_8) \big) |\mathfrak{h}, \\ 
&\hspace{.5in}\parbox{4.5in}{where $k(j) =k(\ell ) = 1$; $k(i) =0$ \\ 
if $i \notin \{ j, \ell \}$, $j, \ell \in \{ 1, \ldots , 6 \}$ distinct } \\
\rule{0pt}{12pt} {\widehat{\alpha }}^{j\ell } & = \big( \onehalf \sum^6_{i=1}(-1)^{k(i)}{\varepsilon }_i - 
({\varepsilon }_7 +{\varepsilon }_8) \big) |\mathfrak{h}, \\ 
\rule{0pt}{16pt}&\hspace{.25in}\parbox{4.5in}{where $k(j) = k(\ell ) =0$; $k(i) =1$ if \\
$i \notin \{ j, \ell \}$, $j, \ell  \in \{ 1, \ldots , 6 \}$ distinct } 
\end{array} \right. $ & \\
\rule{0pt}{12pt} & additive relations in $\widehat{\mathcal{R}}$: none & \\
& nilradical structure: $Z_{32}$
\end{tabular}

\section*{Acknowledgments}
I would like to thank Prof. Nick Burgoyne of the University of California at Santa Cruz for checking that 
the nilradical of a very special sandwich algebra is indeed a sandwich and for other very useful remarks. 
I also would like to thank Prof. Arjeh M. Cohen of the Technical University of Eindhoven for his comments 
on an early version of this paper.

\end{document}